%% file: HKM-update.tex
\documentclass[11pt,reqno]{amsart}

\input{macros}

\title{An update on Heisenberg and Kac-Moody categorification}

\author{Jonathan Brundan}
\address[J.B.]{Department of Mathematics, University of Oregon, Eugene, OR, USA}
\urladdr{\href{https://pages.uoregon.edu/brundan}{https://pages.uoregon.edu/brundan}, \textrm{\textit{ORCiD}:} \href{https://orcid.org/0009-0009-2793-216X}{0009-0009-2793-216X}}
\email{brundan@uoregon.edu}

\author{Alistair Savage}
\address[A.S.]{Department of Mathematics and Statistics, University of Ottawa, Ottawa, ON, Canada}
\urladdr{\href{https://alistairsavage.ca}{alistairsavage.ca}, \textrm{\textit{ORCiD}:} \href{https://orcid.org/0000-0002-2859-0239}{orcid.org/0000-0002-2859-0239}}
\email{alistair.savage@uottawa.ca}

\author{Ben Webster}
\address[B.W.]{Department of Pure Mathematics, University of Waterloo \& Perimeter Institute for Theoretical Physics,
Waterloo, ON, Canada}
\urladdr{\href{https://bwebste.github.io}{https://uwaterloo.ca/scholar/b2webste},\textrm{\textit{ORCiD}:} \href{https://orcid.org/0000-0003-1896-5540}{0000-0003-1896-5540}}
\email{ben.webster@uwaterloo.ca}

\thanks{Research supported in part by NSF grant DMS-2348840 (J.B.), NSERC grant 2023-03842 (A.S.) and NSERC grant RGPIN-2024-03760 (B.W.)}



\begin{document}

\begin{abstract}
Heisenberg categories act on many Abelian categories appearing in type A representation theory. There is also a general procedure to construct from a Heisenberg action another action of a Kac-Moody 2-category for some associated Cartan matrix. One of the adjunctions on the Kac-Moody side is matched up in an easy way with adjunctions on the Heisenberg side, but the second adjunction is much harder to describe. In this paper, we derive explicit formulae for this difficult adjunction, 
leading to some further simplifications to the existing theory.
\end{abstract}

%\dedicatory{Dedicated to Dan Nakano on the occasion of his sixtieth birthday}

\maketitle

%\vspace{-2mm}

\input{s1-intro}
\input{s2-KM}
\input{s3-H}

\input{s4-HtoKM}
\input{s5-example}

\input{s6-example}
\input{s7-bubbles}
\input{s8-quantum}

\bibliographystyle{alphaurl}
\bibliography{HKM-update}

\end{document}

%% file: macros.tex
\usepackage{
    amsmath,
    amsfonts,
    amssymb,
    amsthm,
    amscd,
    comment,
    enumitem,
    etoolbox,
    textcomp,
    gensymb,
    mathtools,
    mathdots,
    stmaryrd,
    accents
}
\usepackage[dvipsnames]{xcolor}
\usepackage[all]{xy}

\makeatletter
\@namedef{subjclassname@2020}{
\textup{2020} Mathematics Subject Classification}
\makeatother

\usepackage[T1]{fontenc}
\usepackage{bbm}
\usepackage[colorlinks=true, linkcolor=black, citecolor=black, urlcolor=black, breaklinks=true]{hyperref}

\usepackage[nameinlink]{cleveref}
\crefname{table}{Table}{Tables}
\crefname{cor}{Corollary}{Corollaries}
\crefname{defin}{Definition}{Definitions}
\crefname{conv}{Convention}{Conventions}
\crefname{eg}{Example}{Examples}
\crefname{rem}{Remark}{Remarks}
\crefname{lem}{Lemma}{Lemmas}
\crefname{prop}{Proposition}{Propositions}
\crefname{theo}{Theorem}{Theorems}
\crefname{conj}{Conjecture}{Conjectures}
\crefname{equation}{}{}
\crefname{enumi}{}{}
\crefname{table}{Table}{Tables}
\crefname{subsection}{Subsection}{Subsections}
\crefname{section}{Section}{Sections}

\setenumerate[1]{label=(\arabic*)}          % Only need to reference enumerate items with \ref{}

\newcommand{\NH}{\mathrm{NH}}
\newcommand{\U}{\mathrm{U}}

\newcommand{\M}{\mathrm{M}}
\newcommand{\UU}{\mathfrak{U}}

\newcommand\C{\mathbb{C}}
\newcommand\N{\mathbb{N}}

\newcommand\Q{\mathbb{Q}}

\newcommand\Z{\mathbb{Z}}
\newcommand\kk{\Bbbk}

\newcommand\one{\mathbbm{1}}

\newcommand\catB{\mathbf{B}}
\newcommand\catC{\mathbf{C}}

\newcommand\qHeis{q\text{-}\hspace{-1pt}\mathbf{Heis}}
\newcommand\Heis{\mathbf{Heis}}
\newcommand\heis{\mathrm{H}}
\newcommand\catR{\mathbf{R}}

\newcommand\B{\mathcal{B}}

\newcommand\eps{\varepsilon}

\newcommand{\Dist}{\operatorname{Dist}}
\newcommand{\End}{\operatorname{End}}

\newcommand\Rep{\mathbf{Rep}}

\newcommand\ch{\operatorname{char}}
\newcommand\transpose{\textup{T}}
              
\newcommand\lround{(\!(}
\newcommand\rround{)\!)}
\def\mod#1{{#1\text{-}\mathbf{mod}_{\operatorname{fd}}}}
\DeclareMathOperator{\ob}{ob}
\DeclareMathOperator{\id}{id}

\DeclareMathOperator{\wt}{wt}

\usepackage{tikz}
\usetikzlibrary{arrows,patterns}
\usetikzlibrary{arrows.meta}
\usetikzlibrary{shapes.geometric}
\usetikzlibrary{decorations.markings,decorations.pathreplacing,decorations.pathmorphing}
\usetikzlibrary{calc}
\usepackage{tikz-cd}

\newcommand*\getscale[1]{%
  \begingroup
    \pgfgettransformentries{\scaleA}{\scaleB}{\scaleC}{\scaleD}{\whatevs}{\whatevs}%
    \pgfmathsetmacro{#1}{sqrt(abs(\scaleA*\scaleD-\scaleB*\scaleC))}%
    \expandafter
  \endgroup
  \expandafter\def\expandafter#1\expandafter{#1}%
}

\tikzset{anchorbase/.style={>=To,baseline={([yshift=-0.5ex]current bounding box.center)}}}
\tikzset{ % Syntax: \begin{tikzpicture}[centerzero={0,0.2}]
    centerzero/.style={>=To,baseline={([yshift=-0.5ex]#1)}},
    centerzero/.default={0,0}
}
\tikzset{wipe/.style={white,line width=4pt}}

\tikzset{Q/.style={semithick,violet}}
\tikzset{H/.style={semithick,purple}}
\tikzset{KM/.style={thin,black}}
\tikzset{notchy/.style={thin,black}}

\tikzset{gcolor/.style={green!60!black,semithick}}
\tikzset{pinhead/.style={gray,fill=yellow!40!white}}
\tikzset{fakebubble/.style={fill=red!20!white}}
\tikzset{Internalbubble/.style={fill=blue!10!white}}
\tikzset{internalbubble/.style={fill=green!20!white}}
\tikzset{bulb/.style={fill=yellow!40!white}}
\tikzset{shadow/.style={double,thin}}
\newcommand\braidup{to[out=up,in=down]}
\newcommand\braiddown{to[out=down,in=up]}

\newcommand\notch[2][0]{% \notch[angle]{position}, where angle is the angle to one of the endpoints
    \draw[notchy] (#2)++(#1:0.06) -- ++(#1:-0.12)
}
\newcommand\projcr[1]{
    \draw (#1) node[rotate=45] {$\color{black}\scriptstyle{\Box}$}
}

\newcommand{\circled}[2]{\node[circle,draw,bulb,inner sep=.85pt] at (#1) {$\hspace{0.06pt}\color{black}\scriptscriptstyle #2$}}

\newcommand{\squared}[2]{\node[regular polygon,regular polygon sides=6,draw,bulb,inner sep=.5pt] at (#1) {$\hspace{0.06pt}\color{black}\scriptscriptstyle #2$}}

\newcommand{\pin}[3]{
    \path (#1) node[inner sep=1.6pt] (x) {} to (#2) node[rectangle,rounded corners,draw,pinhead,inner sep=2.5pt](y) {$\color{black}\scriptstyle#3$};
    \draw[Triangle Cap-,thick,gray] (x)--(y);
    \opendot{#1}
}

\newcommand{\pinpin}[4]{
\path (#1) node[inner sep=1.6pt] (x) {} to (#3) node[rectangle,rounded corners,draw,pinhead,inner sep=2.5pt](y) {$\color{black}\scriptstyle#4$};
    \draw[Triangle Cap-,thick,gray] (x)--(#2)--(y);
\opendot{#2}; \opendot{#1}
}

\newcommand{\pinpinpin}[5]{
\path (#1) node[inner sep=1.6pt] (x) {} to (#4) node[rectangle,rounded corners,draw,pinhead,inner sep=2.5pt](y) {$\color{black}\scriptstyle#5$};
    \draw[Triangle Cap-,thick,gray] (x)--(#2)--(#3)--(y);
\opendot{#3}; \opendot{#2}; \opendot{#1}
}

\colorlet{darkg}{green!50!black}

\newcommand\dotlabel[1]{$\color{teal}\scriptstyle{#1}$}
\newcommand\regionlabel[1]{$\color{gray}\scriptstyle{#1}$}
\newcommand\strandlabel[1]{$\color{olive}\scriptstyle{#1}$}
\newcommand\stringlabel[2]{\node at (#1) {\strandlabel{#2}}}

\newcommand\botlabel[1]{node[inner sep=0.5pt,anchor=north] {\strandlabel{#1}}}
\newcommand\toplabel[1]{node[anchor=south,inner sep=0.5pt] {\strandlabel{#1}}}
 
\newcommand\bubblelabel[2]{\node at (#1) {$\scriptstyle{#2}$}} 
\newcommand\region[2]{\node at (#1) {\regionlabel{#2}}} 
\newcommand\strand[2]{\node at (#1) {\strandlabel{#2}}} 

\newcommand\opendot[1]{
    \node at (#1) {\color{white} $\scriptstyle{\bullet}$};\node at (#1) {$\scriptstyle{\circ}$}
}

\newcommand\multopendot[3]{
    \opendot{#1};
    \draw (#1) node[anchor=#2] {\dotlabel{#3}}
}

\newcommand\filledclockwisebubble[1]{%
\draw[-to,fakebubble] (#1)++(0,.2) arc(90:-270:0.2);
}

\newcommand\stretchedanticlockwisebubble[1]{%
\draw[-to,fakebubble] ($(#1)+(0,.21)$) to ($(#1)+(-.11,.21)$) to [out=180,in=180,looseness=2] ($(#1)+(-.11,-.21)$)
to ($(#1)+(.11,-.21)$) to [out=0,in=0,looseness=2] ($(#1)+(.11,.21)$) to ($(#1)+(0,.21)$);
}

\newcommand\stretchedclockwisebubble[1]{%
\draw[-to,fakebubble] ($(#1)+(0,.21)$) to ($(#1)+(.11,.21)$) to [out=0,in=0,looseness=2] ($(#1)+(.11,-.21)$)
to ($(#1)+(-.11,-.21)$) to [out=180,in=180,looseness=2] ($(#1)+(-.11,.21)$) to ($(#1)+(0,.21)$);
}

\newcommand\filledanticlockwisebubble[1]{%
  \draw[-to,fakebubble] (#1)++(0,.2) arc(90:450:0.2);
}

\newcommand\anticlockwiseinternalbubbleL[2][]{%
 \getscale{\scalefactor};
 \draw[-to,internalbubble] ($(#2)+({-.18 / \scalefactor},0)$)  arc(-180:180:{.18/\scalefactor});
 \bubblelabel{#2}{#1}
}

\newcommand\anticlockwiseinternalbubbleNW[2][]{%
 \getscale{\scalefactor};
 \draw[-to,internalbubble] ($(#2)+({-.12 / \scalefactor},0)$)  arc(-225:135:{.12/\scalefactor});
 \bubblelabel{#2}{#1}
}

\newcommand\clockwiseinternalbubbleR[2][]{%
 \getscale{\scalefactor};
 \draw[-to,internalbubble] ($(#2)+({.18 / \scalefactor},0)$) arc(0:-360:{.18/\scalefactor});
  \bubblelabel{#2}{#1}
}

\newcommand\clockwiseinternalbubbleSE[2][]{%
 \getscale{\scalefactor};
 \draw[-to,internalbubble] ($(#2)+({.12 / \scalefactor},0)$) arc(-45:-425:{.12/\scalefactor});
  \bubblelabel{#2}{#1}
}

\newcommand\anticlockwiseinternalbubbleR[2][]{%
 \getscale{\scalefactor};
 \draw[-to,internalbubble] ($(#2)+({.18 / \scalefactor},0)$)  arc(-3600:0:{.18/\scalefactor});
 \bubblelabel{#2}{#1}
}

\newcommand\clockwiseinternalbubbleL[2][]{%
 \getscale{\scalefactor};
 \draw[-to,internalbubble] ($(#2)+({-.18 / \scalefactor},0)$) arc(180:-180:{.18/\scalefactor});
  \bubblelabel{#2}{#1}
}

\newcommand\anticlockwiseInternalbubbleL[2][]{%
 \getscale{\scalefactor};
 \draw[-to,Internalbubble] ($(#2)+({-.18 / \scalefactor},0)$)  arc(-180:180:{.18/\scalefactor});
 \bubblelabel{#2}{#1}
}

\newcommand\anticlockwiseInternalbubbleR[2][]{%
 \getscale{\scalefactor};
 \draw[-to,Internalbubble] ($(#2)+({.18 / \scalefactor},0)$)  arc(-3600:0:{.18/\scalefactor});
 \bubblelabel{#2}{#1}
}

\newcommand\clockwiseInternalbubbleL[2][]{%
 \getscale{\scalefactor};
 \draw[-to,Internalbubble] ($(#2)+({-.18 / \scalefactor},0)$) arc(180:-180:{.18/\scalefactor});
  \bubblelabel{#2}{#1}
}

\newtheorem{theo}{Theorem}[section]

\newtheorem{lem}[theo]{Lemma}
\newtheorem{cor}[theo]{Corollary}

\theoremstyle{definition}
\newtheorem{defin}[theo]{Definition}

\newtheorem{rem}[theo]{Remark}

\numberwithin{equation}{section}
\allowdisplaybreaks

%% file: s1-intro.tex
\setcounter{section}{0}

%=====================================
\section{Introduction}\label{s1-intro}
%=====================================

In 2004,
Chuang and Rouquier introduced the notion
of an {\em $\mathfrak{sl}_2$-categorification}; see \cite{CR}. A few years later, in \cite{Rou}, the definition was extended from $\mathfrak{sl}_2$ to other symmetrizable Kac-Moody algebras. 
The resulting axiomatic framework of {\em Kac-Moody categorifications} was motivated in part by classical examples of Abelian categories arising in the representation theory of symmetric and general linear groups (and related Hecke algebras and quantum groups).
In these examples, the combinatorics is controlled by 
an underlying Kac-Moody algebra, either
$\mathfrak{sl}_\infty$ when the (quantum) characteristic of the ground field is 0, or $\widehat{\mathfrak{sl}}_p$ in characteristic $p>0$. This point of view was taken already by 
Grojnowski \cite{Groj}, who made an important step in the development of the general theory
by explaining the coincidence
between the crystal graph for the basic representation of $\widehat{\mathfrak{sl}}_p$ determined in \cite{MM} and the modular branching graph for representations of symmetric groups from \cite{Klesh}. Kac-Moody categorifications have been used to construct many interesting Morita and derived equivalences, including the ones used
in \cite{CR} to prove Brou\'e's Abelian Defect Conjecture
for symmetric groups.

In \cite{LLT}, Lascoux, Leclerc and Thibon formulated
a precise categorification conjecture relating the representation theory of cyclotomic Hecke algebras of type A at a complex $p$th root of unity to canonical bases of integrable highest weight representations of
$\widehat{\mathfrak{sl}}_p$. This was quickly proved by Ariki \cite{Ariki}. Subsequently, an analogous conjecture relating graded representation theory of the {\em quiver Hecke algebras} from \cite{Rou,KL1} to canonical bases of  quantum groups of other 
symmetric Cartan types was formulated by Khovanov and Lauda, and proved by Varagnolo and Vasserot in \cite{VV}; see also \cite{Rou2}.
In fact, 
the Lascoux-Leclerc-Thibon conjecture can be seen as a special case of
the Khovanov-Lauda conjecture due to the isomorphism between cyclotomic Hecke algebras and cyclotomic quiver Hecke algebras in 
type A established in \cite{BKiso}. 

The general definition of a Kac-Moody categorification rests on the 2-categories introduced in \cite{Rou,KL3} known as {\em Kac-Moody 2-categories} or {\em 2-quantum groups}. Roughly speaking, these are obtained from quiver Hecke algebras by adjoining certain right duals subject to some non-trivial relations, the {\em inversion relations}, which are categorifications of the commutation relations between $e_i$ and $f_j$ in the underlying Kac-Moody algebra.
However, it is surprisingly difficult to fit the classical examples from type A representation theory into the general formalism.
One approach was developed by Rouquier in \cite[Th.~5.27]{Rou}, a result he called ``control by $K_0$''. This is often useful but still requires a lot of preparation. 

In \cite{HKM}, we developed another approach 
to the construction of Kac-Moody categorifications based on the notion of a {\em Heisenberg categorification}. This rests on another family of strict monoidal categories, the (degenerate or quantum) {\em Heisenberg categories} $\Heis_\kappa$, which are obtained from (degenerate or quantum) affine Hecke algebras of type A by adjoining a certain right dual subject to another inversion relation depending on the value of $\kappa \in \Z$, called the {\em central charge}.
This is a fruitful approach because it
is much easier to check that the Abelian categories arising in many sorts of type A representation theory are Heisenberg categorifications.
Then the main construction in \cite{HKM} gives a systematic procedure to pass from a Heisenberg categorification to a Kac-Moody categorification. 
After that, all of the powerful structural results about Kac-Moody categorifications can be applied.

The bridge from Heisenberg to Kac-Moody has as its starting point the
isomorphism between cyclotomic quotients of affine Hecke algebras and quiver Hecke algebras mentioned already. This is then upgraded 
to obtain from the given action of the Heisenberg category a much less obvious action of a Kac-Moody 2-category. On both sides, there are two families of adjunctions, one defined by rightward cups and caps, and the other defined by leftward cups and caps. The rightward cups and caps are matched up in an obvious way, but then one needs to check that the Kac-Moody inversion relations hold on the Heisenberg side.
In \cite{HKM}, this was done via an indirect argument involving the Chinese Remainder Theorem. The argument does not produce any explicit formula relating the leftward cups and caps on the two sides.

In this article, we begin 
in \cref{s2-KM,s3-H,s4-HtoKM} with a survey of the definitions of degenerate Heisenberg and Kac-Moody categorifications and the main ``Heisenberg to Kac-Moody'' construction from \cite{HKM}. In \cref{s5-example,s6-example}, we work through two classic examples illustrating the general theory, discussing the representation theory of the symmetric groups,
and representations of the general linear group scheme and its Frobenius kernels in positive characteristic. With all of this background in place, in \cref{s7-bubbles}, we formulate and prove the main new result of the article. This gives explicit formulae relating the second adjunction on the Kac-Moody side---the one defined by leftward cups and caps---to the second adjunction on the Heisenberg side. 
Finally, in \cref{s8-quantum}, we treat the quantum case.

The interesting works of Bao and Wang \cite{BW18KL} and Ehrig and Stroppel \cite{ES}
point to the existence of an equally rich theory of categorifications arising from classical groups of types B/C/D. The analogues of the degenerate and quantum Heisenberg categories are the affine Brauer category introduced in \cite{RS} and its quantum analogue, the Kauffman category. The analogues of the Kac-Moody 2-categories are the quasi-split 2-iquantum groups introduced in \cite{BWWquasisplit}; we expect the ones of diagonal Cartan type A and of quasi-split Satake type AIII (possibly affine) to play a role. However, the expected bridge from Brauer categorifications to iquantum group categorifications is still to be built. For 2-iquantum groups at this level of generality, there is no satisfactory analogue of the inversion relation, so a different strategy will be 
needed compared to \cite{HKM}.
This was one of the main motivations prompting us to think again in this article about the identification of leftward cups and caps in Heisenberg and Kac-Moody categorification. 

The two sorts of categorification discussed so far can be thought of as the GL and the OSp branches of a wider story.
There are two more classical families of supergroups over algebraically closed fields: the periplectic supergroups P and the isomeric supergroups Q. In the isomeric case, the first two authors are currently working on the development of a theory of {\em isomeric Heisenberg and Kac-Moody categorification}, building in particular on \cite{KKT}.
Significant progress in the periplectic case has been made recently by Nehme and Stroppel in \cite{NS}.

%\pagebreak

%% file: s2-KM.tex
\setcounter{section}{1}

%=====================================
\section{Kac-Moody categorifications}\label{s2-KM}
%=====================================

Suppose that we are given:
\begin{itemize}
\item A symmetrizable generalized Cartan matrix $A = (a_{i,j})_{i, j \in I}$.
This means that $a_{i,i} = 2$ for all $i \in I$,
$a_{i,j} \in -\N$ for $i\neq j$ in $I$,
$a_{i,j} = 0\Leftrightarrow a_{j,i}=0$, and
there are given positive integers $d_i\:(i \in I)$
such that $d_i a_{i,j} = d_j a_{j,i}$ 
for all $i,j \in I$.
We do not insist that the set $I$ is finite, but the number of non-zero entries in each row and column of $A$ should be finite.
\item
A free Abelian group $X$, 
the {\em weight lattice}, containing 
elements $\alpha_i\:(i \in I)$, the {\em simple roots}, and $\varpi_i\:(i \in I)$, the {\em fundamental weights}, and homomorphisms $h_i:X \rightarrow \Z\:(i \in I)$
such that $h_i(\alpha_j) = a_{i,j}$ and $h_i(\varpi_j) = \delta_{i,j}$ for all $i,j \in I$. We note that the fundamental weights are necessarily linearly independent, but the simple roots may be dependent.
\end{itemize}
To this data, we associate a locally unital $\C$-algebra $\dot\U$, which is the modified form of the universal enveloping algebra of a Kac-Moody algebra of this Cartan type. It has mutually orthogonal distinguished idempotents $1_\lambda\:(\lambda \in X)$, and it is generated
by elements
\begin{align}
e_i 1_\lambda &= 1_{\lambda + \alpha_i} e_i,
&
1_\lambda f_i &= f_i 1_{\lambda+\alpha_i}
\end{align}
for $i \in I, \lambda \in X$
subject to idempotented versions of the usual Serre relations.
Let $\dot\U_\Z$ be the $\Z$-form for $\dot\U$ generated by the divided powers $e_i^{(n)} 1_\lambda := e_i^n 1_\lambda / n!$
and $1_\lambda f_i^{(n)} := 1_\lambda f_i^n / n!$.

The main object of interest in this paper is a certain 2-category $\UU$ which is a categorification of the $q$-analogue of
$\dot\U_\Z$; see \cref{notquite}.
It was introduced originally by
Rouquier \cite{Rou} and, independently, Khovanov and Lauda \cite{KL3}; our exposition follows \cite[Sec.~2]{BWWquasisplit}.
Fix an algebraically closed ground field $\kk$ and parameters $Q_{i,j}(x,y) \in \kk[x,y]$ for all $i,j \in I$
such that 
$Q_{i,i}(x,y) = 0$,
$Q_{i,j}(x,y) = Q_{j,i}(y,x)$, and the following hold for $i \neq j$:
\begin{itemize}
\item
The polynomial $Q_{i,j}(x,y)$ is
homogeneous of degree $-2 d_i a_{i,j}$ if $x$ and $y$ are assigned the
degrees $2d_i$ and $2 d_j$, respectively. 
\item
The coefficient $Q_{i,j}(1,0)$ 
is non-zero.
\end{itemize}
We also let 
\begin{equation}
t_{i,j} := \begin{cases}
Q_{i,j}(1,0)&\text{if $i \neq j$}\\
1&\text{if $i=j$.}
\end{cases}
\end{equation}
The {\em Kac-Moody 2-category} $\UU$
with these parameters is the
$\kk$-linear 2-category 
with object set $X$, generating $1$-morphisms
\begin{align}
E_i \one_\lambda =\one_{\lambda+\alpha_i} E_i &\colon \lambda \to \lambda + \alpha_i,&
\one_{\lambda} F_i = F_i \one_{\lambda+\alpha_i} &\colon \lambda+\alpha_i \to \lambda
\end{align}
for $\lambda \in X$ and $i \in I$,
with identity 2-endomorphisms represented by the oriented strings
$\begin{tikzpicture}[KM,anchorbase]
\draw[-to] (0,-0.2)\botlabel{i} -- (0,0.2);
\region{-0.4,0}{\lambda+\alpha_i};
\region{0.2,0}{\lambda};
\end{tikzpicture}$ and
$\begin{tikzpicture}[KM,anchorbase]
\draw[to-] (0,-0.2) \botlabel{i}-- (0,0.2);
\region{-0.23,0}{\lambda};
\region{0.4,0}{\lambda+\alpha_i};
\end{tikzpicture}\ $,
and the four families of generating 2-morphisms
displayed in \cref{table1}.
\begin{table}
\begin{align*}
\begin{array}{|l|c|}
\hline
\hspace{13mm}\text{Generator}&\text{Degree}\\
\hline
\ \ \begin{tikzpicture}[KM,centerzero]
\draw[-to] (0,-0.3) \botlabel{i} -- (0,0.3);
\opendot{0,0};
\region{0.23,0}{\lambda};
\node at (0,.3) {$\phantom x$};
\end{tikzpicture}\ 
\colon E_i \one_\lambda \Rightarrow E_i \one_\lambda&2 d_i\\
\,\begin{tikzpicture}[KM,centerzero,scale=.9]
\draw[-to] (-0.3,-0.3) \botlabel{i} -- (0.3,0.3);
\draw[-to] (0.3,-0.3) \botlabel{j} -- (-0.3,0.3);
\region{0.38,0}{\lambda};
\node at (0,.3) {$\phantom x$};
\end{tikzpicture}
\;\colon E_i E_j \one_\lambda \Rightarrow E_j E_i \one_\lambda
&-d_i a_{i,j}\\
\:\begin{tikzpicture}[KM,centerzero]
\draw[-to] (-0.25,-0.15) \botlabel{i} to [out=90,in=90,looseness=3](0.25,-0.15);
\region{0.45,0.1}{\lambda};
\node at (0,.3) {$\phantom.$};
\node at (0,-.4) {$\phantom.$};
\end{tikzpicture}
\colon E_i F_i \one_\lambda \Rightarrow \one_\lambda
&d_i(1-h_i(\lambda))\\
\:\begin{tikzpicture}[KM,centerzero]
\draw[-to] (-0.25,0.15) \toplabel{i} to[out=-90,in=-90,looseness=3] (0.25,0.15);
\region{0.45,-.1}{\lambda};
\node at (0,.2) {$\phantom.$};\node at (0,-.3) {$\phantom.$};
\end{tikzpicture}
\colon \one_\lambda \Rightarrow F_i E_i \one_\lambda
&d_i(1+h_i(\lambda))
\\
\hline
\end{array}
\end{align*}
\caption{Generating 2-morphisms of the Kac-Moody 2-category}\label{table1}
\end{table}
The generating 2-morphisms are subject to relations still to be explained. We write these down using the following conventions:
\begin{itemize}
\item
We will usually only label one of the 2-cells in a string diagram with a weight---the others can then be worked out implicitly. When we omit
{\em all} labels in 2-cells, it should be understood that we are discussing something that holds
for all possible choices of labels.
\item
We will label strings just at one end. If we omit a label or orientation, it means that we are discussing something that holds for all possibilities.
\item
When a dot is labelled by a multiplicity,
we mean to take its power under vertical composition.
For a polynomial
$f(x) = \sum_{r=0}^n c_r x^r$, we use the
shorthand
\begin{equation*}\label{singlepin}
    \begin{tikzpicture}[KM,centerzero]
        \draw[-] (0,-0.3) -- (0,0.3);
             \pin{0,0}{-.7,0}{f(x)};
    \end{tikzpicture}
    = \begin{tikzpicture}[KM,centerzero]
        \draw[-] (0,-0.3) -- (0,0.3);
             \pin{0,0}{.7,0}{f(x)};
    \end{tikzpicture}\
    :=
    \sum_{r=0}^n c_r\
    \begin{tikzpicture}[KM,centerzero]
        \draw[-] (0,-0.3) -- (0,0.3);
        \multopendot{0,0}{west}{r};
    \end{tikzpicture}
\end{equation*}
to ``pin'' $f(x)$ to a dot on a string.
Similarly, for
$f(x,y) = \sum_{r=0}^n\sum_{s=0}^m c_{r,s} x^r y^s$, we use
\begin{align*}
\begin{tikzpicture}[KM,centerzero]
\draw[-] (0,-0.3) -- (0,0.3);
\draw[-] (0.4,-0.3) -- (0.4,0.3);
\pinpin{.4,0}{0,0}{-.8,0}{f(x,y)};
\end{tikzpicture}
=
\begin{tikzpicture}[KM,centerzero]
\draw[-] (0,-0.3) -- (0,0.3);
\draw[-] (0.4,-0.3) -- (0.4,0.3);
\pinpin{0,0}{0.4,0}{1.2,0}{f(x,y)};
\end{tikzpicture}\ 
& :=
\sum_{r=0}^n\sum_{s=0}^m c_{r,s}
\begin{tikzpicture}[KM,centerzero]
\draw[-] (0,-0.3) -- (0,0.3);
\draw[-] (0.4,-0.3) -- (0.4,0.3);
\multopendot{.4,0}{west}{s};
\multopendot{0,0}{east}{r};
\end{tikzpicture}
\end{align*}
This notation extends to polynomials $f(x,y,z)$ in three variables pinned to three dots, with
the convention that the variables in alphabetic order correspond to the dots ordered by the lexicographic order on their Cartesian coordinates.
Thus, $x$ corresponds to the leftmost dot, and the lowest one if there are
several such dots in the same vertical line.
\item
We use the following shorthand
to denote the composite 2-morphism obtained by ``rotating'' the upward crossing:
\begin{align}
\begin{tikzpicture}[KM,centerzero]
\draw[to-] (0.3,-0.3) \botlabel{i} -- (-0.3,0.3);
\draw[-to] (-0.3,-0.3) \botlabel{j} -- (0.3,0.3);
\region{0.4,0}{\lambda};        
\end{tikzpicture}\ 
&:=
\begin{tikzpicture}[KM,centerzero,scale=1.2]
\draw[-to] (0.1,-0.3)\botlabel{j} \braidup (-0.1,0.3);
\draw[-to] (-0.4,0.3) -- (-0.4,0.1) to[out=down,in=left] (-0.2,-0.2) to[out=right,in=left] (0.2,0.2) to[out=right,in=up] (0.4,-0.1)  -- (0.4,-0.3)\botlabel{i};
\region{0.6,0}{\lambda};        
\end{tikzpicture}
\ .\label{rightpivot}
\end{align}
\end{itemize}
Now for the relations.
There are three families.
First, we have
the \emph{quiver Hecke algebra} relations:
\begin{align}
\begin{tikzpicture}[KM,centerzero]
\draw[-to] (-0.3,-0.3) \botlabel{i} -- (0.3,0.3);
\draw[-to] (0.3,-0.3) \botlabel{j} -- (-0.3,0.3);
\opendot{-0.15,-0.15};
\end{tikzpicture}
-
\begin{tikzpicture}[KM,centerzero]
\draw[-to] (-0.3,-0.3) \botlabel{i} -- (0.3,0.3);
\draw[-to] (0.3,-0.3) \botlabel{j} -- (-0.3,0.3);
\opendot{0.15,0.15};
\end{tikzpicture}
&= \delta_{i,j}  \ 
\begin{tikzpicture}[KM,centerzero]
\draw[-to] (-0.2,-0.3) \botlabel{i} -- (-0.2,0.3);
\draw[-to] (0.2,-0.3) \botlabel{i} -- (0.2,0.3);
\end{tikzpicture} =    \begin{tikzpicture}[KM,centerzero]
\draw[-to] (-0.3,-0.3) \botlabel{i} -- (0.3,0.3);
\draw[-to] (0.3,-0.3) \botlabel{j} -- (-0.3,0.3);
\opendot{-0.15,0.15};
\end{tikzpicture}
-
\begin{tikzpicture}[KM,centerzero]
\draw[-to] (-0.3,-0.3) \botlabel{i} -- (0.3,0.3);
\draw[-to] (0.3,-0.3) \botlabel{j} -- (-0.3,0.3);
\opendot{0.15,-0.15};
\end{tikzpicture} \ ,\label{dotslide}\\\label{quadratic}
\begin{tikzpicture}[KM,centerzero,scale=1.1]
\draw[-to] (-0.2,-0.4) \botlabel{i} to[out=45,in=down] (0.15,0) to[out=up,in=-45] (-0.2,0.4);
\draw[-to] (0.2,-0.4) \botlabel{j} to[out=135,in=down] (-0.15,0) to[out=up,in=225] (0.2,0.4);
\end{tikzpicture}
&=
\begin{tikzpicture}[KM,centerzero,scale=1.1]
\draw[-to] (-0.2,-0.4) \botlabel{i} -- (-0.2,0.4);
\draw[-to] (0.2,-0.4) \botlabel{j} -- (0.2,0.4);
\pinpin{0.2,0}{-0.2,0}{-1.1,0}{Q_{i,j}(x,y)};
\end{tikzpicture}\ ,\\\label{braid}
\begin{tikzpicture}[KM,centerzero,scale=1.1]
\draw[-to] (-0.4,-0.4) \botlabel{i} -- (0.4,0.4);
\draw[-to] (0,-0.4) \botlabel{j} to[out=135,in=down] (-0.32,0) to[out=up,in=225] (0,0.4);
\draw[-to] (0.4,-0.4) \botlabel{k} -- (-0.4,0.4);
\end{tikzpicture}
\ -\
\begin{tikzpicture}[KM,centerzero,scale=1.1]
\draw[-to] (-0.4,-0.4) \botlabel{i} -- (0.4,0.4);
\draw[-to] (0,-0.4) \botlabel{j} to[out=45,in=down] (0.32,0) to[out=up,in=-45] (0,0.4);
\draw[-to] (0.4,-0.4) \botlabel{k} -- (-0.4,0.4);
\end{tikzpicture}
&=\delta_{i,k}\ 
\begin{tikzpicture}[KM,centerzero,scale=1.1]
\draw[-to] (-0.3,-0.4) \botlabel{i} -- (-0.3,0.4);
\draw[-to] (0,-0.4) \botlabel{j} -- (0,0.4);
\draw[-to] (0.3,-0.4) \botlabel{i} -- (0.3,0.4);
\pinpinpin{.3,0}{0,0}{-.3,0}{-1.7,0}{
\frac{Q_{i,j}(x,y)-Q_{i,j}(z,y)}{x-z}};
\end{tikzpicture}\ .
\end{align}
Next, the \emph{right adjunction relations}
implying that $\one_\lambda F_i$ is right dual to $E_i \one_\lambda$:
\begin{align}\label{rightadj}
\begin{tikzpicture}[KM,centerzero,scale=1.2]
\draw[-to] (-0.3,0.4) -- (-0.3,0) arc(180:360:0.15) arc(180:0:0.15) -- (0.3,-0.4) \botlabel{i};
\end{tikzpicture}
&=
\begin{tikzpicture}[KM,centerzero,scale=1.2]
\draw[to-] (0,-0.4) \botlabel{i} -- (0,0.4);
\end{tikzpicture}
\ ,&
\begin{tikzpicture}[KM,centerzero,scale=1.2]
\draw[-to] (-0.3,-0.4) \botlabel{i}-- (-0.3,0) arc(180:0:0.15) arc(180:360:0.15) -- (0.3,0.4);
\end{tikzpicture}
&=
\begin{tikzpicture}[KM,centerzero,scale=1.2]
\draw[-to] (0,-0.4) \botlabel{i} -- (0,0.4);
\end{tikzpicture}\ ,
\end{align}
Finally, we have the 
\emph{inversion relations} which assert that
\begin{equation}\label{thismatrix}
\begin{tikzpicture}[KM,centerzero]
\draw[-to] (-0.3,-0.3) \botlabel{j}-- (0.3,0.3);
\draw[to-] (0.3,-0.3) \botlabel{i} -- (-0.3,0.3);
\region{0.38,0.02}{\lambda};
\end{tikzpicture}:
E_j F_i \one_\lambda \Rightarrow F_i E_j \one_\lambda\end{equation}
is an isomorphism for all $\lambda \in X$ and $i \neq j$, as are the following matrices for all $\lambda$ and $i$:
\begin{equation}
M_{\lambda;i} :=
\begin{dcases}
\begin{pmatrix}
\begin{tikzpicture}[KM,centerzero]
\draw[-to] (-0.25,-0.25) \botlabel{i}-- (0.25,0.25);
\draw[to-] (0.25,-0.25) \botlabel{i} -- (-0.25,0.25);
\region{0.35,0.02}{\lambda};
\end{tikzpicture} &
\begin{tikzpicture}[KM,centerzero]
\draw[-to] (-0.25,0.15) \toplabel{i} to[out=-90,in=-90,looseness=3] (0.25,0.15);
\region{-0.45,-.1}{\lambda};
\node at (0,.2) {$\phantom.$};\node at (0,-.3) {$\phantom.$};
\end{tikzpicture}
&
\begin{tikzpicture}[KM,centerzero]
\draw[-to] (-0.25,0.15) \toplabel{i} to[out=-90,in=-90,looseness=3] (0.25,0.15);
\region{-0.45,-.1}{\lambda};
\node at (0,.2) {$\phantom.$};\node at (0,-.3) {$\phantom.$};
\opendot{0.23,-0.03};
\end{tikzpicture}
&\!\!\cdots\!\!
&
\begin{tikzpicture}[KM,centerzero]
\draw[-to] (-0.25,0.15) \toplabel{i} to[out=-90,in=-90,looseness=3] (0.25,0.15);
\region{-0.45,-.1}{\lambda};
\node at (0,.2) {$\phantom.$};\node at (0,-.3) {$\phantom.$};
\multopendot{0.23,-0.03}{west}{-h_i(\lambda)-1};
\end{tikzpicture}
\end{pmatrix}\phantom{_T}
&\text{if } h_i(\lambda) \leq 0\\
\begin{pmatrix}   
\begin{tikzpicture}[KM,centerzero]
\draw[-to] (-0.25,-0.25) \botlabel{i}-- (0.25,0.25);
\draw[to-] (0.25,-0.25)\botlabel{i}  -- (-0.25,0.25);
\region{0.38,0.02}{\lambda};
\end{tikzpicture} &
\begin{tikzpicture}[KM,centerzero]
\draw[-to] (-0.25,-0.15) \botlabel{i} to [out=90,in=90,looseness=3](0.25,-0.15);
\region{0.45,0.1}{\lambda};
\node at (0,.3) {$\phantom.$};
\node at (0,-.4) {$\phantom.$};
\end{tikzpicture}
&
\begin{tikzpicture}[KM,centerzero]
\draw[-to] (-0.25,-0.15) \botlabel{i} to [out=90,in=90,looseness=3](0.25,-0.15);
\region{0.45,0.1}{\lambda};
\node at (0,.3) {$\phantom.$};
\node at (0,-.4) {$\phantom.$};
\opendot{-0.23,.03};
\end{tikzpicture}
&
\!\!\!\cdots\!\!\!
&
\begin{tikzpicture}[KM,centerzero]
\draw[-to] (-0.25,-0.15) \botlabel{i} to [out=90,in=90,looseness=3](0.25,-0.15);
\region{0.45,0.1}{\lambda};
\node at (0,.3) {$\phantom.$};
\node at (0,-.4) {$\phantom.$};
\multopendot{-0.23,.03}{east}{h_i(\lambda)-1};
\end{tikzpicture}
\end{pmatrix}^\transpose\phantom{_T}
&\text{if } h_i(\lambda) > 0.
\end{dcases}
\end{equation}

There are also leftward cups and caps 
such that
\begin{align}\label{leftpivots}
\begin{tikzpicture}[KM,centerzero,scale=1.2]
\draw[to-] (-0.3,0.4) -- (-0.3,0) arc(180:360:0.15) arc(180:0:0.15) -- (0.3,-0.4)\botlabel{i};
\end{tikzpicture}
&=
\begin{tikzpicture}[KM,centerzero,scale=1.2]
\draw[-to] (0,-0.4)\botlabel{i} -- (0,0.4);
\end{tikzpicture}
\ ,&
\begin{tikzpicture}[KM,centerzero,scale=1.2]
\draw[to-] (-0.3,-0.4)\botlabel{i} -- (-0.3,0) arc(180:0:0.15) arc(180:360:0.15) -- (0.3,0.4);
\end{tikzpicture}
&=
\begin{tikzpicture}[KM,centerzero,scale=1.2]
\draw[to-] (0,-0.4) \botlabel{i} -- (0,0.4);
\end{tikzpicture}\ .
\end{align}
They define a {\em second adjunction}
making $\one_\lambda F_i$ into a left dual of $E_i \one_\lambda$.
%We have found it helpful to incorporate some extra degree of freedom into the definition of these, which we do by making additional
%choices of {\em normalization homomorphisms}
%\begin{equation}\label{great}
%c_i:X \rightarrow \kk^\times
%\end{equation}
%for all $i \in I$, with $c_i(\alpha_i) = 1$.
%Then we define the following additional 2-morphisms:
They are defined as follows:
\begin{itemize}
\item
Let
$\begin{tikzpicture}[KM,centerzero]
\draw[-to] (0.25,-0.25) \botlabel{j} -- (-0.25,0.25);
\draw[to-] (-0.25,-0.25) \botlabel{i} -- (0.25,0.25);
\region{0.33,0}{\lambda};        
\end{tikzpicture}$
be
$\Big(\begin{tikzpicture}[KM,centerzero]
\draw[to-] (0.25,-0.25)\botlabel{i}  -- (-0.25,0.25);
\draw[-to] (-0.25,-0.25) \botlabel{j} -- (0.25,0.25);
\region{0.33,0}{\lambda};        
\end{tikzpicture}\Big)^{-1}$ if $i \neq j$, 
or the first entry of the matrix $-M_{\lambda;i}^{-1}$ if $i=j$.
\item
Let
$\begin{tikzpicture}[KM,centerzero]
\draw[to-] (-0.25,-0.15) \botlabel{i} to [out=90,in=90,looseness=3](0.25,-0.15);
\region{0.45,0.1}{\lambda};
\node at (0,.3) {$\phantom.$};
\node at (0,-.4) {$\phantom.$};
\end{tikzpicture}$
be the last entry of $%c_i(\lambda)^{-1}
M_{\lambda;i}^{-1}$ if $h_i(\lambda) < 0$ or $-%c_i(\lambda)^{-1}
\left( \begin{tikzpicture}[KM,centerzero,scale=1.2]
 \draw[to-] (-.2,-.3)\botlabel{i} to [out=90,in=-90,looseness=1] (.15,.15) to [out=90,in=90,looseness=1.7] (-.15,.15) to [out=-90,in=90,looseness=1] (.2,-.3);
 \multopendot{-0.14,0.15}{east}{h_i(\lambda)};
\region{0.4,0}{\lambda};
\end{tikzpicture}\right)$ if $h_i(\lambda) \geq 0$.
\item Let $\begin{tikzpicture}[KM,centerzero]
\draw[to-] (-0.25,0.15) \toplabel{i} to[out=-90,in=-90,looseness=3] (0.25,0.15);
\region{0.45,-.1}{\lambda};
\node at (0,.2) {$\phantom.$};\node at (0,-.3) {$\phantom.$};
\end{tikzpicture}$
be the last entry of $%c_i(\lambda)
M_{\lambda;i}^{-1}$ if $h_i(\lambda) > 0$ or
$%c_i(\lambda) 
\begin{tikzpicture}[KM,centerzero,scale=1.2]
\draw[to-] (-.2,.3)\toplabel{i} to [out=-90,in=90,looseness=1] (.15,-.15) to [out=-90,in=-90,looseness=1.7] (-.15,-.15) to [out=90,in=-90,looseness=1] (.2,.3);
\multopendot{0.14,-0.15}{west}{-h_i(\lambda)};
\region{-0.4,0}{\lambda};
\end{tikzpicture}$ if $h_i(\lambda) \leq 0$.
\item
Let $\begin{tikzpicture}[KM,centerzero,scale=1.1]
\draw[to-] (0,-0.4)\botlabel{i} -- (0,0.4);
\opendot{0,0};
\region{0.2,0}{\lambda};
\end{tikzpicture}
:=
\begin{tikzpicture}[KM,centerzero,scale=1.1]
\draw[to-] (0.3,-0.4) \botlabel{i}-- (0.3,0) arc(0:180:0.15) arc(360:180:0.15) -- (-0.3,0.4);
\opendot{0,0};
\region{0.5,0}{\lambda};
\end{tikzpicture}=
 \begin{tikzpicture}[KM,centerzero,scale=1.1]
\draw[-to] (0.3,0.4) -- (0.3,0) arc(360:180:0.15) arc(0:180:0.15) -- (-0.3,-0.4)\botlabel{i};
\opendot{0,0};
\region{0.45,0}{\lambda};
\end{tikzpicture}$.
\item Let
$\begin{tikzpicture}[KM,centerzero]
\draw[to-] (0.3,-0.3) \botlabel{j}-- (-0.3,0.3) ;
\draw[to-] (-0.3,-0.3) \botlabel{i}-- (0.3,0.3) ;
\region{0.4,0}{\lambda};        
\end{tikzpicture} :=%r_{i,j}\
t_{i,j}^{-1}\ 
\begin{tikzpicture}[KM,anchorbase,scale=.7]
\draw[-to] (1.3,.4) to (1.3,-1.2)\botlabel{j};
\draw[-] (-1.3,-.4) to (-1.3,1.2);
\draw[-] (.5,1.1) to [out=0,in=90,looseness=1] (1.3,.4);
\draw[-] (-.35,.4) to [out=90,in=180,looseness=1] (.5,1.1);
\draw[-] (-.5,-1.1) to [out=180,in=-90,looseness=1] (-1.3,-.4);
\draw[-] (.35,-.4) to [out=-90,in=0,looseness=1] (-.5,-1.1);
\draw[-] (.35,-.4) to [out=90,in=-90,looseness=1] (-.35,.4);
\draw[-] (-0.35,-.5) to[out=0,in=180,looseness=1] (0.35,.5);
\draw[-to] (0.35,.5) to[out=0,in=90,looseness=1] (0.8,-1.2)\botlabel{i};
\draw[-] (-0.35,-.5) to[out=180,in=-90,looseness=1] (-0.8,1.2);
\region{1.55,0}{\lambda};
\end{tikzpicture}=
%r_{j,i}\ 
t_{j,i}^{-1}\,
\begin{tikzpicture}[KM,anchorbase,scale=.7]
\draw[-to] (-1.3,.4) to (-1.3,-1.2) \botlabel{i};
\draw[-] (1.3,-.4) to (1.3,1.2);
\draw[-] (-.5,1.1) to [out=180,in=90,looseness=1] (-1.3,.4);
\draw[-] (.35,.4) to [out=90,in=0,looseness=1] (-.5,1.1);
\draw[-] (.5,-1.1) to [out=0,in=-90,looseness=1] (1.3,-.4);
\draw[-] (-.35,-.4) to [out=-90,in=180,looseness=1] (.5,-1.1);
\draw[-] (-.35,-.4) to [out=90,in=-90,looseness=1] (.35,.4);
\draw[-] (0.35,-.5) to[out=180,in=0,looseness=1] (-0.35,.5);
\draw[-to] (-0.35,.5) to[out=180,in=90,looseness=1] (-0.8,-1.2) \botlabel{j};
\draw[-] (0.35,-.5) to[out=0,in=-90,looseness=1] (0.8,1.2);
\region{1.55,0}{\lambda};
\end{tikzpicture}$.
%where
%$r_{i,j} := 
%\begin{cases}
%\frac{c_i(\alpha_j)}{t_{i,j}}
%&\text{if $i \neq j$}\\
%1&\text{if $i=j$.}
%\end{cases}$
\end{itemize}
The equalities in the definitions of the downward dot and crossing just given are by no means obvious; they are justified in \cite{Brundan}.
Equivalently, we have that
\begin{align}\label{lotsmore}
\begin{tikzpicture}[KM,centerzero,scale=1.2]
\draw[-to] (-0.25,0.15) \toplabel{i} to[out=-90,in=-90,looseness=3] (0.25,0.15);
\opendot{-0.22,-0.05};
\end{tikzpicture}
&=
\begin{tikzpicture}[KM,centerzero,scale=1.2]
\draw[-to] (-0.25,0.15) \toplabel{i} to[out=-90,in=-90,looseness=3] (0.25,0.15);
\opendot{0.23,-0.05};
\end{tikzpicture}\ ,&
\begin{tikzpicture}[KM,centerzero,scale=1.2]
\draw[to-] (-0.25,0.15) \toplabel{i} to[out=-90,in=-90,looseness=3] (0.25,0.15);
\opendot{-0.22,-0.05};
\end{tikzpicture}
&=
\begin{tikzpicture}[KM,centerzero,scale=1.2]
\draw[to-] (-0.25,0.15) \toplabel{i} to[out=-90,in=-90,looseness=3] (0.25,0.15);
\opendot{0.23,-0.05};
\end{tikzpicture}\ ,&
  \begin{tikzpicture}[KM,centerzero,scale=1.2]
\draw[-to] (-0.25,-0.15) \botlabel{i} to[out=90,in=90,looseness=3] (0.25,-0.15);
\opendot{-0.22,0.05};
\end{tikzpicture}
&=
\begin{tikzpicture}[KM,centerzero,scale=1.2]
\draw[-to] (-0.25,-0.15) \botlabel{i} to[out=90,in=90,looseness=3] (0.25,-0.15);
\opendot{0.23,0.05};
\end{tikzpicture}\ ,&
\begin{tikzpicture}[KM,centerzero,scale=1.2]
\draw[to-] (-0.25,-0.15) \botlabel{i} to[out=90,in=90,looseness=3] (0.25,-0.15);
\opendot{-0.22,0.05};
\end{tikzpicture}
&=
\begin{tikzpicture}[KM,centerzero,scale=1.2]
\draw[to-] (-0.25,-0.15) \botlabel{i} to[out=90,in=90,looseness=3] (0.25,-0.15);
\opendot{0.23,0.05};
\end{tikzpicture}\ ,\\\label{ruby}
\begin{tikzpicture}[KM,anchorbase,scale=1.2]
\draw[-to] (-0.25,0.15) \toplabel{i} to[out=-90,in=-90,looseness=3] (0.25,0.15);
\draw[-to] (-0.3,-0.4) to[out=up,in=down] (0,0.15)\toplabel{j};
\end{tikzpicture}
&=
\begin{tikzpicture}[KM,anchorbase,scale=1.2]
\draw[-to] (-0.25,0.15) \toplabel{i} to[out=-90,in=-90,looseness=3] (0.25,0.15);
\draw[-to] (0.3,-0.4)to[out=up,in=down] (0,0.15)\toplabel{j};
\end{tikzpicture}\ ,&
\begin{tikzpicture}[KM,anchorbase,scale=1.2]
\draw[to-] (-0.25,0.15) \toplabel{i} to[out=-90,in=-90,looseness=3] (0.25,0.15);
\draw[to-] (-0.3,-0.4) to[out=up,in=down] (0,0.15)\toplabel{j};
\end{tikzpicture}
&=
\begin{tikzpicture}[KM,anchorbase,scale=1.2]
\draw[to-] (-0.25,0.15) \toplabel{i} to[out=-90,in=-90,looseness=3] (0.25,0.15);
\draw[to-] (0.3,-0.4)to[out=up,in=down] (0,0.15)\toplabel{j};
\end{tikzpicture}\ ,
&
\begin{tikzpicture}[KM,anchorbase,scale=1.2]
\draw[-to] (-0.25,0.15) \toplabel{i} to[out=-90,in=-90,looseness=3] (0.25,0.15);
\draw[to-] (-0.3,-0.4) to[out=up,in=down] (0,0.15)\toplabel{j};
\end{tikzpicture}
&=%r_{j,i}\ 
t_{j,i}^{-1}\ 
\begin{tikzpicture}[KM,anchorbase,scale=1.2]
\draw[-to] (-0.25,0.15) \toplabel{i} to[out=-90,in=-90,looseness=3] (0.25,0.15);
\draw[to-] (0.3,-0.4)to[out=up,in=down] (0,0.15)\toplabel{j};
\end{tikzpicture}\ ,&
\begin{tikzpicture}[KM,anchorbase,scale=1.2]
\draw[to-] (-0.25,0.15) \toplabel{i} to[out=-90,in=-90,looseness=3] (0.25,0.15);
\draw[-to] (-0.3,-0.4) to[out=up,in=down] (0,0.15)\toplabel{j};
\end{tikzpicture}
&=
%r_{i,j}^{-1}\ 
t_{i,j}\ 
\begin{tikzpicture}[KM,anchorbase,scale=1.2]
\draw[to-] (-0.25,0.15) \toplabel{i} to[out=-90,in=-90,looseness=3] (0.25,0.15);
\draw[-to] (0.3,-0.4)to[out=up,in=down] (0,0.15)\toplabel{j};
\end{tikzpicture}\ ,\\\label{wax}
\begin{tikzpicture}[KM,anchorbase,scale=1.2]
\draw[-to] (-0.25,-0.15) \botlabel{i} to[out=90,in=90,looseness=3] (0.25,-0.15);
\draw[to-] (-0.3,0.4) \braiddown (0,-0.15)\botlabel{j};
\end{tikzpicture}
&=
\begin{tikzpicture}[KM,anchorbase,scale=1.2]
\draw[-to] (-0.25,-0.15) \botlabel{i} to[out=90,in=90,looseness=3] (0.25,-0.15);
\draw[to-] (0.3,0.4) \braiddown (0,-0.15)\botlabel{j};
\end{tikzpicture}\ ,&
\begin{tikzpicture}[KM,baseline=-1mm,scale=1.2]
\draw[to-] (-0.25,-0.15) \botlabel{i} to[out=90,in=90,looseness=3] (0.25,-0.15);
\draw[-to] (-0.3,0.4) \braiddown (0,-0.15)\botlabel{j};
\end{tikzpicture}
&=
\begin{tikzpicture}[KM,baseline=-1mm,scale=1.2]
\draw[to-] (-0.25,-0.15) \botlabel{i} to[out=90,in=90,looseness=3] (0.25,-0.15);
\draw[-to] (0.3,0.4) \braiddown (0,-0.15)\botlabel{j};
\end{tikzpicture}\ ,
&\begin{tikzpicture}[KM,baseline=-1mm,scale=1.2]
\draw[-to] (-0.25,-0.15) \botlabel{i} to[out=90,in=90,looseness=3] (0.25,-0.15);
\draw[-to] (-0.3,0.4) \braiddown (0,-0.15)\botlabel{j};
\end{tikzpicture}
&=%r_{j,i}^{-1}\ 
t_{j,i}\ 
\begin{tikzpicture}[KM,baseline=-1mm,scale=1.2]
\draw[-to] (-0.25,-0.15) \botlabel{i} to[out=90,in=90,looseness=3] (0.25,-0.15);
\draw[-to] (0.3,0.4) \braiddown (0,-0.15)\botlabel{j};
\end{tikzpicture}\ ,&
\begin{tikzpicture}[KM,anchorbase,scale=1.2]
\draw[to-] (-0.25,-0.15) \botlabel{i} to[out=90,in=90,looseness=3] (0.25,-0.15);
\draw[to-] (-0.3,0.4) \braiddown (0,-0.15)\botlabel{j};
\end{tikzpicture}
&=%r_{i,j}\ 
t_{i,j}^{-1}\ 
\begin{tikzpicture}[KM,anchorbase,scale=1.2]
\draw[to-] (-0.25,-0.15) \botlabel{i} to[out=90,in=90,looseness=3] (0.25,-0.15);
\draw[to-] (0.3,0.4)\braiddown (0,-0.15)\botlabel{j};
\end{tikzpicture}\ .
\end{align}
%One can always choose the normalization homomorphisms so that $c_i = 1$ for all $i$, in which case $r_{i,j} = t_{i,j}^{-1}$ for $i \neq j$.
%However,
%the ``pitchfork relations'' just recorded are more reasonable
%if one can choose the normalization homomorphisms so that
%$c_i(\alpha_j) = t_{i,j}$ for all $i \neq j$, for then $r_{i,j} =1$ always and the string diagrams are invariant under planar isotopy.
%Unfortunately, it might not be possible to find such a choice (a sufficient condition is that the simple roots are linearly independent); this subtlety was not explained correctly in \cite{HKM}.

Like in the general theory of symmetric functions, when working with $\UU$, it is often helpful to work in terms of generating functions, which in general will be formal Laurent series in auxiliary variables $u^{-1}, v^{-1},\dots$.
For such a generating function $f(u)$,
we use the notation
$[f(u)]_{u:n}$ to denote its $u^{n}$-coefficient,
$[f(u)]_{u;\geq 0}$ for its polynomial part, 
and so on.
For a polynomial $f(x)$, we have that
\begin{align}\label{trick}
\left[\frac{f(u)}{u-x}\right]_{u:-1}&=f(x),&
\left[\frac{f(u)}{u-x}\right]_{u:<0}&=\frac{f(x)}{u-x}.
\end{align}
A trivial but perhaps countertuitive consequence is that
\begin{equation}\label{trickconsequence}
\prod_{i=1}^n \left[\frac{f_i(u)}{u-x}\right]_{u:-1}
= \left[\frac{\prod_{i=1}^n f_i(u)}{u-x}
\right]_{u:-1}
\end{equation}
for $f_1(x),\dots,f_n(x) \in \kk[x]$.
We view the series
\begin{align}\label{often}
\frac{1}{u-x} &= \sum_{r \geq 0} x^r u^{-r-1}.
\end{align}
as a generating function for multiple dots on a string, introducing the additional shorthand
\begin{align}\label{dgf}
\begin{tikzpicture}[KM,centerzero,scale=1.1]
\draw[-] (0,-0.3) -- (0,0.3);
\circled{0,0}{u};
\end{tikzpicture}
&:=
\begin{tikzpicture}[KM,centerzero,scale=1.1]
\draw[-] (0,-0.3) -- (0,0.3);
\pin{0,0}{.7,0}{\frac{1}{u-x}};
\end{tikzpicture}\
.
\end{align}
The following is a consequence of \cref{dotslide}:
\begin{align}\label{gendotslide}
\begin{tikzpicture}[KM,centerzero,scale=1.2]
\draw[-to] (-0.3,-0.3) \botlabel{i} -- (0.3,0.3);
\draw[-to] (0.3,-0.3) \botlabel{j} -- (-0.3,0.3);
\circled{-0.15,-0.15}{u};
\end{tikzpicture}
-
\begin{tikzpicture}[KM,centerzero,scale=1.3]
\draw[-to] (-0.3,-0.3) \botlabel{i} -- (0.3,0.3);
\draw[-to] (0.3,-0.3) \botlabel{j} -- (-0.3,0.3);
\circled{0.15,0.15}{u};
\end{tikzpicture}
&= \delta_{i,j}  \ 
\begin{tikzpicture}[KM,centerzero,scale=1.3]
\draw[-to] (-0.2,-0.3) \botlabel{i} -- (-0.2,0.3);
\draw[-to] (0.2,-0.3) \botlabel{i} -- (0.2,0.3);
\circled{-0.2,0}{u};
\circled{0.2,0}{u};
\end{tikzpicture} =    
\begin{tikzpicture}[KM,centerzero,scale=1.3]
\draw[-to] (-0.3,-0.3) \botlabel{i} -- (0.3,0.3);
\draw[-to] (0.3,-0.3) \botlabel{j} -- (-0.3,0.3);
\circled{-0.15,0.15}{u};
\end{tikzpicture}
-
\begin{tikzpicture}[KM,centerzero,scale=1.3]
\draw[-to] (-0.3,-0.3) \botlabel{i} -- (0.3,0.3);
\draw[-to] (0.3,-0.3) \botlabel{j} -- (-0.3,0.3);
\circled{0.15,-0.15}{u};
\end{tikzpicture} \ .
\end{align}
There are also the extremely useful {\em bubble generating functions}
\begin{align}
\label{bubblegeneratingfunction1}
\begin{tikzpicture}[KM,baseline=-1mm]
\draw[to-] (-0.25,0) arc(180:-180:0.25);
\node at (0,-.4) {\strandlabel{i}};
\region{0.95,0}{\lambda};
\node at (.55,0) {$(u)$};
\end{tikzpicture}
&\in %c_i(\lambda)^{-1}
u^{h_i(\lambda)}\id_{\one_\lambda}
+ u^{h_i(\lambda)-1} \kk\llbracket u^{-1}\rrbracket\End_{\UU}(\one_\lambda),\\
\begin{tikzpicture}[KM,baseline=-1mm]
\draw[-to] (-0.25,0) arc(180:-180:0.25);
\node at (0,-.4) {\strandlabel{i}};
\region{0.95,0}{\lambda};
\node at (.55,0) {$(u)$};
\end{tikzpicture}
&\in %c_i(\lambda) 
u^{-h_i(\lambda)}\id_{\one_\lambda}
+ u^{-h_i(\lambda)-1} \kk\llbracket u^{-1}\rrbracket\End_{\UU}(\one_\lambda),
\label{bubblegeneratingfunction2}
\end{align}
which are the unique formal Laurent series with leading coefficients as in \cref{bubblegeneratingfunction1,bubblegeneratingfunction2}
such that
\begin{align}\label{miles}
\left[\ \begin{tikzpicture}[KM,baseline=-1mm]
\draw[to-] (-0.25,0) arc(180:-180:0.25);
\node at (0,-.4) {\strandlabel{i}};
\region{0.95,0}{\lambda};
\node at (.55,0) {$(u)$};
\end{tikzpicture}\right]_{u:<0}\!\!
&= 
\begin{tikzpicture}[KM,baseline=-1mm]
\draw[to-] (-0.25,0) arc(180:-180:0.25);
\node at (0,-0.4) {\strandlabel{i}};
\region{0.6,0}{\lambda};
\circled{.25,0}{u};
\end{tikzpicture},
&
\left[\ \begin{tikzpicture}[KM,baseline=-1mm]
\draw[-to] (-0.25,0) arc(180:-180:0.25);
\node at (0,-.4) {\strandlabel{i}};
\region{0.95,0}{\lambda};
\node at (.55,0) {$(u)$};
\end{tikzpicture}\right]_{u:<0}\!\!
&= 
\begin{tikzpicture}[KM,baseline=-1mm]
\draw[-to] (0.25,0) arc(360:0:0.25);
\node at (0,-0.4) {\strandlabel{i}};
\region{0.5,0}{\lambda};
\circled{-.25,0}{u};
\end{tikzpicture},
\end{align}
and the {\em infinite Grassmannian relation}
\begin{align}
\label{infgrass}
\begin{tikzpicture}[KM,centerzero,scale=1]
\draw[to-] (-0.68,0) arc(180:-180:0.25);
\node[black] at (0.1,0) {$(u)$};
\node at (-.4,-.43) {\strandlabel{i}};
\region{.6,0}{\lambda};
\end{tikzpicture}\: 
\begin{tikzpicture}[KM,centerzero,scale=1]
\draw[-to] (-.25,0) arc(180:-180:0.25);
\node[black] at (0.54,0) {$(u)$};
\node at (0,-0.4) {\strandlabel{i}};
\end{tikzpicture} &= 1
\end{align}
holds.
If $h_i(\lambda) \geq 0$ then the clockwise bubble generating function is completely determined by \cref{bubblegeneratingfunction2,miles}, and the counterclockwise bubble generating function is 
the two-sided inverse of the clockwise one.
On the other hand if $h_i(\lambda) \leq 0$ then the counterclockwise bubble generating function is completely determined by \cref{bubblegeneratingfunction1,miles}, and the clockwise one is its inverse.
It is a non-trivial consequence of the defining relations of $\UU$ that such series exist; see \cite{Brundan}.
Here are some further relations involving the bubble generating functions, which explains their importance:
\begin{align}
\label{bubslide}
\begin{tikzpicture}[KM,anchorbase,scale=.9]
\draw[-to] (-0.6,-0.5)\botlabel{j} to (-0.6,0.5);
\draw[to-] (-0.25,0) arc(180:-180:0.25);
\node at (0,-.4) {\strandlabel{i}};
\node at (.55,0) {$(u)$};
\end{tikzpicture}
&=
\begin{tikzpicture}[KM,anchorbase,scale=.9]
\draw[to-] (-0.25,0) arc(180:-180:0.25);
\node at (0,-.4) {\strandlabel{i}};
\node at (.55,0) {$(u)$};
\draw[-to] (1,-0.5)\botlabel{j} to (1,0.5);
\pin{1,0}{2.1,0}{R_{i, j}(u,x)};
\end{tikzpicture}\ ,
&\begin{tikzpicture}[KM,anchorbase,scale=.9]
\draw[-to] (-0.25,0) arc(180:-180:0.25);
\node at (0,-0.4) {\strandlabel{i}};
\node at (.55,0) {$(u)$};
\draw[-to] (1,-0.5)\botlabel{j} to (1,0.5);
\end{tikzpicture}
&=
\begin{tikzpicture}[KM,anchorbase,scale=.9]
\draw[-to] (-0.6,-0.5)\botlabel{j} to (-0.6,0.5);
\draw[-to] (-0.25,0) arc(180:-180:0.25);
\node at (0,-.4) {\strandlabel{i}};
\node at (.55,0) {$(u)$};
\pin{-.6,0}{-1.7,0}{R_{i, j}(u,x)};\end{tikzpicture}
\\\intertext{where $
R_{i,j}(x,y) := \begin{cases}
%r_{i,j}
t_{i,j}^{-1} Q_{i,j}(x,y)
&\text{if $i \neq j$}\\
(x-y)^{-2}&\text{if $i=j$,}
\end{cases}
$}
\begin{tikzpicture}[KM,anchorbase,scale=1.1]
\draw[-to] (0,-0.5)\botlabel{i} to[out=up,in=180] (0.3,0.2) to[out=0,in=up] (0.45,0) to[out=down,in=0] (0.3,-0.2) to[out=180,in=down] (0,0.5);
\circled{.42,0}{u};
\end{tikzpicture}
&=
-\left[\ 
\begin{tikzpicture}[KM,anchorbase,scale=1.1]
\draw[-to] (-0.8,-0.5)\botlabel{i} -- (-0.8,0.5);
\circled{-0.8,0}{u};
\draw[-to] (-.4,0) arc(180:-180:0.2);
\node at (0.26,0) {$(u)$};
\node at (-.17,-.33) {\strandlabel{i}};
\end{tikzpicture}
\right]_{u:< 0},
&
\begin{tikzpicture}[KM,anchorbase,scale=1.1]
\draw[-to] (0,-0.5)\botlabel{i} to[out=up,in=0] (-0.3,0.2) to[out=180,in=up] (-0.45,0) to[out=down,in=180] (-0.3,-0.2) to[out=0,in=down] (0,0.5);
\circled{-.42,0}{u};
\end{tikzpicture}
&=
\left[
\ \begin{tikzpicture}[KM,anchorbase,scale=1.1]
\draw[-to] (1.2,-0.5)\botlabel{i} -- (1.2,0.5);
\circled{1.2,0}{u};
\draw[to-] (0,0) arc(180:-180:0.2);
\node at (0.65,0) {$(u)$};
\node at (0.2,-.33) {\strandlabel{i}};
\end{tikzpicture}\ 
\right]_{u:< 0},
\label{curlrels}\\
\label{altquadratic}
\begin{tikzpicture}[KM,centerzero,scale=1.2]
\draw[-to] (-0.2,-0.4) \botlabel{i} to[out=45,in=down] (0.15,0) to[out=up,in=-45] (-0.2,0.4);
\draw[to-] (0.2,-0.4) \botlabel{j} to[out=135,in=down] (-0.15,0) to[out=up,in=225] (0.2,0.4);
\end{tikzpicture}
&=
(-1)^{\delta_{i,j}}\begin{tikzpicture}[KM,centerzero,scale=1.4]
\draw[-to] (-0.14,-0.3) \botlabel{i} -- (-0.14,0.3);
\draw[to-] (0.14,-0.3) \botlabel{j} -- (0.14,0.3);
\end{tikzpicture}
+
\delta_{i,j}
\left[\,\begin{tikzpicture}[KM,centerzero,scale=1.6]
\draw[-to] (-0.2,-0.3) \botlabel{i} to [looseness=2.2,out=90,in=90] (0.2,-0.3);
\draw[-to] (0.2,0.3) to [looseness=2.2,out=-90,in=-90] (-0.2,0.3)\toplabel{i};
\circled{-.15,.12}{u};
\circled{-.15,-.12}{u};
\draw[to-] (0.27,0) arc(180:-180:0.142);
\node at (0.73,0) {$(u)$};
\node at (0.42,-.23) {\strandlabel{i}};
\end{tikzpicture}\,\right]_{u:-1}\!\!\!\!\!,
&
\begin{tikzpicture}[KM,centerzero,scale=1.2]
\draw[to-] (-0.2,-0.4) \botlabel{i} to[out=45,in=down] (0.15,0) to[out=up,in=-45] (-0.2,0.4);
\draw[-to] (0.2,-0.4) \botlabel{j} to[out=135,in=down] (-0.15,0) to[out=up,in=225] (0.2,0.4);
\end{tikzpicture}
&=
(-1)^{\delta_{i,j}}\begin{tikzpicture}[KM,centerzero,scale=1.4]
\draw[to-] (-0.14,-0.3) \botlabel{i} -- (-0.14,0.3);
\draw[-to] (0.14,-0.3) \botlabel{j} -- (0.14,0.3);
\end{tikzpicture}
+
\delta_{i,j}
\left[\,\begin{tikzpicture}[KM,centerzero,scale=1.6]
\draw[-to] (0.2,-0.3)  to [looseness=2.2,out=90,in=90] (-0.2,-0.3)\botlabel{i};
\draw[-to] (-0.2,0.3)\toplabel{i} to [looseness=2.2,out=-90,in=-90] (0.2,0.3);
\circled{.15,.12}{u};
\circled{.15,-.12}{u};
\draw[to-] (-0.8,0) arc(-180:180:0.142);
\node at (-0.34,0) {$(u)$};
\node at (-0.638,-.23) {\strandlabel{i}};
\end{tikzpicture}\,\right]_{u:-1}\!\!\!\!\!\!,
\end{align}
Again, these follow from relations derived in \cite{Brundan} (also \cite{KL3}), although it is some work to translate into the generating function form (see \cite{HKM}).

\begin{rem}\label{notquite}
There is a $\Z$-grading on $\UU$ defined by declaring that the degrees of the generating 2-morphisms are as listed in \cref{table1}. 
Taking it into account, one can introduce a free $\Z[q,q^{-1}]$-algebra $K_0(\UU_q)$, the split Grothendieck ring of the graded Karoubi envelope of $\UU$. It is 
isomorphic to Lusztig's modified integral form $\dot \U_{\Z[q,q^{-1}]}$ of the quantized enveloping algebra of the same type as $\dot \U$. This result is due to Khovanov and Lauda \cite{KL3}, although it depends on the Non-degeneracy Conjecture formulated therein which was proved later; see \cite[Sec.~3]{unfurling} or \cite[Sec.~5.2]{BWWquasisplit}.
The graded version of $\UU$ will not play a role in the remainder of this article.
\end{rem}

\begin{defin}\label{kmcatdef}
A {\em Kac-Moody categorification} with the type and parameters fixed above
is a locally finite\footnote{This means that objects have finite length and morphism spaces are finite-dimensional as vector spaces over $\kk$.} $\kk$-linear 
Abelian category $\catR$ with a decomposition
\begin{equation}
\catR = \bigoplus_{\lambda \in X} \catR_\lambda
\end{equation}
as the internal direct sum\footnote{If $\catC$ is any additive category
and $\catC_i\:(i \in I)$ is a family of 
full subcategories, $\catC$ is their {\em internal direct sum}, denoted
$\catC = \bigoplus_{i \in I} \catC_i$,
if the subcategories $\catC_i$ are mutually orthogonal
and, for each $X \in \ob\catC$, there are objects 
$X_i \in \ob\catC_i\:(i \in I)$, all but finitely many of which are zero, such that $X$ is a biproduct of $X_i\:(i \in I)$.
Denoting the biproduct inclusions and projections by
$\iota_{i,X}:X_i \hookrightarrow X$ and
$\pi_{i,X}:X \twoheadrightarrow X_i$, there is an associated 
projection functor $\operatorname{Pr}_i:\catC \rightarrow \catC_i$ mapping object $X$ to $X_i$ and morphism $f:X \rightarrow Y$ to $\pi_{i,Y} \circ f\circ \iota_{i,X}$. It is convenient to assume that $\iota_{i,X_i} = \pi_{i,X_i} = \id_{X_i}$ so that $\operatorname{Pr}_i \circ \operatorname{Inc}_i = \id_{\catC_i}$, where $\operatorname{Inc}_i:\catC_i\rightarrow \catC$ is the inclusion functor.}
of Serre subcategories $\catR_\lambda\:(\lambda \in X)$,
the {\em weight subcategories} of $\catR$,
plus adjoint pairs $(E_i, F_i)$ of $\kk$-linear endofunctors 
for each $i \in I$,
such that:
\begin{enumerate}
\item[(KM0)]
The functor $E_i$ takes objects of $\catR_\lambda$ to $\catR_{\lambda+\alpha_i}$; equivalently, $F_i$ takes objects of $\catR_{\lambda+\alpha_i}$ to $\catR_{\lambda}$.
\item[(KM1)] The adjoint pair $(E_i, F_i)$ has a prescribed 
adjunction
with unit and counit of adjunction denoted
 $\begin{tikzpicture}[KM,centerzero,scale=.8]
\draw[-to] (-0.25,0.15) \toplabel{i} to[out=-90,in=-90,looseness=3] (0.25,0.15);
\node at (0,.2) {$\phantom.$};\node at (0,-.3) {$\phantom.$};
\end{tikzpicture}:\id_\catR \Rightarrow F_i \circ E_i$
and $\begin{tikzpicture}[KM,centerzero,scale=.8]
\draw[-to] (-0.25,-0.15) \botlabel{i} to [out=90,in=90,looseness=3](0.25,-0.15);
\node at (0,.3) {$\phantom.$};
\node at (0,-.4) {$\phantom.$};
\end{tikzpicture}:E_i \circ F_i \Rightarrow \id_\catR$.
\item[(KM2)]\label{km2} There are
given natural transformations
$\begin{tikzpicture}[KM,centerzero]
\draw[-to] (0,-0.2) \botlabel{i} -- (0,0.2);
\opendot{0,0};
\end{tikzpicture}:E_i \Rightarrow E_i$ and $\begin{tikzpicture}[KM,centerzero,scale=.9]
\draw[-to] (-0.2,-0.2) \botlabel{i} -- (0.2,0.2);
\draw[-to] (0.2,-0.2) \botlabel{j} -- (-0.2,0.2);
\end{tikzpicture}:E_i \circ E_j \Rightarrow E_j \circ E_i$
satisfying the QHA relations \cref{dotslide,quadratic,braid}.
\item[(KM3)]
Defining $\begin{tikzpicture}[KM,centerzero]
\draw[to-] (0.2,-0.2) \botlabel{i} -- (-0.2,0.2);
\draw[-to] (-0.2,-0.2) \botlabel{j} -- (0.2,0.2);
\end{tikzpicture}$
as in \cref{rightpivot}, the natural transformations
$\begin{tikzpicture}[KM,centerzero]
\draw[to-] (0.2,-0.2) \botlabel{i} -- (-0.2,0.2);
\draw[-to] (-0.2,-0.2) \botlabel{j} -- (0.2,0.2);
\end{tikzpicture}:E_j \circ F_i \Rightarrow F_i \circ E_j$
are isomorphisms for all $i \neq j$, as are the matrices
\begin{align*}
\begin{pmatrix}
\begin{tikzpicture}[KM,centerzero]
\draw[-to] (-0.25,-0.25) \botlabel{i}-- (0.25,0.25);
\draw[to-] (0.25,-0.25) \botlabel{i} -- (-0.25,0.25);
\end{tikzpicture} &
\begin{tikzpicture}[KM,centerzero]
\draw[-to] (-0.25,0.15) \toplabel{i} to[out=-90,in=-90,looseness=3] (0.25,0.15);
\node at (0,.2) {$\phantom.$};\node at (0,-.3) {$\phantom.$};
\end{tikzpicture}
&
\begin{tikzpicture}[KM,centerzero]
\draw[-to] (-0.25,0.15) \toplabel{i} to[out=-90,in=-90,looseness=3] (0.25,0.15);
\node at (0,.2) {$\phantom.$};\node at (0,-.3) {$\phantom.$};
\opendot{0.23,-0.03};
\end{tikzpicture}
&\!\!\cdots\!\!
&
\begin{tikzpicture}[KM,centerzero]
\draw[-to] (-0.25,0.15) \toplabel{i} to[out=-90,in=-90,looseness=3] (0.25,0.15);
\node at (0,.2) {$\phantom.$};\node at (0,-.3) {$\phantom.$};
\multopendot{0.23,-0.03}{west}{h_i(-\lambda)-1};
\end{tikzpicture}
\end{pmatrix}\phantom{_T}&:E_i \circ F_i|_{\catR_\lambda}
\oplus \id_{\catR_\lambda}^{\oplus h_i(-\lambda)}
\Rightarrow F_i \circ E_i|_{\catR_\lambda}
&&\text{if } h_i(\lambda) \leq 0,\text{ or}\\
\begin{pmatrix}   
\begin{tikzpicture}[KM,centerzero]
\draw[-to] (-0.25,-0.25) \botlabel{i}-- (0.25,0.25);
\draw[to-] (0.25,-0.25)\botlabel{i}  -- (-0.25,0.25);
\end{tikzpicture} &
\begin{tikzpicture}[KM,centerzero]
\draw[-to] (-0.25,-0.15) \botlabel{i} to [out=90,in=90,looseness=3](0.25,-0.15);
\node at (0,.3) {$\phantom.$};
\node at (0,-.4) {$\phantom.$};
\end{tikzpicture}
&
\begin{tikzpicture}[KM,centerzero]
\draw[-to] (-0.25,-0.15) \botlabel{i} to [out=90,in=90,looseness=3](0.25,-0.15);
\node at (0,.3) {$\phantom.$};
\node at (0,-.4) {$\phantom.$};
\opendot{-0.23,.03};
\end{tikzpicture}
&
\!\!\!\cdots\!\!\!
&
\begin{tikzpicture}[KM,centerzero]
\draw[-to] (-0.25,-0.15) \botlabel{i} to [out=90,in=90,looseness=3](0.25,-0.15);
\node at (0,.3) {$\phantom.$};
\node at (0,-.4) {$\phantom.$};
\multopendot{-0.23,.03}{east}{h_i(\lambda)-1};
\end{tikzpicture}
\end{pmatrix}^\transpose\phantom{_T}
&:E_i \circ F_i|_{\catR_\lambda}
\Rightarrow F_i \circ E_i|_{\catR_\lambda}
\oplus \id_{\catR_\lambda}^{\oplus h_i(\lambda)}
&&\text{if } h_i(\lambda) > 0,
\end{align*}
for all $i \in I$ and $\lambda \in X$.
\end{enumerate}
This data is equivalent to that of
a strict $\kk$-linear 2-functor
from $\UU$ to the $\kk$-linear 2-category of locally finite $\kk$-linear 
Abelian categories.
\end{defin}

\cref{kmcatdef} may seem quite complicated.
However, there are many examples. Perhaps the most celebrated one 
comes from 
finite-dimensional 
representations of the symmetric groups $S_n\:(n \geq 0)$ with $E_i$ and
$F_i$ being the so-called $i$-induction and $i$-restriction
functors. This fits naturally into the general framework of
Heisenberg categorification introduced in the next section, and will be discussed further in
\cref{s5-example}. Other examples arising from rational
representations of the general linear group are developed in
\cref{s6-example}.

The work of Chuang and Rouquier \cite{CR}, and Rouquier's subsequent work \cite{Rou} which formulated the definition in the manner just stated, establishes some impressive and useful
structural results about Kac-Moody categorifications. 
Since it provides the essential motivation for the material in the rest of the article, we conclude the section by
listing some of these properties. 
A more detailed account can be found in the survey article \cite{BD}, to which we defer for the appropriate citations to the literature.
Let $\catR$ be a Kac-Moody categorification as above.
\begin{enumerate}
\item \label{noods1} For $n \geq 1$, the axiom (KM2) implies that there is a homomorphism $\NH_n \rightarrow \End_\kk(E_i^n)$, where $\NH_n$ is the nil-Hecke algebra over $\kk$. Letting $e_n \in \NH_n$ be some preferred choice of a homogeneous primitive idempotent, its image under this homomorphism defines a summand $E_i^{(n)}$ of the functor $E_i^n$ such that
$E_i^n \cong (E_i^{(n)})^{\oplus n!}$.
Similarly, there is a summand $F_i^{(n)}$ of $F_i^n$ such that
$F_i^n \cong (F_i^{(n)})^{\oplus n!}$, which may be chosen so that the given adjunction $(E_i, F_i)$ induces an adjunction $(E_i^{(n)}, F_i^{(n)})$.
\item \label{noods2} As well as being right adjoint to $E_i^{(n)}$, the functor $F_i^{(n)}$ is also left adjoint to $E_i^{(n)}$ with unit and counit induced by the second adjunction described above. 
Hence, $E_i^{(n)}$ and $F_i^{(n)}$ are exact and they preserve injective and projective objects.
\item \label{noods3} Let $G_0(\catR)$ be the Grothendieck group of the Abelian category $\catR$. It is naturally a module over the algebra $\dot\U_\Z$ with weight decomposition $G_0(\catR) = \bigoplus_{\lambda \in X} G_0(\catR_\lambda)$. The actions of $e_i^{(n)} 1_\lambda$ and $f_i^{(n)} 1_\lambda$ are induced by the exact functors $E_i^{(n)}|_{\catR_\lambda}$ and $F_i^{(n)}|_{\catR_\lambda}$. In fact, $G_0(\catR)$ is an {\em integrable} $\dot\U_\Z$-module: for each $V \in \ob\catR$ there is $n \in \N$ such that $E_i^{(n)} V = F_i^{(n)} V = 0$. Moreover, for each $V \in \ob\catR$, we have that $E_i V = F_i V = 0$ for all but finitely many $i \in I$.
\end{enumerate}
For the remaining results in this summary, we assume 
for simplicity that $\catR$ is a {\em nilpotent} 
categorification\footnote{
There are similar results without this nilpotency assumption, but one needs to ``unfurl'' the Cartan matrix to obtain a larger Kac-Moody algebra. See \cite[Sec.~4]{unfurling} which discusses the general case.},
that is, the endomorphism of $E_i V$ defined by
$\begin{tikzpicture}[KM,centerzero]
\draw[-to] (0,-0.2) \botlabel{i} -- (0,0.2);
\opendot{0,0};
\end{tikzpicture}$
is nilpotent for all $V \in \ob\catR$ and $i \in I$;
equivalently, the endomorphism of $F_i V$ defined by
$\begin{tikzpicture}[KM,centerzero]
\draw[to-] (0,-0.2) \botlabel{i} -- (0,0.2);
\opendot{0,0};
\end{tikzpicture}$
is nilpotent for all $V$ and $i$.
\begin{enumerate}[resume]
\item \label{noods4} Let $\B = \bigsqcup_{\lambda \in X} \B_\lambda$
be a set such that $\B_\lambda$ indexes a full set $L(b)\:(b \in \B_\lambda)$ of pairwise inequivalent irreducible objects of $\catR_\lambda$. For $b \in \B$, 
let $\eps_i(b)$ and $\phi_i(b)$
be the nilpotency degrees of the endomorphisms
$x:E_i L(b)\rightarrow E_i L(b)$ and $y:F_i L(b)\rightarrow F_i L(b)$ defined by the natural transformations
$\begin{tikzpicture}[KM,centerzero]
\draw[-to] (0,-0.2) \botlabel{i} -- (0,0.2);
\opendot{0,0};
\end{tikzpicture}$
and $\begin{tikzpicture}[KM,centerzero]
\draw[to-] (0,-0.2) \botlabel{i} -- (0,0.2);
\opendot{0,0};
\end{tikzpicture}$, respectively.
Then the endomorphism algebras
$\End_{\catR}(E_i L(b))$ and $\End_{\catR}(F_i L(b))$ are 
isomorphic to $\kk[x] / \left(x^{\eps_i(b)}\right)$ and $\kk[y] / \left(y^{\phi_i(b)}\right)$, respectively.
In particular,
$E_i L(b)$ is non-zero if and only if
$\eps_i(b) \neq 0$, and
$F_i L(b)$ is non-zero if and only if $\phi_i(b) \neq 0$.
\item \label{noods5}
For $b \in \B_\lambda$,
Schur's Lemma implies that the bubble generating functions
act on $L(b)$ by series in $\kk\lround u^{-1}\rround$.
In fact, $\ \begin{tikzpicture}[KM,centerzero,scale=.8]
\draw[to-] (-0.68,0) arc(180:-180:0.25);
\node[black] at (0.16,0) {$(u)$};
\node at (-.4,-.41) {\strandlabel{i}};
\end{tikzpicture}\!$ acts on $L(b)$ by $%c_i(\lambda)^{-1} 
u^{\phi_i(b)-\eps_i(b)}$
and
$\ \begin{tikzpicture}[KM,centerzero,scale=.8]
\draw[-to] (-0.68,0) arc(180:-180:0.25);
\node[black] at (0.16,0) {$(u)$};
\node at (-.4,-.41) {\strandlabel{i}};
\end{tikzpicture}\!$ acts by $%c_i(\lambda) 
u^{\eps_i(b)-\phi_i(b)}$.
Hence:
\begin{equation}\label{abhilash}
h_i(\lambda) = \phi_i(b)-\eps_i(b).
\end{equation}
\item \label{noods6}
For $b \in \B$ with $\eps_i(b) > 0$,
$E_i L(b)$ is an indecomposable object
with irreducible socle and head
isomorphic to $L\big(\tilde e_i(b)\big)$ for some $\tilde e_i(b) \in \B$; moreover, $\eps_i(\tilde e_i(b)) = \eps_i(b) - 1$
and $\phi_i(\tilde e_i(b)) = \phi_i(b) + 1$.
Similarly, if $\phi_i(b) > 0$ then $F_i L(b)$ is indecomposable with irreducible socle and head isomorphic to $L\big(\tilde f_i(b)\big)$ for some $\tilde f_i(b) \in \B$; 
moreover, $\phi_i(\tilde f_i(b)) = \phi_i(b) - 1$
and $\eps_i(\tilde f_i(b)) = \eps_i(b) + 1$.
This defines mutually inverse bijections
$$
\begin{tikzcd}
\{b \in \B\:|\:\eps_i(b) > 0\}
\arrow[r,"\tilde e_i" above,yshift=3pt]&
\arrow[l,"\tilde f_i" below,yshift=-3pt]
\{b \in \B\:|\:\phi_i(b) > 0\}.
\end{tikzcd}
$$
Also defining
$\wt:\B \rightarrow X$ so that $\wt(b) = \lambda\Leftrightarrow b \in \B_\lambda$,
the datum 
$(\B, \tilde e_i, \tilde f_i, \eps_i, \phi_i, \wt)$ is a {\em normal crystal} in the sense\footnote{Note though that we prefer to set things up so that domains of $\tilde e_i$ and $\tilde f_i$ are subsets of $\B$, leaving them undefined on the elements $b \in \B$ with $\eps_i(b) = 0$ or $\phi_i(b)=0$, respectively.} 
of \cite[Sec.~7.6]{Kashcrystal}.
\item \label{noods7} For $b \in \B_\lambda$ and $0 \leq m \leq n := \eps_i(b)$,
the object $E_i^{(m)} L(b)$ is indecomposable with
simple socle and head isomorphic to $L\big(\tilde e_i^m(b)\big)$.
The composition multiplicity $[E_i^{(m)} L(b): L\big(\tilde e_i^m(b)\big)]$ equals $\binom{n}{m}$, and
all other composition factors are of the form $L(b')$ for $b' \in \B_{\lambda + m \alpha_i}$
with $\eps_i(b') < n - m$. Also 
the endomorphism algebra $\End_{\catR}(E_i^{(m)} L(b))$ is explicitly determined; it is isomorphic to the cohomology algebra
$H^*(\operatorname{Gr}_{m,n},\kk)$ of the Grassmannian. 
There are analogous statements about $F_i^{(m)} L(b)$
for $0 \leq m \leq n:=\phi_i(b)$.
\item \label{noods8}
The isomorphism classes $[L(b)]\:(b \in \B)$
give a {\em perfect basis} for the integrable $\dot\U_\Z$-module
$G_0(\catR)$ in the sense of \cite[Def.~5.49]{BeK}.
\end{enumerate}
Let $W$ be the Weyl group associated to the Cartan matrix $A$.
It is generated by simple reflections $s_i\:(i \in I)$,
and acts on $X$ so that 
$s_i(\lambda) = \lambda - h_i(\lambda) \alpha_i$.
Since $G_0(\catR)$ is an integrable $\dot\U_\Z$-module,
the set of weights $\lambda \in X$ such that $\B_\lambda \neq \varnothing$ is a union of $W$-orbits.
\begin{enumerate}[resume]
\item \label{noods9}
Suppose that $\lambda \in X$ and $i \in I$ such that
$\eps_i(b) = 0$ for all $b \in \B_\lambda$.
By \cref{abhilash}, we have that $n := h_i(\lambda) = \phi_i(b)$ for any $b \in \B_\lambda$.
Then the functor $F_i^{(n)}:\catR_\lambda \rightarrow \catR_{s_i(\lambda)}$ is an equivalence of categories with quasi-inverse defined by $E_i^{(n)}$. 
\item \label{noods10} Take $\lambda \in X$ and $i \in I$
such that $m := h_i(\lambda) \geq 0$.
There is a complex of functors
$$
\cdots \longrightarrow E_i^{(2)} F_i^{(m+2)} \longrightarrow E_i F_i^{(m+1)} \longrightarrow F_i^{(m)} \longrightarrow 0
$$
which induces a triangulated functor $D^b(\catR_\lambda)
\rightarrow D^b(\catR_{s_i(\lambda)})$ between bounded derived categories. In \cite{CR}, Chuang and Rouquier proved that this is an equivalence of categories; in fact, there is a perverse equivalence
between $\catR_\lambda$ and $\catR_{s_i(\lambda)}$ in the sense of \cite{CRother}.
Hence, weight subcategories of $\catR$ for weights lying in the same $W$-orbit are derived equivalent.
\end{enumerate}

\begin{rem}\label{schurian}
We have formulated the definition of Kac-Moody categorification in the context of locally finite $\kk$-linear 
Abelian categories. There is a parallel theory 
in which $\catR$ is merely a $\kk$-linear 
Karoubian category with finite-dimensional morphism spaces. 
There are analogous results to the above in this finite-dimensional Karoubian setting. One needs to replace $G_0(\catR)$ with the split Grothendieck group $K_0(\catR)$, which again has the induced structure of an integrable $\dot\U_\Z$-module. The isomorphism classes of indecomposable objects of $\catR$ define a {\em dual perfect basis} for $K_0(\catR)$ in the sense of \cite{KKKS}. This is discussed further in \cite{BD}.
\end{rem}

%% file: s3-H.tex
\setcounter{section}{2}

%=====================================
\section{Heisenberg categorifications}\label{s3-H}
%=====================================

Recall that $\kk$ is an algebraically closed field.
Fix a {\em central charge} $\kappa \in \Z$. The {\em degenerate Heisenberg category} $\Heis_\kappa$ is a strict $\kk$-linear monoidal category defined by generators and relations.
It was introduced originally by
Khovanov \cite{K} for $\kappa=\pm 1$, then the definition was extended to all non-zero central charges in \cite{MS}, while for $\kappa=0$ it is the
affine oriented Brauer category from \cite{BCNR}.
We define it here following the approach of
\cite[Theorem 1.2]{Bheis}, which is similar in spirit to the definition of $\UU$ given in the previous section.
It has two
generating objects $E$ and $F$, whose identity endomorphisms are represented by the oriented strings $\ \begin{tikzpicture}[H,anchorbase]
\draw[-to] (0,-0.2) -- (0,0.2);
\end{tikzpicture}\ $ and
$\ \begin{tikzpicture}[H,anchorbase]
\draw[to-] (0,-0.2) -- (0,0.2);
\end{tikzpicture}\ $,
and the following four generating morphisms:
\begin{align*}
\begin{tikzpicture}[H,centerzero]
\draw[-to] (0,-0.3) -- (0,0.3);
\opendot{0,0};
\end{tikzpicture}\,
&\colon E \rightarrow E,
&
\begin{tikzpicture}[H,centerzero,scale=.9]
\draw[-to] (-0.3,-0.3) -- (0.3,0.3);
\draw[-to] (0.3,-0.3) -- (-0.3,0.3);
\end{tikzpicture}
&\colon E\otimes E  \rightarrow E \otimes E,
&
\begin{tikzpicture}[H,baseline=1mm]
\draw[-to] (-0.15,-0.05) to [out=90,in=90,looseness=3](0.35,-0.05);
\end{tikzpicture}\ \ 
&\colon E\otimes F \rightarrow \one,
&
\begin{tikzpicture}[H,baseline=-1.3mm]
\draw[-to] (-0.25,0.15) to[out=-90,in=-90,looseness=3] (0.25,0.15);
\end{tikzpicture}\ \ &
\colon \one \rightarrow F\otimes E.
\end{align*}
These are subject to certain relations.
First, we have
the \emph{degenerate affine Hecke algebra} relations:
\begin{align}
\begin{tikzpicture}[H,centerzero]
\draw[-to] (-0.3,-0.3) -- (0.3,0.3);
\draw[-to] (0.3,-0.3) -- (-0.3,0.3);
\opendot{-0.15,-0.15};
\end{tikzpicture}
-
\begin{tikzpicture}[H,centerzero]
\draw[-to] (-0.3,-0.3) -- (0.3,0.3);
\draw[-to] (0.3,-0.3) -- (-0.3,0.3);
\opendot{0.15,0.15};
\end{tikzpicture}
&= 
\begin{tikzpicture}[H,centerzero]
\draw[-to] (-0.2,-0.3) -- (-0.2,0.3);
\draw[-to] (0.2,-0.3)  -- (0.2,0.3);
\end{tikzpicture} =    
\begin{tikzpicture}[H,centerzero]
\draw[-to] (-0.3,-0.3)  -- (0.3,0.3);
\draw[-to] (0.3,-0.3) -- (-0.3,0.3);
\opendot{-0.15,0.15};
\end{tikzpicture}
-
\begin{tikzpicture}[H,centerzero]
\draw[-to] (-0.3,-0.3) -- (0.3,0.3);
\draw[-to] (0.3,-0.3) -- (-0.3,0.3);
\opendot{0.15,-0.15};
\end{tikzpicture} \ ,\label{Hdotslide}\\\label{Hquadratic}
\begin{tikzpicture}[H,centerzero,scale=1.1]
\draw[-to] (-0.2,-0.4) to[out=45,in=down] (0.15,0) to[out=up,in=-45] (-0.2,0.4);
\draw[-to] (0.2,-0.4) to[out=135,in=down] (-0.15,0) to[out=up,in=225] (0.2,0.4);
\end{tikzpicture}
&=
\begin{tikzpicture}[H,centerzero,scale=1.1]
\draw[-to] (-0.15,-0.4) -- (-0.15,0.4);
\draw[-to] (0.15,-0.4) -- (0.15,0.4);
\end{tikzpicture}\ ,\\\label{Hbraid}
\begin{tikzpicture}[H,centerzero,scale=1.1]
\draw[-to] (-0.4,-0.4) -- (0.4,0.4);
\draw[-to] (0,-0.4) to[out=135,in=down] (-0.32,0) to[out=up,in=225] (0,0.4);
\draw[-to] (0.4,-0.4)  -- (-0.4,0.4);
\end{tikzpicture}
&=
\begin{tikzpicture}[H,centerzero,scale=1.1]
\draw[-to] (-0.4,-0.4) -- (0.4,0.4);
\draw[-to] (0,-0.4) to[out=45,in=down] (0.32,0) to[out=up,in=-45] (0,0.4);
\draw[-to] (0.4,-0.4) -- (-0.4,0.4);
\end{tikzpicture}\ .
\end{align}
Next, the \emph{right adjunction relations}
implying that $F$ is right dual to $E$:
\begin{align}\label{Hrightadj}
\begin{tikzpicture}[H,centerzero,scale=1.2]
\draw[-to] (-0.3,0.4) -- (-0.3,0) arc(180:360:0.15) arc(180:0:0.15) -- (0.3,-0.4);
\end{tikzpicture}
&=
\begin{tikzpicture}[H,centerzero,scale=1.2]
\draw[to-] (0,-0.4) -- (0,0.4);
\end{tikzpicture}
\ ,&
\begin{tikzpicture}[H,centerzero,scale=1.2]
\draw[-to] (-0.3,-0.4) -- (-0.3,0) arc(180:0:0.15) arc(180:360:0.15) -- (0.3,0.4);
\end{tikzpicture}
&=
\begin{tikzpicture}[H,centerzero,scale=1.2]
\draw[-to] (0,-0.4) -- (0,0.4);
\end{tikzpicture}\ .
\end{align}
Finally, we have the 
\emph{inversion relation} which asserts that the 
morphism
\begin{equation}\label{rybnikov}
M_\kappa :=
\begin{dcases}
\begin{pmatrix}
\begin{tikzpicture}[H,anchorbase]
\draw[-to] (-0.25,-0.25) -- (0.25,0.25);
\draw[to-] (0.25,-0.25) -- (-0.25,0.25);
\end{tikzpicture} &
\begin{tikzpicture}[H,anchorbase]
\draw[-to] (-0.25,0.15) to[out=-90,in=-90,looseness=3] (0.25,0.15);
\node at (0,.2) {$\phantom.$};\node at (0,-.3) {$\phantom.$};
\end{tikzpicture}
&
\begin{tikzpicture}[H,anchorbase]
\draw[-to] (-0.25,0.15) to[out=-90,in=-90,looseness=3] (0.25,0.15);
\node at (0,.2) {$\phantom.$};\node at (0,-.3) {$\phantom.$};
\opendot{0.23,-0.03};
\end{tikzpicture}
&\!\!\cdots\!\!
&
\begin{tikzpicture}[H,anchorbase]
\draw[-to] (-0.25,0.15) to[out=-90,in=-90,looseness=3] (0.25,0.15);
\node at (0,.2) {$\phantom.$};\node at (0,-.3) {$\phantom.$};
\multopendot{0.23,-0.03}{west}{-\kappa-1};
\end{tikzpicture}
\end{pmatrix}\phantom{_T}:E \otimes F \oplus \one^{\oplus (-\kappa)}
\rightarrow F \otimes E
&\text{if } \kappa \leq 0\\
\begin{pmatrix}   
\begin{tikzpicture}[H,centerzero]
\draw[-to] (-0.25,-0.25) -- (0.25,0.25);
\draw[to-] (0.25,-0.25)  -- (-0.25,0.25);
\end{tikzpicture} &
\begin{tikzpicture}[H,centerzero]
\draw[-to] (-0.25,-0.15) to [out=90,in=90,looseness=3](0.25,-0.15);
\node at (0,.3) {$\phantom.$};
\node at (0,-.4) {$\phantom.$};
\end{tikzpicture}
&
\begin{tikzpicture}[H,centerzero]
\draw[-to] (-0.25,-0.15) to [out=90,in=90,looseness=3](0.25,-0.15);
\node at (0,.3) {$\phantom.$};
\node at (0,-.4) {$\phantom.$};
\opendot{-0.23,.03};
\end{tikzpicture}
&
\!\!\!\cdots\!\!\!
&
\begin{tikzpicture}[H,centerzero]
\draw[-to] (-0.25,-0.15) to [out=90,in=90,looseness=3](0.25,-0.15);
\node at (0,.3) {$\phantom.$};
\node at (0,-.4) {$\phantom.$};
\multopendot{-0.23,.03}{east}{\kappa-1};
\end{tikzpicture}\ 
\end{pmatrix}^\transpose\phantom{_T}:E \otimes F \rightarrow 
F \otimes E \oplus \one^{\oplus \kappa}
&\text{if } \kappa > 0
\end{dcases}
\end{equation}
in the additive envelope of $\Heis_\kappa$ is invertible, where the rightward crossing is defined by
\begin{align}
\begin{tikzpicture}[H,centerzero]
\draw[to-] (0.3,-0.3) -- (-0.3,0.3);
\draw[-to] (-0.3,-0.3) -- (0.3,0.3);
\end{tikzpicture}\ 
&:=
\begin{tikzpicture}[H,centerzero,scale=1.2]
\draw[-to] (0.1,-0.3) \braidup (-0.1,0.3);
\draw[-to] (-0.4,0.3) -- (-0.4,0.1) to[out=down,in=left] (-0.2,-0.2) to[out=right,in=left] (0.2,0.2) to[out=right,in=up] (0.4,-0.1)  -- (0.4,-0.3);
\end{tikzpicture}
\ .\label{Hrightpivot}
\end{align}

\begin{rem}\label{muckycats}
Recall that the {\em Demazure operator}
$\partial_{xy}:\kk[x,y] \rightarrow \kk[x,y]$ is the linear map with
\begin{equation}\label{demazure}
\partial_{xy} f(x,y) := \frac{f(x,y)-f(y,x)}{x-y}.
\end{equation}
The dot slide relation \cref{Hdotslide} implies that
\begin{align} \label{naha}
\begin{tikzpicture}[H,centerzero]
\draw[-to] (-0.4,-0.4)  -- (0.4,0.4);
\draw[-to] (0.4,-0.4)  -- (-0.4,0.4);
\pinpin{0.25,-0.25}{-0.25,-0.25}{-1.1,-.25}{f(x,y)};
\end{tikzpicture}\ -\ 
\begin{tikzpicture}[H,centerzero]
\draw[-to] (-0.4,-0.4) -- (0.4,0.4);
\draw[-to] (0.4,-0.4) -- (-0.4,0.4);
\pinpin{-0.25,0.25}{0.25,0.25}{1.1,.25}{f(y,x)};
\end{tikzpicture}
&=
\begin{tikzpicture}[H,centerzero]
\draw[-to] (-0.3,-0.4) -- (-0.3,0.4);
\draw[-to] (0.3,-0.4) -- (0.3,0.4);
\pinpin{0.3,0}{-0.3,0}{-1.3,0}{\partial_{xy} f(x,y)};
\end{tikzpicture}\ 
\end{align}
for any $f(x,y) \in \kk[x,y]$.
\end{rem}

Again like in the previous section,
there is also a
{\em second adjunction}
making $F$ into a left dual of $E$, that is, 
there are leftward cups and caps 
such that
\begin{align}\label{Hleftpivots}
\begin{tikzpicture}[H,centerzero,scale=1.2]
\draw[to-] (-0.3,0.4) -- (-0.3,0) arc(180:360:0.15) arc(180:0:0.15) -- (0.3,-0.4);
\end{tikzpicture}
&=
\begin{tikzpicture}[H,centerzero,scale=1.2]
\draw[-to] (0,-0.4) -- (0,0.4);
\end{tikzpicture}
\ ,&
\begin{tikzpicture}[H,centerzero,scale=1.2]
\draw[to-] (-0.3,-0.4) -- (-0.3,0) arc(180:0:0.15) arc(180:360:0.15) -- (0.3,0.4);
\end{tikzpicture}
&=
\begin{tikzpicture}[H,centerzero,scale=1.2]
\draw[to-] (0,-0.4) -- (0,0.4);
\end{tikzpicture}\ .
\end{align}
These, and some other useful shorthands, are defined as follows:
\begin{itemize}
\item
Let
$\begin{tikzpicture}[H,centerzero]
\draw[-to] (0.25,-0.25)  -- (-0.25,0.25);
\draw[to-] (-0.25,-0.25)-- (0.25,0.25);
\end{tikzpicture}$
be the first entry of the inverse matrix $M_\kappa^{-1}$.
\item
Let
$\ \begin{tikzpicture}[H,centerzero]
\draw[to-] (-0.25,-0.15)  to [out=90,in=90,looseness=3](0.25,-0.15);
\end{tikzpicture}\ $
be the last entry of $
M_\kappa^{-1}$ if $\kappa < 0$ or $
 \begin{tikzpicture}[H,centerzero,scale=1.2]
 \draw[to-] (-.2,-.3) to [out=90,in=-90,looseness=1] (.15,.15) to [out=90,in=90,looseness=1.7] (-.15,.15) to [out=-90,in=90,looseness=1] (.2,-.3);
 \multopendot{-0.14,0.15}{east}{\kappa};
\end{tikzpicture}$ if $\kappa \geq 0$.
\item Let $\begin{tikzpicture}[H,centerzero]
\draw[to-] (-0.25,0.15) to[out=-90,in=-90,looseness=3] (0.25,0.15);
\end{tikzpicture}$
be the last entry of $-M_{\kappa}^{-1}$ if $\kappa> 0$ or
$\begin{tikzpicture}[H,centerzero,scale=1.2]
\draw[to-] (-.2,.3) to [out=-90,in=90,looseness=1] (.15,-.15) to [out=-90,in=-90,looseness=1.7] (-.15,-.15) to [out=90,in=-90,looseness=1] (.2,.3);
\multopendot{0.14,-0.15}{west}{-\kappa};
\end{tikzpicture}$ if $\kappa \leq 0$.
\item
Let $\begin{tikzpicture}[H,centerzero,scale=1.1]
\draw[to-] (0,-0.4) -- (0,0.4);
\opendot{0,0};
\end{tikzpicture}
:=
\begin{tikzpicture}[H,centerzero,scale=1.1]
\draw[to-] (0.3,-0.4) -- (0.3,0) arc(0:180:0.15) arc(360:180:0.15) -- (-0.3,0.4);
\opendot{0,0};
\end{tikzpicture}=
 \begin{tikzpicture}[H,centerzero,scale=1.1]
\draw[-to] (0.3,0.4) -- (0.3,0) arc(360:180:0.15) arc(0:180:0.15) -- (-0.3,-0.4);
\opendot{0,0};
\end{tikzpicture}\ $.
\item Let
$\begin{tikzpicture}[H,centerzero]
\draw[to-] (0.3,-0.3) -- (-0.3,0.3);
\draw[to-] (-0.3,-0.3) -- (0.3,0.3);
\end{tikzpicture} :=
\begin{tikzpicture}[H,anchorbase,scale=.7]
\draw[-to] (1.3,.4) to (1.3,-1.2);
\draw[-] (-1.3,-.4) to (-1.3,1.2);
\draw[-] (.5,1.1) to [out=0,in=90,looseness=1] (1.3,.4);
\draw[-] (-.35,.4) to [out=90,in=180,looseness=1] (.5,1.1);
\draw[-] (-.5,-1.1) to [out=180,in=-90,looseness=1] (-1.3,-.4);
\draw[-] (.35,-.4) to [out=-90,in=0,looseness=1] (-.5,-1.1);
\draw[-] (.35,-.4) to [out=90,in=-90,looseness=1] (-.35,.4);
\draw[-] (-0.35,-.5) to[out=0,in=180,looseness=1] (0.35,.5);
\draw[-to] (0.35,.5) to[out=0,in=90,looseness=1] (0.8,-1.2);
\draw[-] (-0.35,-.5) to[out=180,in=-90,looseness=1] (-0.8,1.2);
\end{tikzpicture}=
\begin{tikzpicture}[H,anchorbase,scale=.7]
\draw[-to] (-1.3,.4) to (-1.3,-1.2);
\draw[-] (1.3,-.4) to (1.3,1.2);
\draw[-] (-.5,1.1) to [out=180,in=90,looseness=1] (-1.3,.4);
\draw[-] (.35,.4) to [out=90,in=0,looseness=1] (-.5,1.1);
\draw[-] (.5,-1.1) to [out=0,in=-90,looseness=1] (1.3,-.4);
\draw[-] (-.35,-.4) to [out=-90,in=180,looseness=1] (.5,-1.1);
\draw[-] (-.35,-.4) to [out=90,in=-90,looseness=1] (.35,.4);
\draw[-] (0.35,-.5) to[out=180,in=0,looseness=1] (-0.35,.5);
\draw[-to] (-0.35,.5) to[out=180,in=90,looseness=1] (-0.8,-1.2);
\draw[-] (0.35,-.5) to[out=0,in=-90,looseness=1] (0.8,1.2);
\end{tikzpicture}\ $.
\end{itemize}
The equalities in the definitions of the downward dot and crossing  are justified in \cite{Bheis}.
It follows that string diagrams for morphisms in $\Heis_\kappa$ are
invariant under planar isotopy. Hence, $\Heis_\kappa$ is strictly pivotal.

We use the pin notation and generating functions as explained in the previous section. In particular, we have the dot generating function
\begin{align}\label{Hdgf}
\begin{tikzpicture}[H,centerzero,scale=1.1]
\draw[-] (0,-0.3) -- (0,0.3);
\circled{0,0}{u};
\end{tikzpicture}
&:=
\begin{tikzpicture}[H,centerzero,scale=1.1]
\draw[-] (0,-0.3) -- (0,0.3);
\pin{0,0}{.7,0}{\frac{1}{u-x}};
\end{tikzpicture}\ ,
\end{align}
which satisfies
\begin{align}\label{Hgendotslide}
\begin{tikzpicture}[H,centerzero,scale=1.2]
\draw[-to] (-0.3,-0.3) -- (0.3,0.3);
\draw[-to] (0.3,-0.3) -- (-0.3,0.3);
\circled{-0.15,-0.15}{u};
\end{tikzpicture}
-
\begin{tikzpicture}[H,centerzero,scale=1.3]
\draw[-to] (-0.3,-0.3) -- (0.3,0.3);
\draw[-to] (0.3,-0.3) -- (-0.3,0.3);
\circled{0.15,0.15}{u};
\end{tikzpicture}
&= 
\begin{tikzpicture}[H,centerzero,scale=1.3]
\draw[-to] (-0.2,-0.3) -- (-0.2,0.3);
\draw[-to] (0.2,-0.3) -- (0.2,0.3);
\circled{-0.2,0}{u};
\circled{0.2,0}{u};
\end{tikzpicture} =    
\begin{tikzpicture}[H,centerzero,scale=1.3]
\draw[-to] (-0.3,-0.3) -- (0.3,0.3);
\draw[-to] (0.3,-0.3) -- (-0.3,0.3);
\circled{-0.15,0.15}{u};
\end{tikzpicture}
-
\begin{tikzpicture}[H,centerzero,scale=1.3]
\draw[-to] (-0.3,-0.3) -- (0.3,0.3);
\draw[-to] (0.3,-0.3) -- (-0.3,0.3);
\circled{0.15,-0.15}{u};
\end{tikzpicture} \ .
\end{align}
In the Heisenberg setting, 
the
{\em bubble generating functions}\footnote{Note that the clockwise bubble generating function here is $-1$ times the one in \cite{HKM}---with hindsight we think this is a more consistent convention.} are the unique formal Laurent series
\begin{align}
\label{Hbubblegeneratingfunction1}
\begin{tikzpicture}[H,baseline=-1mm]
\draw[to-] (-0.25,0) arc(180:-180:0.25);
\node at (.55,0) {\color{black}$(u)$};
\end{tikzpicture}
&\in 
u^{\kappa}\id_{\one}
+ u^{\kappa-1} \kk\llbracket u^{-1}\rrbracket\End_{\Heis_\kappa}(\one),\\
\begin{tikzpicture}[H,baseline=-1mm]
\draw[-to] (-0.25,0) arc(180:-180:0.25);
\node at (.55,0) {$\color{black}(u)$};
\end{tikzpicture}
&\in -u^{-\kappa}\id_{\one}
+ u^{-\kappa-1} \kk\llbracket u^{-1}\rrbracket\End_{\Heis_\kappa}(\one)
\label{Hbubblegeneratingfunction2}
\end{align}
such that
\begin{align}\label{Hmiles}
\left[\ \begin{tikzpicture}[H,baseline=-1mm]
\draw[to-] (-0.25,0) arc(180:-180:0.25);
\node at (.55,0) {$\color{black}(u)$};
\end{tikzpicture}\right]_{u:<0}\!\!
&= 
\begin{tikzpicture}[H,baseline=-1mm]
\draw[to-] (-0.25,0) arc(180:-180:0.25);
\circled{.25,0}{u};
\end{tikzpicture},
&
\left[\ \begin{tikzpicture}[H,baseline=-1mm]
\draw[-to] (-0.25,0) arc(180:-180:0.25);
\node at (.55,0) {$\color{black}(u)$};
\end{tikzpicture}\right]_{u:<0}\!\!
&= 
\begin{tikzpicture}[H,baseline=-1mm]
\draw[-to] (0.25,0) arc(360:0:0.25);
\circled{-.25,0}{u};
\end{tikzpicture}\ ,&
\begin{tikzpicture}[H,centerzero,scale=1]
\draw[to-] (-0.68,0) arc(180:-180:0.25);
\node[black] at (0.1,0) {$\color{black}(u)$};
\end{tikzpicture} 
\begin{tikzpicture}[H,centerzero,scale=1]
\draw[-to] (-.25,0) arc(180:-180:0.25);
\node[black] at (0.54,0) {$\color{black}(u)$};
\end{tikzpicture} &= -1.
\end{align}
Existence of these series is justified in \cite{Bheis}.
The Heisenberg analogues of \cref{bubslide,curlrels,altquadratic} are
\begin{align}
\label{Hbubslide}
\begin{tikzpicture}[H,anchorbase,scale=.9]
\draw[-to] (-0.6,-0.5) to (-0.6,0.5);
\draw[to-] (-0.25,0) arc(180:-180:0.25);
\node at (.55,0) {$\color{black}(u)$};
\end{tikzpicture}
&=
\begin{tikzpicture}[H,anchorbase,scale=.9]
\draw[to-] (-0.25,0) arc(180:-180:0.25);
\node at (.55,0) {$\color{black}(u)$};
\draw[-to] (1,-0.5) to (1,0.5);
\pin{1,0}{1.9,0}{R(u,x)};
\end{tikzpicture}\ ,
&\begin{tikzpicture}[H,anchorbase,scale=.9]
\draw[-to] (-0.25,0) arc(180:-180:0.25);
\node at (.55,0) {$\color{black}(u)$};
\draw[-to] (1,-0.5) to (1,0.5);
\end{tikzpicture}
&=
\begin{tikzpicture}[H,anchorbase,scale=.9]
\draw[-to] (-0.6,-0.5) to (-0.6,0.5);
\draw[-to] (-0.25,0) arc(180:-180:0.25);
\node at (.55,0) {$\color{black}(u)$};
\pin{-.6,0}{-1.5,0}{R(u,x)};\end{tikzpicture}
\ \\\intertext{where $R(x,y) := 
1-(x-y)^{-2}=
\frac{(x-y-1)(x-y+1)}{(x-y)^2}$,}
\begin{tikzpicture}[H,anchorbase,scale=1.1]
\draw[-to] (0,-0.5) to[out=up,in=180] (0.3,0.2) to[out=0,in=up] (0.45,0) to[out=down,in=0] (0.3,-0.2) to[out=180,in=down] (0,0.5);
\circled{.42,0}{u};
\end{tikzpicture}
&=
-\left[\ 
\begin{tikzpicture}[H,anchorbase,scale=1.1]
\draw[-to] (-0.8,-0.5) -- (-0.8,0.5);
\circled{-0.8,0}{u};
\draw[-to] (-.4,0) arc(180:-180:0.2);
\node at (0.26,0) {$\color{black}(u)$};
\end{tikzpicture}
\right]_{u:< 0},
&
\begin{tikzpicture}[H,anchorbase,scale=1.1]
\draw[-to] (0,-0.5) to[out=up,in=0] (-0.3,0.2) to[out=180,in=up] (-0.45,0) to[out=down,in=180] (-0.3,-0.2) to[out=0,in=down] (0,0.5);
\circled{-.42,0}{u};
\end{tikzpicture}
&=
\left[
\ \begin{tikzpicture}[H,anchorbase,scale=1.1]
\draw[-to] (1.2,-0.5) -- (1.2,0.5);
\circled{1.2,0}{u};
\draw[to-] (0,0) arc(180:-180:0.2);
\node at (0.65,0) {$\color{black}(u)$};
\end{tikzpicture}\ 
\right]_{u:< 0},
\label{Hcurlrels}\\
\label{Haltquadratic}
\begin{tikzpicture}[H,centerzero,scale=1.2]
\draw[-to] (-0.2,-0.4) to[out=45,in=down] (0.15,0) to[out=up,in=-45] (-0.2,0.4);
\draw[to-] (0.2,-0.4) to[out=135,in=down] (-0.15,0) to[out=up,in=225] (0.2,0.4);
\end{tikzpicture}
&=\begin{tikzpicture}[H,centerzero,scale=1.4]
\draw[-to] (-0.14,-0.3) -- (-0.14,0.3);
\draw[to-] (0.14,-0.3) -- (0.14,0.3);
\end{tikzpicture}
+
\left[\ \begin{tikzpicture}[H,centerzero,scale=1.6]
\draw[-to] (-0.2,-0.3) to [looseness=2.2,out=90,in=90] (0.2,-0.3);
\draw[-to] (0.2,0.3) to [looseness=2.2,out=-90,in=-90] (-0.2,0.3);
\circled{-.15,.12}{u};
\circled{-.15,-.12}{u};
\draw[to-] (0.27,0) arc(180:-180:0.142);
\node at (0.73,0) {$\color{black}(u)$};
\end{tikzpicture}\ \right]_{u:-1}\!\!\!\!\!,
&
\begin{tikzpicture}[H,centerzero,scale=1.2]
\draw[to-] (-0.2,-0.4) to[out=45,in=down] (0.15,0) to[out=up,in=-45] (-0.2,0.4);
\draw[-to] (0.2,-0.4) to[out=135,in=down] (-0.15,0) to[out=up,in=225] (0.2,0.4);
\end{tikzpicture}
&=\begin{tikzpicture}[H,centerzero,scale=1.4]
\draw[to-] (-0.14,-0.3) -- (-0.14,0.3);
\draw[-to] (0.14,-0.3) -- (0.14,0.3);
\end{tikzpicture}
+
\left[\ \begin{tikzpicture}[H,centerzero,scale=1.6]
\draw[-to] (0.2,-0.3)  to [looseness=2.2,out=90,in=90] (-0.2,-0.3);
\draw[-to] (-0.2,0.3) to [looseness=2.2,out=-90,in=-90] (0.2,0.3);
\circled{.15,.12}{u};
\circled{.15,-.12}{u};
\draw[to-] (-0.8,0) arc(-180:180:0.142);
\node at (-0.34,0) {$\color{black}(u)$};
\end{tikzpicture}\ \right]_{u:-1}\!\!\!\!\!.
\end{align}
These identities are copied from \cite[Sec.~3.1]{HKM}.

\begin{rem}
Let $K_0(\Heis_k)$ be the Grothendieck ring of the additive Karoubi envelope of $\Heis_\kappa$.
When the characteristic of the ground field is zero, $K_0(\Heis_\kappa)$
is isomorphic to the {\em Heisenberg ring}
$\heis_\kappa$, that is, the ring
generated by elements $\{h_n^+, e_n^-\:|\:n \geq 0\}$ subject to
the relations
\begin{align*}
h_0^+&=e_0^-=1,&
h_m^+ h_n^+ &= h_n^+ h_m^+,&
e_m^- e_n^- &= e_n^- e_m^-,&
  h_m^+ e_n^- &=
\sum_{r=0}^{\min(m,n)}
                  \binom{\,\kappa\,}{\,r\,}\:e_{n-r}^- h_{m-r}^+.
\end{align*}
This ring is a $\Z$-form
for the universal enveloping algebra of the infinite-dimensional
Heisenberg Lie algebra specialized at central charge $\kappa$.
The existence of an isomorphism
$K_0(\Heis_\kappa) \cong \heis_\kappa$ was conjectured by Khovanov
in \cite{K} (for $\kappa=\pm 1$), and proved in \cite{K0}.
Under the isomorphism,
the classes $[E], [F] \in K_0(\Heis_\kappa)$
correspond to $h_1^+, e_1^- \in \heis_\kappa$.
More generally,
$h_n^+, e_n^-$ are the isomorphism classes of the summands of $E^n$ and $F^n$ defined by
idempotents arising from the trivial and sign representations
of the symmetric group $S_n$.
\end{rem}

\begin{defin}\label{hcatdef}
A {\em degenerate Heisenberg categorification} of central charge $\kappa$ 
is a locally finite $\kk$-linear Abelian category $\catR$
plus an adjoint pair $(E, F)$ of $\kk$-linear endofunctors
such that:
\begin{itemize}
\item[(H1)]
The adjoint pair $(E,F)$ has a prescribed adjunction with unit 
and counit of adjunction denoted
$\;\begin{tikzpicture}[H,baseline=-4pt,scale=.8]
\draw[-to] (-0.25,0.15) to[out=-90,in=-90,looseness=3] (0.25,0.15);
\end{tikzpicture}\;:\id_\catR \Rightarrow F \circ E$
and
$\;\begin{tikzpicture}[H,baseline=-2pt,scale=.8]
\draw[-to] (-0.25,-0.15) to [out=90,in=90,looseness=3](0.25,-0.15);
\end{tikzpicture}\;:E\circ F \Rightarrow \id_\catR$.
\item[(H2)]
There are given natural transformations
$\begin{tikzpicture}[H,centerzero]
\draw[-to] (0,-0.2) -- (0,0.2);
\opendot{0,0};
\end{tikzpicture}:E \Rightarrow E$ and $\begin{tikzpicture}[H,centerzero,scale=.9]
\draw[-to] (-0.2,-0.2) -- (0.2,0.2);
\draw[-to] (0.2,-0.2) -- (-0.2,0.2);
\end{tikzpicture}:E \circ E \Rightarrow E \circ E$
satisfying the dAHA relations \cref{Hdotslide,Hquadratic,Hbraid}.
\item[(H3)]
Defining the rightward crossing \begin{tikzpicture}[H,centerzero,scale=1]
\draw[to-] (0.2,-0.2) -- (-0.2,0.2);
\draw[-to] (-0.2,-0.2) -- (0.2,0.2);
\end{tikzpicture}
according to \cref{Hrightpivot},
the matrix
\begin{align*}
\begin{pmatrix}
\begin{tikzpicture}[H,anchorbase]
\draw[-to] (-0.25,-0.25) -- (0.25,0.25);
\draw[to-] (0.25,-0.25) -- (-0.25,0.25);
\end{tikzpicture} &
\begin{tikzpicture}[H,anchorbase]
\draw[-to] (-0.25,0.15) to[out=-90,in=-90,looseness=3] (0.25,0.15);
\node at (0,.2) {$\phantom.$};\node at (0,-.3) {$\phantom.$};
\end{tikzpicture}
&
\begin{tikzpicture}[H,anchorbase]
\draw[-to] (-0.25,0.15) to[out=-90,in=-90,looseness=3] (0.25,0.15);
\node at (0,.2) {$\phantom.$};\node at (0,-.3) {$\phantom.$};
\opendot{0.23,-0.03};
\end{tikzpicture}
&\!\!\cdots\!\!
&
\begin{tikzpicture}[H,anchorbase]
\draw[-to] (-0.25,0.15) to[out=-90,in=-90,looseness=3] (0.25,0.15);
\node at (0,.2) {$\phantom.$};\node at (0,-.3) {$\phantom.$};
\multopendot{0.23,-0.03}{west}{-\kappa-1};
\end{tikzpicture}
\end{pmatrix}\phantom{_T}&:E \circ F \oplus \id_\catR^{\oplus (-\kappa)}
\Rightarrow F \circ E
&&\text{if } \kappa \leq 0,\text{ or}\\
\begin{pmatrix}   
\begin{tikzpicture}[H,centerzero]
\draw[-to] (-0.25,-0.25) -- (0.25,0.25);
\draw[to-] (0.25,-0.25)  -- (-0.25,0.25);
\end{tikzpicture} &
\begin{tikzpicture}[H,centerzero]
\draw[-to] (-0.25,-0.15) to [out=90,in=90,looseness=3](0.25,-0.15);
\node at (0,.3) {$\phantom.$};
\node at (0,-.4) {$\phantom.$};
\end{tikzpicture}
&
\begin{tikzpicture}[H,centerzero]
\draw[-to] (-0.25,-0.15) to [out=90,in=90,looseness=3](0.25,-0.15);
\node at (0,.3) {$\phantom.$};
\node at (0,-.4) {$\phantom.$};
\opendot{-0.23,.03};
\end{tikzpicture}
&
\!\!\!\cdots\!\!\!
&
\begin{tikzpicture}[H,centerzero]
\draw[-to] (-0.25,-0.15) to [out=90,in=90,looseness=3](0.25,-0.15);
\node at (0,.3) {$\phantom.$};
\node at (0,-.4) {$\phantom.$};
\multopendot{-0.23,.03}{east}{\kappa-1};
\end{tikzpicture}\ 
\end{pmatrix}^\transpose\phantom{_T}&:E \circ F \Rightarrow 
F \circ E \oplus \one^{\oplus \kappa}
&&\text{if } \kappa > 0,
\end{align*}
is
invertible.
\end{itemize}
Thus, $\catR$ is a strict left $\Heis_\kappa$-module category.
\end{defin}

There is an analogy between the properties (H1)--(H3) in \cref{hcatdef}
and (KM1)--(KM3) from \cref{kmcatdef}.
However, the definition of Heisenberg categorification is simpler than that of Kac-Moody categorification---there is no underlying root datum or weight decomposition, and the inversion relation (H3) depends only on 
the fixed central charge $\kappa$ whereas 
(KM3) depends locally on the weight.
There is also no analogue of (KM0), although we will introduce a variant (H0) of that later on when we discuss the examples in \cref{s5-example,s6-example}.

%% file: s4-HtoKM.tex
\setcounter{section}{3}

%=====================================
\section{Heisenberg to Kac-Moody}\label{s4-HtoKM}
%=====================================

Now that we have explained \cref{kmcatdef,hcatdef}, 
we can recall the main construction from \cite[Sec.~4]{HKM}
giving the bridge from a degenerate Heisenberg categorification to a
Kac-Moody categorification; see also \cref{s8-quantum} where we
discuss the quantum variant.
Let $\catR$ be a degenerate Heisenberg categorification of central charge $\kappa \in \Z$.
Let $L(b)\:(b \in \B)$ be a full set of pairwise inequivalent irreducible objects of $\catR$.
For $b \in \B$, we let $m_b(x)$ (resp., $n_b(x)$) be the (monic) minimal polynomial of the endomorphism of $E L(b)$ (resp., $F L(b)$) defined by $\begin{tikzpicture}[H,centerzero]
\draw[-to] (0,-0.2) -- (0,0.2);
\opendot{0,-0.02};
\end{tikzpicture}$ (resp., by $\begin{tikzpicture}[H,centerzero]
\draw[to-] (0,-0.2) -- (0,0.2);
\opendot{0,0.03};
\end{tikzpicture}$).

Let $Z(\catR)$ be the center of the category $\catR$, that is, the
endomorphism algebra of the identity functor $\id_\catR:\catR\rightarrow\catR$.
The bubble generating function 
$\ \begin{tikzpicture}[H,baseline=-1mm]
\draw[to-] (-0.2,0) arc(180:-180:0.2);
\node at (.48,0) {$\color{black}(u)$}; 
\end{tikzpicture}$ defines an
element of $Z(\catR)\lround u^{-1} \rround$.
By Schur's Lemma, 
this evaluates to a formal Laurent series $\chi_b(u) \in \kk\lround u^{-1} \rround$ on the irreducible object $L(b)\:(b \in \B)$. The coefficients of $\chi_b(u)$ encode 
useful information about the central character of $L(b)$.
The first key result, which is \cite[Lem.~4.4]{HKM},
is that
\begin{align}\label{mainfact}
\chi_b(u) = n_b(u) / m_b(u).
\end{align}
In particular, $\chi_b(u)$ is a rational function.

The {\em spectrum} of $\catR$ is the subset $I \subseteq \kk$ consisting of the roots of the minimal polynomials $m_b(x)$ for all $b \in \B$; equivalently, by adjunction, it is the set of roots of the minimal polynomials $n_b(x)$ for all $b \in \B$.
According to \cite[Lem.~4.6]{HKM}, the set $I$ is closed under the operations $i \mapsto i \pm 1$. We define a symmetric Cartan matrix $A = (a_{i,j})_{i,j \in I}$ by setting
\begin{equation}\label{cartanmx}
a_{i,j} := \begin{cases}
2&\text{if $i = j$}\\
-1&\text{if $i = j \pm 1$ and $\ch \kk \neq 2$}\\
-2&\text{if $i = j+1 = j-1$}\\
0&\text{otherwise.}
\end{cases}
\end{equation}
The connected components of $A$ are all of type $A_\infty$ if 
$\ch\kk = 0$
or $A_{p-1}^{(1)}$ if 
$\ch\kk = p > 0$.
We also fix a root datum of this Cartan type with root lattice 
\begin{equation}\label{wtlattice}
X := \bigoplus_{i \in I} \Z \varpi_i
\end{equation}
spanned by the fundamental weights $\varpi_i\:(i \in I)$,
and simple roots defined by $\alpha_j := \sum_{i \in I} a_{i,j} \varpi_i$.
The unique homorphisms $h_i:X \rightarrow \Z$ with $h_i(\varpi_j) = \delta_{i,j}$ have the required property that
$h_i(\alpha_j) = a_{i,j}$ for all $i,j \in I$.

Since $\kk$ is algebraically closed, we can write
\begin{align}\label{epsphi}
m_b(x) &= \prod_{i \in I} (x-i)^{\eps_i(b)},
&
n_b(x) &= \prod_{i \in I} (x-i)^{\phi_i(b)}
\end{align}
for some functions $\eps_i, \phi_i: \B \rightarrow \N$.
We also define $\wt:\B \rightarrow X$ by setting
\begin{equation}\label{wtdef}
\wt(b) := \sum_{i \in I} (\phi_i(b) - \eps_i(b)) \varpi_i \in X.
\end{equation}
By \cref{mainfact}, 
the coefficient $h_i(\wt(b))$ of $\varpi_i$ in $\wt(b)$
is the multiplicity of $i$ as a zero or pole of the rational function $\chi_b(u)$. Thus, $\wt(b)$ encodes the same central character information as the formal series $\chi_b(u)$.
In view of \cref{Hbubblegeneratingfunction1}, $\chi_b(u)$ has leading term $u^{\kappa}$, so we have that
\begin{equation}\label{dunking}
\kappa = \sum_{i \in I} h_i(\wt(b))
\end{equation}
for all $b \in \B$.

For $\lambda \in X$, let $\B_\lambda := \{b \in \B\:|\:\wt(b) = \lambda\}$, and definte the {\em weight subcategory} 
$\catR_\lambda$ be the Serre subcategory of $\catR$ generated by the irreducible objects $L(b)$ for $b \in \B_\lambda$.
By central character considerations, the category $\catR$ decomposes as
the internal direct sum
\begin{equation}\label{blockdec}
\catR = \bigoplus_{\lambda \in X} \catR_\lambda.
\end{equation}
By \cref{dunking}, the weight subcategory $\catR_\lambda$ is zero unless $\sum_{i \in I} h_i(\lambda) = \kappa$.

The functors $E$ and $F$ decompose as
\begin{align}\label{eigenfunctors}
E &= \bigoplus_{i \in I} E_i,&
F &= \bigoplus_{i \in I} F_i
\end{align}
where $E_i$ is the subfunctor of $E$ such that $E_i V$ is the generalized $i$-eigenspace of the endomorphism of $EV$ defined by
$\begin{tikzpicture}[H,centerzero]
\draw[-to] (0,-0.2) -- (0,0.2);
\opendot{0,-0.02};
\end{tikzpicture}$
for each $V \in \ob\catR$,
and $F_i$ is the subfunctor of $F$ defined similarly by taking the generalized $i$-eigenspace of
$\begin{tikzpicture}[H,centerzero]
\draw[to-] (0,-0.2) -- (0,0.2);
\opendot{0,0.03};
\end{tikzpicture}$.
The definition of the spectrum $I$ ensures that all of the functors $E_i$ and $F_i$ appearing in
\cref{eigenfunctors} are non-zero.

To complete the picture, we are going to use the usual string calculus for the 2-category of $\kk$-linear categories to introduce some further natural transformations between these functors. In such diagrams, an unlabelled 2-cell should be interpreted as the category $\catR$, and a 2-cell labelled by $\lambda$ indicates the category $\catR_\lambda$.
The identity endomorphisms of $E:\catR\rightarrow \catR$ and $F:\catR\rightarrow\catR$ will be 
represented by the strings
$\ \begin{tikzpicture}[H,anchorbase]
\draw[-to] (0,-0.2) -- (0,0.2);
\end{tikzpicture}\ $ and
$\ \begin{tikzpicture}[H,anchorbase]
\draw[to-] (0,-0.2) -- (0,0.2);
\end{tikzpicture}\ $.
The same strings with the rightmost 
2-cell labelled $\lambda$ 
denote the identity endomorphisms of the functors obtained by pre-composing 
with the inclusion $\operatorname{Inc}_\lambda:\catR_\lambda \rightarrow \catR$, while if the 
the leftmost 2-cell is labelled $\mu$
we mean the identity endomorphisms of the functors obtained by 
post-composing with the projection $\operatorname{Pr}_\mu:\catR\rightarrow \catR_\mu$.
Similar conventions apply to $E_i:\catR\rightarrow\catR$ and $F_i:\catR\rightarrow\catR$, whose identity endomorphisms
will be represented using the labelled strings
$\begin{tikzpicture}[KM,centerzero]
\draw[-to] (0,-0.2) \botlabel{i} -- (0,0.2);
\end{tikzpicture}$
and
$\begin{tikzpicture}[KM,centerzero]
\draw[to-] (0,-0.2) \botlabel{i} -- (0,0.2);
\end{tikzpicture}$, respectively.
\iffalse
Note by \cref{firstaxiom} that
$\begin{tikzpicture}[KM,centerzero]
\draw[-to] (0,-0.2) \botlabel{i} -- (0,0.2);
\region{0.3,0}{\lambda};
\region{-0.3,0}{\mu};
\end{tikzpicture} = \begin{tikzpicture}[KM,centerzero]
\draw[-to] (0,-0.2) \botlabel{i} -- (0,0.2);
\region{-0.3,0}{\lambda};
\region{0.3,0}{\mu};
\end{tikzpicture} =
0$ if $\mu \neq \lambda+\alpha_i$.
Hence, by \cref{eigenfunctors}, we have that
$\begin{tikzpicture}[H,centerzero]
\draw[-to] (0,-0.2) -- (0,0.2);
\region{0.3,0}{\lambda};
\region{-0.3,0}{\mu};
\end{tikzpicture} = \begin{tikzpicture}[H,centerzero]
\draw[-to] (0,-0.2) -- (0,0.2);
\region{-0.3,0}{\lambda};
\region{0.3,0}{\mu};
\end{tikzpicture} =
0$ if $\mu \neq \lambda+\alpha_i$ for some $i \in I$.
\fi
We use the 
string diagrams
\begin{align}\label{projectors}
\begin{tikzpicture}[KM,centerzero]
\draw (0.08,-.3) \botlabel{i} to (.08,0);
\draw[H,-to] (.08,0) to (0.08,.3);
\notch{.08,0};
\end{tikzpicture}
:E_i &\Rightarrow E,
&
\begin{tikzpicture}[KM,centerzero]
\draw[to-] (0.08,-.3) \botlabel{i} to (0.08,0);
\draw[H] (.08,0) to (0.08,.3);
\notch{.08,0};
\end{tikzpicture}
:F_i &\Rightarrow F,
&
\begin{tikzpicture}[KM,centerzero]
\draw[H] (0.08,-.3) to (0.08,0);
\draw[-to] (.08,0) to (0.08,.3)\toplabel{i};
\notch{.08,0};
\end{tikzpicture}
:E &\Rightarrow E_i,
&
\begin{tikzpicture}[KM,centerzero]
\draw[to-,H] (0.08,-.3) to (0.08,0);
\draw (.08,0) to (.08,.3)\toplabel{i};
\notch{.08,0};
\end{tikzpicture}
:F &\Rightarrow F_i
\end{align}
to denote the various inclusions and projections between $E$ and $F$ and their summands $E_i$ and $F_i$.
We then have that
\begin{align}\label{orthogonality}
\begin{tikzpicture}[KM,centerzero]
\draw[-to] (0.08,-.3) \botlabel{i} to (0.08,.3)\toplabel{j};
\draw[H] (0.08,-.1) to (0.08,.1);
\notch{.08,.1};
\notch{.08,-.1};
\end{tikzpicture}
&= \delta_{i,j}
\begin{tikzpicture}[KM,centerzero]
\draw[-to] (0.08,-.3) \botlabel{i} to (0.08,.3);
\end{tikzpicture}\ ,&
\begin{tikzpicture}[KM,centerzero]
\draw[to-] (0.08,-.3) \botlabel{i} to (0.08,.3)\toplabel{j};
\draw[H] (0.08,-.07) to (0.08,.13);
\notch{.08,.1};
\notch{.08,-.1};
\end{tikzpicture}
&= \delta_{i,j}
\begin{tikzpicture}[KM,centerzero]
\draw[to-] (0.08,-.3) \botlabel{i} to (0.08,.3);
\end{tikzpicture}\ .
\end{align}

Composing the rightward cup and cap with projectors 
gives natural transformations
\begin{align}\label{flimsy}
\begin{tikzpicture}[KM,anchorbase,scale=1.1]
\draw[to-] (0.3,0.4)\toplabel{j} to (0.3,.2);
\draw[-] (-0.1,.2) to (-0.1,0.4)\toplabel{i};
\draw[H] (.3,.2) to[out=-90,in=-90,looseness=1.5] (-.1,.2);
\notch{-.1,.2};
\notch{.3,.2};
\end{tikzpicture}
:&\id_\catR \Rightarrow F_i \circ E_j,&
\begin{tikzpicture}[KM,anchorbase,scale=1.1]
\draw[to-] (0.3,-0.4)\botlabel{j} to (0.3,-.2);
\draw[-] (-0.1,-.2) to (-0.1,-0.4)\botlabel{i};
\draw[H] (.3,-.2) to[out=90,in=90,looseness=1.5] (-.1,-.2);
\notch{-.1,-.2};
\notch{.3,-.2};
\end{tikzpicture}:&E_i \circ F_j \Rightarrow \id_\catR.
\end{align}
These are zero if $i \neq j$, as may be proved by a similar argument to the proof of \cite[Lem.~4.1]{HKM} using that dots slide across rightward cups and caps in $\Heis_\kappa$. It follows that
\begin{align}\label{icecreamtruck}
\begin{tikzpicture}[KM,centerzero,scale=1.3]
\draw[-,H] (-0.3,.1) arc(180:0:0.15);
\draw[-,H] (0,-.1) arc(180:360:0.15);
\draw (0,-.1) to (0,.1);
\draw (-.3,0.1) to (-.3,-.4)\botlabel{i};
\draw[-to] (.3,-0.1) to (.3,.4)\toplabel{i};
\notch{0,-.1};
\notch{.3,-.1};
\notch{-.3,.1};
\notch{0,.1};
\stringlabel{0.11,0.02}{i};
\end{tikzpicture}
&=
\sum_{j \in I} \begin{tikzpicture}[KM,centerzero,scale=1.3]
\draw[-,H] (-0.3,.1) arc(180:0:0.15);
\draw[-,H] (0,-.1) arc(180:360:0.15);
\draw (0,-.1) to (0,.1);
\draw (-.3,0.1) to (-.3,-.4)\botlabel{i};
\draw[-to] (.3,-0.1) to (.3,.4)\toplabel{i};
\notch{0,-.1};
\notch{0,.1};
\notch{.3,-.1};
\notch{-.3,.1};
\stringlabel{0.11,0.02}{j};
\end{tikzpicture}
=
\begin{tikzpicture}[KM,centerzero,scale=1.3]
\draw[-,H] (-0.3,.1) arc(180:0:0.15);
\draw[-,H] (0,-.1) arc(180:360:0.15);
\draw[H] (0,-.1) to (0,.1);
\draw (-.3,0.1) to (-.3,-.4)\botlabel{i};
\draw[-to] (.3,-0.1) to (.3,.4)\toplabel{i};
\notch{.3,-.1};
\notch{-.3,.1};
\end{tikzpicture}
\stackrel{\cref{Hrightadj}}{=}\begin{tikzpicture}[KM,centerzero,scale=1.3]
\draw[-to] (0,-0.4)\botlabel{i} -- (0,0.4)\toplabel{i};
\draw[H] (0,.15) to (0,-.15);
\notch{0,.15};
\notch{0,-.15};
\end{tikzpicture}
\stackrel{\cref{orthogonality}}{=}
\begin{tikzpicture}[KM,centerzero,scale=1.3]
\draw[-to] (0,-0.4)\botlabel{i} -- (0,0.4);
\end{tikzpicture}
\ ,&
\begin{tikzpicture}[KM,centerzero,scale=1.3]
\draw[-,H] (-0.3,-.1) arc(180:360:0.15);
\draw[-,H] (0,.1) arc(180:0:0.15);
\draw (0,-.1) to (0,.1);
\draw (-.3,-0.1) to (-.3,.4)\toplabel{i};
\draw[-to] (.3,0.1) to (.3,-.4)\botlabel{i};
\notch{0,.1};
\notch{0,-.1};
\notch{.3,.1};
\notch{-0.3,-.1};
\stringlabel{0.11,0.02}{i};
\end{tikzpicture}
&=
\sum_{j \in I} \begin{tikzpicture}[KM,centerzero,scale=1.3]
\draw[-,H] (-0.3,-.1) arc(180:360:0.15);
\draw[-,H] (0,.1) arc(180:0:0.15);
\draw (0,-.1) to (0,.1);
\draw (-.3,-0.1) to (-.3,.4)\toplabel{i};
\draw[-to] (.3,0.1) to (.3,-.4)\botlabel{i};
\notch{0,.1};
\notch{0,-.1};
\notch{.3,.1};
\notch{-.3,-.1};
\stringlabel{0.11,0.02}{j};
\end{tikzpicture}
=
\begin{tikzpicture}[KM,centerzero,scale=1.3]
\draw[-,H] (-0.3,-.1) arc(180:360:0.15);
\draw[-,H] (0,.1) arc(180:0:0.15);
\draw[H] (0,-.1) to (0,.1);
\draw (-.3,-0.1) to (-.3,.4)\toplabel{i};
\draw[-to] (.3,0.1) to (.3,-.4)\botlabel{i};
\notch{.3,.1};
\notch{-.3,-.1};
\end{tikzpicture}
\stackrel{\cref{Hrightadj}}{=}
\begin{tikzpicture}[KM,centerzero,scale=1.3]
\draw[to-] (0,-0.4)\botlabel{i} -- (0,0.4)\toplabel{i};
\draw[H] (0,.15) to (0,-.15);
\notch{0,.15};
\notch{0,-.15};
\end{tikzpicture}
\stackrel{\cref{orthogonality}}{=}
\begin{tikzpicture}[KM,centerzero,scale=1.3]
\draw[to-] (0,-0.4)\botlabel{i} -- (0,0.4);
\end{tikzpicture}
\ .
\end{align}
This shows that the {\em rightward Kac-Moody cups and caps}
defined by
\begin{align}
\begin{tikzpicture}[KM,centerzero]
\draw[-to] (-0.25,0.15) \toplabel{i} to[out=-90,in=-90,looseness=2.5] (0.25,0.15);
\end{tikzpicture}\ &:=\begin{tikzpicture}[KM,anchorbase,scale=1.1]
\draw[to-] (0.3,0.4)\toplabel{i} to (0.3,.2);
\draw[-] (-0.1,.2) to (-0.1,0.4)\toplabel{i};
\draw[H] (.3,.2) to[out=-90,in=-90,looseness=1.5] (-.1,.2);
\notch{-.1,.2};
\notch{.3,.2};
\end{tikzpicture}\ ,&
\begin{tikzpicture}[KM,centerzero]
\draw[-to] (-0.25,-0.15) \botlabel{i} to [out=90,in=90,looseness=2.5](0.25,-0.15);
\end{tikzpicture}\ 
&:=
\begin{tikzpicture}[KM,anchorbase,scale=1.1]
\draw[to-] (0.3,-0.4)\botlabel{i} to (0.3,-.2);
\draw[-] (-0.1,-.2) to (-0.1,-0.4)\botlabel{i};
\draw[H] (.3,-.2) to[out=90,in=90,looseness=1.5] (-.1,-.2);
\notch{-.1,-.2};
\notch{.3,-.2};
\end{tikzpicture}\label{gofigure}\end{align}
are the unit and counit of an adjunction making
$(E_i, F_i)$ into an adjoint pair for each $i \in I$.
Using the vanishing properties of \cref{flimsy}, 
it is also easy to check that
\begin{align}\label{sliders}
\begin{tikzpicture}[KM,anchorbase,scale=1.1]
\draw[-] (-0.1,.2) to (-0.1,0.4)\toplabel{i};
\draw[to-,H] (.3,.4) to (.3,.2) to[out=-90,in=-90,looseness=1.5] (-.1,.2);
\notch{-.1,.2};
\end{tikzpicture}\ 
&=
\ \begin{tikzpicture}[KM,anchorbase,scale=1.1]
\draw[to-,H] (0.3,.4) to (0.3,0.2);
\draw[-] (.3,.2) to[out=-90,in=-90,looseness=1.5] (-.1,.2) to (-.1,.4)\toplabel{i};
\notch{.3,.2};
\end{tikzpicture}\ ,&
\begin{tikzpicture}[KM,anchorbase,scale=1.1]
\draw[-,H] (-0.1,.2) to (-0.1,0.4);
\draw[to-] (.3,.4) \toplabel{i}to (.3,.2) to[out=-90,in=-90,looseness=1.5] (-.1,.2);
\notch{-0.1,.2};
\end{tikzpicture}\ 
&=
\ \begin{tikzpicture}[KM,anchorbase,scale=1.1]
\draw[to-] (0.3,.4) \toplabel{i}to (0.3,0.2);
\draw[-,H] (.3,.2) to[out=-90,in=-90,looseness=1.5] (-.1,.2) to (-.1,.4);
\notch{.3,.2};
\end{tikzpicture}\ ,&
\begin{tikzpicture}[KM,anchorbase,scale=1.1]
\draw[-] (-0.1,-.2) to (-0.1,-0.4)\botlabel{i};
\draw[to-,H] (.3,-.4) to (.3,-.2) to[out=90,in=90,looseness=1.5] (-.1,-.2);
\notch{-.1,-.2};
\end{tikzpicture}\ 
&=
\ \begin{tikzpicture}[KM,anchorbase,scale=1.1]
\draw[to-,H] (0.3,-.4) to (0.3,-0.2);
\draw[-] (.3,-.2) to[out=90,in=90,looseness=1.5] (-.1,-.2) to (-.1,-.4)\botlabel{i};
\notch{.3,-.2};
\end{tikzpicture}\ ,&
\begin{tikzpicture}[KM,anchorbase,scale=1.1]
\draw[-,H] (-0.1,-.2) to (-0.1,-0.4);
\draw[to-] (.3,-.4) \botlabel{i}to (.3,-.2) to[out=90,in=90,looseness=1.5] (-.1,-.2);
\notch{-.1,-.2};
\end{tikzpicture}\ 
&=
\ \begin{tikzpicture}[KM,anchorbase,scale=1.1]
\draw[to-] (0.3,-.4) \botlabel{i}to (0.3,-0.2);
\draw[-,H] (.3,-.2) to[out=90,in=90,looseness=1.5] (-.1,-.2) to (-.1,-.4);
\notch{.3,-.2};
\end{tikzpicture}\ .
\end{align}
In a similar way, the leftward cup and cap in $\Heis_\kappa$ can be projected to obtain the unit and counit of an adjunction $(F_i, E_i)$, but we will not introduce any special diagrammatic notation for these natural transformations---they are {\em not} simply equal to the leftward Kac-Moody cups and caps
$\begin{tikzpicture}[KM,centerzero,scale=.8]
\draw[to-] (-0.25,0.15) \toplabel{i} to[out=-90,in=-90,looseness=2.5] (0.25,0.15);
\end{tikzpicture}$
and 
$\begin{tikzpicture}[KM,centerzero,scale=.8]
\draw[to-] (-0.25,-0.15) \botlabel{i} to [out=90,in=90,looseness=2.5](0.25,-0.15);
\end{tikzpicture}$ appearing in \cref{average} below.

For the projected dots, we incorporate an additional change-of-origin, defining the {\em Kac-Moody dots}
\begin{align}\label{prodots}
\begin{tikzpicture}[KM,centerzero]
\draw[-to] (0.08,-.35) \botlabel{i} to (.08,.35);
\opendot{0.08,0};
\end{tikzpicture}
&:=
\begin{tikzpicture}[H,centerzero]
\draw[KM,-to] (0.08,-.35) \botlabel{i} to (.08,.35);
\draw (.08,-.15) to (0.08,.15);
\notch{0.08,.15};
\notch{0.08,-.15};
\pin{0.08,0}{.7,0}{x-i};
\end{tikzpicture}\ ,&
\begin{tikzpicture}[KM,centerzero]
\draw[to-] (0.08,-.35) \botlabel{i} to (.08,.35);
\opendot{0.08,0};
\end{tikzpicture}
&:=
\begin{tikzpicture}[H,centerzero]
\draw[KM,to-] (0.08,-.35) \botlabel{i} to (.08,.35);
\draw (.08,-.15) to (0.08,.15);
\notch{0.08,.15};
\notch{0.08,-.15};
\pin{0.08,0}{.7,0}{x-i};
\end{tikzpicture}\ .
\end{align}
The point of this shift is that these natural transformations are {\em locally nilpotent}, that is, they evaluate to nilpotent endomorphisms on any object of $\catR$. Consequently, we can extend the pin notation by allowing any formal power series in
$\kk\llbracket x \rrbracket$ (rather than merely a polynomial in $\kk[x]$) to be pinned to the Kac-Moody dots
$\begin{tikzpicture}[KM,centerzero]
\draw[-to] (0.08,-.2) \botlabel{i} to (.08,.2);
\opendot{0.08,0};
\end{tikzpicture}$
or 
$\begin{tikzpicture}[KM,centerzero]
\draw[to-] (0.08,-.2) \botlabel{i} to (.08,.2);
\opendot{0.08,0};
\end{tikzpicture}$.
From \cref{orthogonality}, it follows that
\begin{align}\label{pan}
\begin{tikzpicture}[KM,centerzero]
\draw (0.08,-.3) \botlabel{i} to (.08,0.1);
\draw[H,-to] (.08,0.1) to (0.08,.3);
\notch{.08,.1};
\opendot{.08,-.1};
\end{tikzpicture}
&=
\begin{tikzpicture}[H,centerzero]
\draw[KM] (0.08,-.3) \botlabel{i} to (.08,-0.1);
\draw[-to] (.08,-0.1) to (0.08,.3);
\notch{.08,-.1};
\pin{.08,.1}{.7,.1}{x-i};
\end{tikzpicture}\ ,&
\begin{tikzpicture}[KM,centerzero]
\draw[to-] (0.08,-.3) \botlabel{i} to (.08,0.1);
\draw[H] (.08,0.1) to (0.08,.3);
\notch{0.08,.1};
\opendot{.08,-.1};
\end{tikzpicture}
&=
\begin{tikzpicture}[H,centerzero]
\draw[KM,to-] (0.08,-.3) \botlabel{i} to (.08,-0.1);
\draw (.08,-0.1) to (0.08,.3);
\notch{0.08,-.1};
\pin{.08,.1}{.7,.1}{x-i};
\end{tikzpicture}\ ,&
\begin{tikzpicture}[H,centerzero]
\draw (0.08,-.3)  to (.08,0.1);
\draw[KM,-to] (.08,0.1) to (0.08,.3) \toplabel{i};
\notch{.08,.1};
\pin{.08,-.1}{.7,-.1}{x-i};
\end{tikzpicture}
&=
\begin{tikzpicture}[KM,centerzero]
\draw[H] (0.08,-.3)  to (.08,-0.1);
\draw[-to] (.08,-0.1) to (0.08,.3)\toplabel{i};
\notch{0.08,-.1};
\opendot{.08,.1};
\end{tikzpicture}\ ,&
\begin{tikzpicture}[H,centerzero]
\draw[to-](0.08,-.3)  to (.08,0.1);
\draw[KM] (.08,0.1) to (0.08,.3) \toplabel{i};
\notch{.08,.1};
\pin{.08,-.1}{.7,-.1}{x-i};
\end{tikzpicture}
&=
\begin{tikzpicture}[KM,centerzero]
\draw[H,to-] (0.08,-.3)  to (.08,-0.1);
\draw (.08,-0.1) to (0.08,.3)\toplabel{i};
\notch{.08,-.1};
\opendot{.08,.1};
\end{tikzpicture}\ .
\end{align}
Using these and the fact that Heisenberg dots slide across Heisenberg cups and caps, it is easy to deduce that the Kac-Moody dots slide across the rightward Kac-Moody cups and caps.

The analysis of projected crossings is more interesting. We will denote them by
\begin{align}\label{polly}
\begin{tikzpicture}[centerzero,KM]
\draw[-to] (0.28,-.28) \botlabel{j} to (-0.28,.28)\toplabel{i'};
\draw[-to] (-0.28,-.28) \botlabel{i} to (0.28,.28)\toplabel{j'};
\projcr{0,0};
\end{tikzpicture}
&:=
\begin{tikzpicture}[centerzero,KM]
\draw[-to] (0.28,-.28) \botlabel{j} to (-0.28,.28)\toplabel{i'};
\draw[-to] (-0.28,-.28) \botlabel{i} to (0.28,.28)\toplabel{j'};
\draw[H] (0.1,-.1) to (-0.1,.1);
\draw[H] (-0.1,-.1) to (0.1,.1);
\notch[45]{.1,-.1};
\notch[45]{-.1,.1};
\notch[-45]{.1,.1};
\notch[-45]{-.1,-.1};
\end{tikzpicture}\ ,&
\begin{tikzpicture}[centerzero,KM]
\draw[to-] (0.28,-.28) \botlabel{j'} to (-0.28,.28)\toplabel{i};
\draw[-to] (-0.28,-.28) \botlabel{j} to (0.28,.28)\toplabel{i'};
\projcr{0,0};
\end{tikzpicture}
&:=
\begin{tikzpicture}[centerzero,KM]
\draw[to-] (0.28,-.28) \botlabel{j'} to (-0.28,.28)\toplabel{i};
\draw[-to] (-0.28,-.28) \botlabel{j} to (0.28,.28)\toplabel{i'};
\draw[H] (0.1,-.1) to (-0.1,.1);
\draw[H] (-0.1,-.1) to (0.1,.1);
\notch[45]{.1,-.1};
\notch[45]{-.1,.1};
\notch[-45]{.1,.1};
\notch[-45]{-.1,-.1};
\end{tikzpicture}\ ,
&\begin{tikzpicture}[centerzero,KM]
\draw[to-] (0.28,-.28) \botlabel{i'} to (-0.28,.28)\toplabel{j};
\draw[to-] (-0.28,-.28) \botlabel{j'} to (0.28,.28)\toplabel{i};
\projcr{0,0};
\end{tikzpicture}
&:=
\begin{tikzpicture}[centerzero,KM]
\draw[to-] (0.28,-.28) \botlabel{i'} to (-0.28,.28)\toplabel{j};
\draw[to-] (-0.28,-.28) \botlabel{j'} to (0.28,.28)\toplabel{i};
\draw[H] (0.1,-.1) to (-0.1,.1);
\draw[H] (-0.1,-.1) to (0.1,.1);
\notch[45]{.1,-.1};
\notch[45]{-.1,.1};
\notch[-45]{.1,.1};
\notch[-45]{-.1,-.1};
\end{tikzpicture}\ ,&
\begin{tikzpicture}[centerzero,KM]
\draw[-to] (0.28,-.28) \botlabel{i} to (-0.28,.28)\toplabel{j'};
\draw[to-] (-0.28,-.28) \botlabel{i'} to (0.28,.28)\toplabel{j};
\projcr{0,0};
\end{tikzpicture}
&:=
\begin{tikzpicture}[centerzero,KM]
\draw[-to] (0.28,-.28) \botlabel{i} to (-0.28,.28)\toplabel{j'};
\draw[to-] (-0.28,-.28) \botlabel{i'} to (0.28,.28)\toplabel{j};
\draw[H] (0.1,-.1) to (-0.1,.1);
\draw[H] (-0.1,-.1) to (0.1,.1);
\notch[45]{.1,-.1};
\notch[45]{-.1,.1};
\notch[-45]{.1,.1};
\notch[-45]{-.1,-.1};
\end{tikzpicture}\ .
\end{align}
In \cite[Lem.~4.1]{HKM}, it is shown that these natural transformations are zero unless $\{i,j\} = \{i',j'\}$. Also, by \cite[Lem.~4.2]{HKM},
we have that
\begin{align}
\begin{tikzpicture}[centerzero,KM]
\draw[-to] (0.28,-.28) \botlabel{j} to (-0.28,.28)\toplabel{i};
\draw[-to] (-0.28,-.28) \botlabel{i} to (0.28,.28)\toplabel{j};
\projcr{0,0};
\end{tikzpicture}
&=
\begin{tikzpicture}[centerzero,KM]
\draw[-to] (0.2,-.28) \botlabel{j} to (0.2,.28);
\draw[-to] (-0.2,-.28) \botlabel{i} to (-0.2,.28);
\pinpin{-.2,0}{.2,0}{1.5,0}{(i-j+x-y)^{-1}};
\end{tikzpicture}\label{eggs}\ ,&
\begin{tikzpicture}[centerzero,KM]
\draw[to-] (0.28,-.28) \botlabel{j} to (-0.28,.28)\toplabel{i};
\draw[to-] (-0.28,-.28) \botlabel{i} to (0.28,.28)\toplabel{j};
\projcr{0,0};
\end{tikzpicture}
&=
-\:\begin{tikzpicture}[centerzero,KM]
\draw[to-] (0.2,-.28) \botlabel{j} to (0.2,.28);
\draw[to-] (-0.2,-.28) \botlabel{i} to (-0.2,.28);
\pinpin{-.2,0}{.2,0}{1.5,0}{(i-j+x-y)^{-1}};
\end{tikzpicture}
\end{align}
when $i \neq j$; the expression on the right hand side  makes sense because $i-j+x-y \in \kk\llbracket x, y \rrbracket^\times$ and $x, y$ are locally nilpotent.
In the following lemma, we record some
further properties of the projected crossings. The proof makes a good exercise to get better acquainted with the diagrammatics.
(Perhaps this is a good time to remind again of our labelling convention for multipins: alphabetical order of variables corresponds to lexicographic order of Cartesian coordinates.)

\begin{lem}\label{muckypups}
Recall that $R(x,y) = 1-(x-y)^{-2}$,
and $\partial_{xy}$ denotes the Demazure operator from \cref{muckycats}.
The following hold for $i,j,k \in I$:
\begin{align}\label{muckypups1}
\begin{tikzpicture}[centerzero,KM]
\draw[-to] (0.28,-.28) \botlabel{j} to (-0.28,.28)\toplabel{j};
\draw[-to] (-0.28,-.28) \botlabel{i} to (0.28,.28)\toplabel{i};
\opendot{-.18,-.18};
\projcr{0,0};
\end{tikzpicture}
-\begin{tikzpicture}[centerzero,KM]
\draw[-to] (0.28,-.28) \botlabel{j} to (-0.28,.28)\toplabel{j};
\draw[-to] (-0.28,-.28) \botlabel{i} to (0.28,.28)\toplabel{i};
\opendot{.15,.15};
\projcr{0,0};
\end{tikzpicture}
&=
\delta_{i,j}\begin{tikzpicture}[centerzero,KM]
\draw[-to] (0.15,-.28) \botlabel{i} to (0.15,.28);
\draw[-to] (-0.15,-.28) \botlabel{i} to (-0.15,.28);
\end{tikzpicture}=
\begin{tikzpicture}[centerzero,KM]
\draw[-to] (0.28,-.28) \botlabel{j} to (-0.28,.28)\toplabel{j};
\draw[-to] (-0.28,-.28) \botlabel{i} to (0.28,.28)\toplabel{i};
\opendot{-.15,.15};
\projcr{0,0};
\end{tikzpicture}
-
\begin{tikzpicture}[centerzero,KM]
\draw[-to] (0.28,-.28) \botlabel{j} to (-0.28,.28)\toplabel{j};
\draw[-to] (-0.28,-.28) \botlabel{i} to (0.28,.28)\toplabel{i};
\opendot{.18,-.18};
\projcr{0,0};
\end{tikzpicture}\ ,\\
\label{cold}
\begin{tikzpicture}[KM,centerzero]
\draw[-to] (-0.2,-0.6) \botlabel{i} \braidup (0.2,0) \braidup (-0.2,0.6) \toplabel{i};
\draw[-to] (0.2,-0.6) \botlabel{j} \braidup (-0.2,0) \braidup (0.2,0.6) \toplabel{j};
\projcr{0,-0.3};
\projcr{0,0.3};
\strand{0.35,0}{i};
\strand{-0.35,0}{j};
\end{tikzpicture}
&=
\begin{dcases}
\begin{tikzpicture}[KM,centerzero]
\draw[-to] (-0.1,-0.3) \botlabel{i} -- (-0.1,0.3);
\draw[-to] (0.2,-0.3) \botlabel{j} -- (0.2,0.3);
\pinpin{-0.1,0}{0.2,0}{1.3,0}{R(x+i,y+j)};
\end{tikzpicture}
&\text{if $i \ne j$}\\
\ \begin{tikzpicture}[KM,centerzero]
\draw[-to] (-0.1,-0.3) \botlabel{i} -- (-0.1,0.3);
\draw[-to] (0.2,-0.3) \botlabel{i} -- (0.2,0.3);
\end{tikzpicture}
&\text{if $i=j$,}
\end{dcases}\\
\begin{tikzpicture}[KM,centerzero,scale=.8]
\draw[-to] (-0.5,-1) \botlabel{i} \braidup (0.5,1) \toplabel{i};
\draw[-to] (0,-1) \botlabel{j} \braidup (-0.5,0) \braidup (0,1) \toplabel{j};
\draw[-to] (0.5,-1) \botlabel{k} \braidup (-0.5,1) \toplabel{k};
\projcr{-0.32,0.42};
\projcr{0,0};
\projcr{-0.32,-0.42};
\node at (-0.07,-0.28) {\strandlabel{i}};
\node at (-0.05,0.31) {\strandlabel{k}};
\node at (-0.63,0) {\strandlabel{j}};
\end{tikzpicture}
-
\begin{tikzpicture}[KM,centerzero,scale=.8]
\draw[-to] (-0.5,-1) \botlabel{i} \braidup (0.5,1) \toplabel{i};
\draw[-to] (0,-1) \botlabel{j} \braidup (0.5,0) \braidup (0,1) \toplabel{j};
\draw[-to] (0.5,-1) \botlabel{k} \braidup (-0.5,1) \toplabel{k};
\projcr{0.33,0.42};
\projcr{0,0};
\projcr{0.33,-0.42};
\node at (0.07,-0.29) {\strandlabel{k}};
\node at (0.05,0.3) {\strandlabel{i}};
\node at (0.63,0) {\strandlabel{j}};
\end{tikzpicture}
&= 
\begin{dcases}
\begin{tikzpicture}[KM,centerzero,scale=.8]
\draw[-to] (-0.4,-.5) \botlabel{i} -- (-0.4,.5);
\draw[-to] (0,-.5) \botlabel{j} -- (0,.5);
\draw[-to] (0.4,-.5) \botlabel{i} -- (0.4,.5);
\pinpinpin{-0.4,0}{0,0}{0.4,0}{2,0}{\partial_{xz} R(x+i,y+j)};
\end{tikzpicture}
&\text{if $i=k \neq j$}\\
0&\text{otherwise.}
\end{dcases}\label{hot}
\end{align}
\end{lem}

Now we introduce the {\em upward Kac-Moody crossings}
\begin{equation}\label{upwardmagic}
\begin{tikzpicture}[KM,centerzero,scale=.9]
\draw[-to] (-0.4,-0.4) \botlabel{i} -- (0.4,0.4);
\draw[-to] (0.4,-0.4) \botlabel{j} -- (-0.4,0.4);
\end{tikzpicture}
:=
\begin{dcases}
\begin{tikzpicture}[KM,centerzero,scale=.9]
\draw[-to] (0.4,-.4)\botlabel{i} to (-0.4,.4)\toplabel{i};
\draw[-to] (-0.4,-.4)\botlabel{i} to (0.4,.4)\toplabel{i};
\projcr{0,0};
\pinpin{-.25,-.25}{.25,-.25}{1.5,-.25}{(1+x-y)^{-1}};
\end{tikzpicture}
+\begin{tikzpicture}[KM,centerzero,scale=.9]
\draw[-to] (0.5,-.4)\botlabel{i} to (0.5,.4);
\draw[-to] (0.1,-.4)\botlabel{i} to (0.1,.4);
\pinpin{.1,0}{.5,0}{1.7,0}{(1+x-y)^{-1}};
\end{tikzpicture}
&\text{if $i=j$}\\
\begin{tikzpicture}[KM,centerzero,scale=.9]
\draw[-to] (0.4,-.4)\botlabel{i-1} to (-0.4,.4)\toplabel{i-1};
\draw[-to] (-0.4,-.4)\botlabel{i} to (0.4,.4)\toplabel{i};
\projcr{0,0};
\pinpin{-.25,-.25}{.25,-.25}{1.25,-.25}{1+x-y};
\end{tikzpicture}
&\text{if $i=j+1$}\\
-\ \begin{tikzpicture}[KM,centerzero,scale=.9]
\draw[-to] (0.4,-.4) \botlabel{j} to (-0.4,.4)\toplabel{j};
\draw[-to] (-0.4,-.4) \botlabel{i} to (0.4,.4)\toplabel{i};
\projcr{0,0};
\pinpin{-.25,-.25}{.25,-.25}{2.85,-.25}{(i-j+x-y)(i-j-1+x-y)^{-1}};
\end{tikzpicture}
&\text{if $i \neq j,j+1$.}
\end{dcases}
\end{equation}
By \cref{muckypups1} and \cref{muckycats}, this definition can be written equivalently as
\begin{equation}\label{upwardmagicalt}
\begin{tikzpicture}[KM,centerzero,scale=.9]
\draw[-to] (-0.4,-0.4) \botlabel{i} -- (0.4,0.4);
\draw[-to] (0.4,-0.4) \botlabel{j} -- (-0.4,0.4);
\end{tikzpicture}
=
\begin{dcases}
\begin{tikzpicture}[KM,centerzero,scale=.9]
\draw[-to] (0.4,-.4)\botlabel{i} to (-0.4,.4)\toplabel{i};
\draw[-to] (-0.4,-.4)\botlabel{i} to (0.4,.4)\toplabel{i};
\projcr{0,0};
\pinpin{.25,.25}{.25,-.25}{1.5,-.25}{(1+y-x)^{-1}};
\end{tikzpicture}
&\text{if $i=j$}\\
\begin{tikzpicture}[KM,centerzero,scale=.9]
\draw[-to] (0.4,-.4)\botlabel{i-1} to (-0.4,.4)\toplabel{i-1};
\draw[-to] (-0.4,-.4)\botlabel{i} to (0.4,.4)\toplabel{i};
\projcr{0,0};
\pinpin{.25,.25}{.25,-.25}{1.25,-.25}{1+y-x};
\end{tikzpicture}
&\text{if $i=j+1$}\\
- \begin{tikzpicture}[KM,centerzero,scale=.9]
\draw[-to] (0.4,-.4) \botlabel{j} to (-0.4,.4)\toplabel{j};
\draw[-to] (-0.4,-.4) \botlabel{i} to (0.4,.4)\toplabel{i};
\projcr{0,0};
\pinpin{.25,.25}{.25,-.25}{2.85,-.25}{(i-j+y-x)(i-j-1+y-x)^{-1}};
\end{tikzpicture}
&\text{if $i \neq j,j+1$.}
\end{dcases}
\end{equation}
From the $i=j$ case of this and  \cref{lotsmore,ruby,wax}, we also have that
\begin{equation}
\begin{tikzpicture}[KM,centerzero,scale=.9]
\draw[to-] (0.4,-.4)\botlabel{i} to (-0.4,.4)\toplabel{i};
\draw[to-] (-0.4,-.4)\botlabel{i} to (0.4,.4)\toplabel{i};
\projcr{0,0};
\end{tikzpicture}=
\begin{tikzpicture}[KM,centerzero,scale=.9]
\draw[to-] (-0.4,-0.4) \botlabel{i} -- (0.4,0.4);
\draw[to-] (0.4,-0.4) \botlabel{i} -- (-0.4,0.4);
\pinpin{-.25,-.25}{-.25,.25}{-1.5,.25}{1+x-y};
\end{tikzpicture}\ .\label{billy}
\end{equation}

\begin{theo}[{\cite[Th.~4.11]{HKM}}]\label{maintheorem}
Let $\catR$ be a degenerate Heisenberg categorification of central charge $\kappa$.
The definitions just explained make $\catR$ into a Kac-Moody categorification with the parameters defined from
$$
Q_{i,j}(x,y) := 
\begin{cases}
0&\text{if $i=j$}\\
y-x&\text{if $i = j+ 1 \neq j-1$}\\
x-y&\text{if $i = j- 1 \neq j+1$}\\
(y-x)(x-y)&\text{if $i = j + 1 = j-1$}\\
1&\text{if $i \neq j,j \pm 1$.}
\end{cases}
$$
\end{theo}

\begin{proof}[Sketch proof]
We have in our hands all of the required data, 
and it just remains to 
verify the conditions of (KM0)--(KM3) from \cref{kmcatdef}.
The axiom (KM0) follows by \cite[Lem.~4.7]{HKM}; the proof is an application of the bubble slide relation \cref{Hbubslide}.
We already constructed the adjunctions required for (KM1).
The axiom (KM2) follows from
the isomorphism established in \cite{BKiso}.
Alternatively, one can simply verify the QHA relations from scratch---with \cref{muckypups} in view, and using string diagrams rather than following the algebraic approach of \cite{BKiso}, it becomes a routine calculation. More conceptual approaches to the calculation are available in the literature, but they do not save any time compared to this most naive direct assault!
This leaves the hardest axiom (KM3), the inversion relation.
Using \cref{upwardmagicalt,Hrightpivot}, the rightward Kac-Moody crossing defined via \cref{rightpivot} is
\begin{equation}\label{rightwardmagic}
\begin{tikzpicture}[KM,centerzero,scale=.9]
\draw[-to] (-0.4,-0.4) \botlabel{j} -- (0.4,0.4);
\draw[to-] (0.4,-0.4) \botlabel{i} -- (-0.4,0.4);
\end{tikzpicture}
=
\begin{dcases}
\begin{tikzpicture}[KM,centerzero,scale=.9]
\draw[to-] (0.4,-.4)\botlabel{i} to (-0.4,.4)\toplabel{i};
\draw[-to] (-0.4,-.4)\botlabel{i} to (0.4,.4)\toplabel{i};
\projcr{0,0};
\pinpin{-.25,-.25}{.25,-.25}{1.5,-.25}{(1+y-x)^{-1}};
\end{tikzpicture}
&\text{if $j=i$}\\
\begin{tikzpicture}[KM,centerzero,scale=.9]
\draw[to-] (0.4,-.4)\botlabel{i} to (-0.4,.4)\toplabel{i-1};
\draw[-to] (-0.4,-.4)\botlabel{i-1} to (0.4,.4)\toplabel{i};
\projcr{0,0};
\pinpin{-.25,-.25}{.25,-.25}{1.25,-.25}{1+y-x};
\end{tikzpicture}
&\text{if $i=j+1$}\\
-\ \begin{tikzpicture}[KM,centerzero,scale=.9]
\draw[to-] (0.4,-.4) \botlabel{i} to (-0.4,.4)\toplabel{i};
\draw[-to] (-0.4,-.4) \botlabel{j} to (0.4,.4)\toplabel{j};
\projcr{0,0};
\pinpin{-.25,-.25}{.25,-.25}{2.85,-.25}{(i-j+y-x)(i-j-1+y-x)^{-1}};
\end{tikzpicture}
&\text{if $i \neq j,j+1$.}
\end{dcases}
\end{equation}
If $i \neq j$, it is easy to show that the rightward crossing
$\begin{tikzpicture}[KM,centerzero,scale=.7]
\draw[to-] (0.3,-.3)\botlabel{i} to (-0.3,.3)\toplabel{i};
\draw[-to] (-0.3,-.3)\botlabel{j} to (0.3,.3)\toplabel{j};
\projcr{0,0};
\end{tikzpicture}$ is invertible, with two-sided inverse given by the leftward crossing
$\begin{tikzpicture}[KM,centerzero,scale=.7]
\draw[-to] (0.3,-.3)\botlabel{j} to (-0.3,.3)\toplabel{j};
\draw[to-] (-0.3,-.3)\botlabel{i} to (0.3,.3)\toplabel{i};
\projcr{0,0};
\end{tikzpicture}$; see \cite[Lem.~4.8]{HKM}. From this and \cref{rightwardmagic},
it follows easily that
$\begin{tikzpicture}[KM,centerzero,scale=.6]
\draw[-to] (-0.3,-0.3) \botlabel{j} -- (0.3,0.3);
\draw[to-] (0.3,-0.3) \botlabel{i} -- (-0.3,0.3);
\end{tikzpicture}$ is invertible when $i \neq j$, checking the inversion relation in these cases.
The inversion relation when $i=j$ follows instead from \cite[Lem.~4.9, Lem.~4.10]{HKM}, which are more difficult to prove.
\end{proof}

\begin{rem}
We have explained \cref{maintheorem} in its simplest vanilla
flavor. If more is known about the category $\catR$, 
the construction can be modified so as 
to incorporate extra information into the choice of the
weight lattice $X$. We will not attempt to formulate such
modifications formally, preferring to explain the idea in the context
of the examples in characteristic $p > 0$ described the next two
sections. In these, we will use an
additional {\em degree decomposition} $\catR = \bigoplus_{n \in \Z} \catR_n$
to upgrade $X$ from the weight lattice of the Kac-Moody algebra
$\widehat{\mathfrak{sl}}_p'$ to that of $\widehat{\mathfrak{sl}}_p$.
\end{rem}

\begin{rem}
Another point we would like to make is that we have formulated the definition of Heisenberg categorification in the context of locally finite $\kk$-linear Abelian categories. There is a parallel theory 
in which $\catR$ is a $\kk$-linear {\em Schurian category}, that is, it is equivalent to the category of locally finite-dimensional modules over a locally finite-dimensional locally unital $\kk$-algebra.
The analogue of \cref{maintheorem} in the Schurian
setting makes the full subcategory of $\catR$ consisting of
finitely generated projective objects into a finite-dimensional
Karoubian Kac-Moody categorification as in \cref{schurian}. This is discussed further in \cite{HKM}.
\end{rem}

%% file: s5-example.tex
\setcounter{section}{4}

%=====================================
\section{Example: representations of symmetric groups}\label{s5-example}
%=====================================

In this section, we explain the first example of a Heisenberg
categorification, which was also the original motivation for Khovanov's definition of the Heisenberg category in \cite{K}. As usual, 
let $\kk$ be an algebraically closed field of characteristic $p \geq 0$.
Let $\kk S := \bigoplus_{n \geq 0} \kk S_n$ be the direct sum of the group algebras of all symmetric groups. It is a locally unital algebra with
distinguished idempotents $1_n\:(n \in \N)$ that are the identity elements of the groups $S_n$.
We are interested in the locally finite $\kk$-linear
Abelian category $\mod{\kk S}$ of finite-dimensional left
$\kk S$-modules $V = \bigoplus_{n \geq 0} 1_n V$.
Identifying $\mod{\kk S_n}$ with the full subcategory of $\mod{\kk S}$ consisting of the modules $V$ such that $1_n V = V$, we have the internal direct sum decomposition
\begin{equation}\label{degreedec}
\mod{\kk S} = \bigoplus_{n \geq 0} \mod{\kk S_n}
\end{equation}
We make $\mod{\kk S}$ into a degenerate Heisenberg categorification of central charge $\kappa := -1$ in the following way:
\begin{itemize}
\item
Let $E:\mod{\kk S} \rightarrow \mod{\kk S}$ be the $\kk$-linear functor defined on $\mod{\kk S_n}$
by the induction functor $\operatorname{ind}_{S_n}^{S_{n+1}} := \kk S_{n+1} \otimes_{\kk S_n} -$ (the embedding of $S_n$ into $S_{n+1}$ is the usual one fixing the last letter $n+1$).
This functor has the natural right adjoint $F$ defined 
on $\mod{\kk S_{n+1}}$ 
by the restriction functor $\operatorname{res}^{S_{n+1}}_{S_{n}}$,
and taking objects of $\mod{\kk S_0}$ to $\{0\}$.
\item The endomorphism
$\begin{tikzpicture}[H,centerzero]
\draw[-to] (0,-0.2) -- (0,0.2);
\opendot{0,0};
\end{tikzpicture}:E \Rightarrow E$ is defined on objects of $\mod{\kk S_n}$ by the endomorphism of $\operatorname{ind}_{S_n}^{S_{n+1}}$ arising from the endomorphism of the $(\kk S_{n+1}, \kk S_n)$-bimodule $\kk S_{n+1}$
given by right multiplication by the Jucys-Murphy element $x_{n+1} := \sum_{i=1}^n (i\:\:n\!+\!1)$.
\item The endomorphism
$\begin{tikzpicture}[H,centerzero,scale=.9]
\draw[-to] (-0.2,-0.2) -- (0.2,0.2);
\draw[-to] (0.2,-0.2) -- (-0.2,0.2);
\end{tikzpicture}:E\circ E \Rightarrow E\circ E$ is defined on objects of $\mod{\kk S_n}$ by the endomorphism of $\operatorname{ind}_{S_{n+1}}^{S_{n+2}}
\circ \operatorname{ind}_{S_n}^{S_{n+1}}
\cong \operatorname{ind}_{S_n}^{S_{n+2}}$ arising from the endomorphism of the $(\kk S_{n+2}, \kk S_n)$-bimodule $\kk S_{n+2}$
given by right multiplication by the transposition
$(n\!+\!1\:\:n\!+\!2)$.
\end{itemize}
We still need to check (H2) and (H3) from \cref{hcatdef}.
For (H2), the relations \cref{Hdotslide,Hquadratic,Hbraid} are straightforward to verify.
For example, \cref{Hdotslide} holds because
$(n\!+\!1\:\:n\!+\!2) x_{n+2} - x_{n+1} (n\!+\!1\:\:n\!+\!2) = 1_{n+2}$
in $\kk S_{n+2}$.
For the inversion relation (H3), we need to show that
$$
\begin{pmatrix}
\ \begin{tikzpicture}[H,anchorbase]
\draw[-to] (-0.25,-0.25) -- (0.25,0.25);
\draw[to-] (0.25,-0.25) -- (-0.25,0.25);
\end{tikzpicture} &
\begin{tikzpicture}[H,anchorbase]
\draw[-to] (-0.25,0.15) to[out=-90,in=-90,looseness=3] (0.25,0.15);
\node at (0,.2) {$\phantom.$};\node at (0,-.3) {$\phantom.$};
\end{tikzpicture}\ 
\end{pmatrix}:
E \circ F \oplus \id_{\mod{\kk S}} \Rightarrow F \circ E
$$
is an isomorphism.
On calculating the rightward crossing explicitly, this reduces to showing that the $(\kk S_n, \kk S_n)$-bimodule homomorphism
\begin{align*}
\theta:\kk S_{n} \otimes_{\kk S_{n-1}} \kk S_{n} \oplus \kk S_n
&\rightarrow \kk S_{n+1},&
\left(\begin{smallmatrix}g_1 \otimes g_2\\g_3\end{smallmatrix}\right) &\mapsto g_1 (n\:\:n\!+\!1) g_2 +g_3
\end{align*}
is an isomorphism for each $n \geq 0$; in the case $n=0$, the first summand should be interpreted as zero. The vectors 
$\left(\begin{smallmatrix}(i\:\:i+1\:\:\cdots\:\:n) \otimes 1_n\\0\end{smallmatrix}\right)$ for $1 \leq i \leq n$ together with $\left(\begin{smallmatrix}0\\1_n\end{smallmatrix}\right)$
give a basis for 
$\kk S_{n} \otimes_{\kk S_{n-1}} \kk S_{n} \oplus \kk S_n$
as a free right $\kk S_n$-module. The homomorphism $\theta$ maps them to
$(i\:\:i\!+\!1\:\:\cdots\:\:n\!+\!1)\in\kk S_{n+1}$ for $1 \leq i \leq n+1$,
which give a basis for $\kk S_{n+1}$
as a free right $\kk S_n$-module. So $\theta$ is indeed an isomorphism.
We have just proved that 
$$
\operatorname{ind}_{S_{n-1}}^{S_n} \circ \operatorname{res}^{S_n}_{S_{n-1}} \oplus \id_{\mod{\kk S_n}} \cong \operatorname{res}^{S_{n+1}}_{S_n} \circ\operatorname{ind}_{S_n}^{S_{n+1}},
$$
which is a special case of the Mackey theorem.

Now we apply the general construction from \cref{s4-HtoKM} to
equip $\mod{\kk S}$ with the structure of a Kac-Moody categorication. Hence, we can freely apply the general results \cref{noods1,noods2,noods3,noods4,noods5,noods6,noods7,noods8,noods9,noods10} about Kac-Moody
categorifications from the end of \cref{s2-KM}. 
Let 
$L(b)\:(b \in \B_n)$ be a full set of pairwise inequivalent irreducible left $\kk S_n$-modules. Thus, the set 
$\B := \B_0 \sqcup \B_1 \sqcup\cdots$ indexes a full set $L(b)\:(b \in \B)$ of irreducibles in $\mod{\kk S}$.
Of course, there is a well-known parametrization by certain partitions, but we do not want
to assume knowledge of any such explicit construction yet.
However, we will
use the notation $\varnothing$ to denote the unique
element of $\B_0$, so $L(\varnothing) \cong \kk S_0$.

The next step is to determine the spectrum $I$.
By a {\em generator} for a Heisenberg categorification
we mean an object $P$ such that the set of objects obtained by applying all finite sequences of the functors 
$E$ and $F$ to $P$ is a generating family 
in the usual sense of Abelian categories.

\begin{lem}\label{speclem}
Let $P$ be a generator for a Heisenberg categorification $\catR$.
The spectrum $I$ of $\catR$ is the closure of the set of roots of the minimal polynomials of the endomorphisms 
$\begin{tikzpicture}[H,anchorbase]
\draw[-to] (0,-0.2) -- (0,0.2);
\opendot{0,0};
\end{tikzpicture}:EP \rightarrow EP$ and
$\begin{tikzpicture}[H,anchorbase]
\draw[to-] (0,-0.2) -- (0,0.2);
\opendot{0,0};
\end{tikzpicture}:FP \rightarrow FP$ under the 
operation $i \mapsto i \pm 1$.
\end{lem}

\begin{proof}
Let $I$ be the spectrum of $\catR$ and $J$ be the closure of the set of roots of the minimal polynomials in the statement of the lemma.
We have already observed that the spectrum is closed under the operation
$i \mapsto i \pm 1$, hence, $I \supseteq J$.
To complete the proof, we assume for a contradiction that $I \neq J$.
Because $P$ is a generator, 
we can find some $i \in I - J$ and a finite composition $G$ of the functors
$E_j, F_j\:(j \in J)$ such that either $E_i G P$ or $F_i G P$ is non-zero.
We have that $a_{i,j} = 0$ for all $j \in J$ so,
by the defining relations for $\UU$, 
$E_i$ and $F_i$ both commute with $G$ (up to isomorphism). We deduce that either $E_i P$ or $F_i P$ is non-zero, hence, $i$ is a root of one of the minimal polynomials used to define $J$. This is a contradiction.
\end{proof}

To apply this to the example in hand, note that
the irreducible module $L(\varnothing)$ is a generator because $\kk S_n \cong E^n (\kk S_0)$ for $n \geq 0$.
Also, the minimal polynomials $m_\varnothing(x)$ 
and $n_\varnothing(x)$ are easily checked to be $x$ and $1$, respectively. Applying \cref{speclem}, we deduce that the spectrum $I$
of $\catR$ is the closure of $\{0\}$, so it is $\Z$ if $p =0$
or the prime subfield $\Z / p \Z$ of $\kk$ if $p > 0$. 
Taking $X$ be the weight lattice exactly as in \cref{wtlattice},
the locally unital $\C$-algebra $\dot\U$ from \cref{s2-KM}
is the modified form of the universal enveloping algebra of the Lie algebra $\mathfrak{sl}_\infty$ if $p=0$, or of the Lie algebra $\widehat{\mathfrak{sl}}_p'$ if $p > 0$. 

Everything is now set up the way we want in the $p=0$ case, but when $p > 0$ it is helpful to modify the general
construction slightly so that $\dot\U$ becomes the modified
form of the universal enveloping algebra of 
$\widehat{\mathfrak{sl}}_p$, which has the additional scaling element $d$ compared to $\widehat{\mathfrak{sl}}'_p$.
This is done by using some extra structure
of $\catR = \mod{\kk S}$:
letting $\catR_n := \mod{\kk S_n}$ if $n \geq 0$, or the full subcategory consisting of all of the zero objects if $n < 0$, we have the decomposition
\begin{equation}\label{degdef}
\catR = \bigoplus_{n \in \N} \catR_n
\end{equation}
as an internal direct sum of Serre subcategories.
Moreover, the following additional axiom, which is analogous to (KM0), is obviously satisfied:
\begin{itemize}
\item[(H0)]\label{h0} $E$ maps objects of $\catR_n$ to $\catR_{n+1}$; equivalently, $F$ maps objects of
$\catR_{n+1}$ to $\catR_{n}$.
\end{itemize}
We incorporate this extra degree information into the Kac-Moody categorification by extending the weight lattice to
the free Abelian group
\begin{equation}\label{enlarged}
X := {\textstyle\frac{1}{2p}}\Z\delta 
\oplus \bigoplus_{i \in I} \Z \varpi_i
\end{equation}
for some additional weight $\delta$, the {\em null root}.
We redefine
the simple coroots $h_i:X\rightarrow \Z$ so that $h_i(\delta) = 0$ and $h_i(\varpi_j) = \delta_{i,j}$ for all $i,j \in I$, and
we redefine the simple roots by setting 
\begin{equation}\label{taps}
\alpha_j := \delta_{j,0} \delta + \sum_{i \in I} a_{i,j} \varpi_i.
\end{equation}
Notice now that the simple roots are linearly independent, whereas before they summed to 0.
The scaling element $d$ is a homomorphism 
$d:X \rightarrow {\textstyle\frac{1}{2p}}\Z$ with 
$d(\delta) = 1$, $d(\varpi_i) = 0$ and $d(\alpha_i) = \delta_{i,0}$.
Let
\begin{align}
\label{sisters}
\deg := pd + \sum_{i=0}^{p-1} \textstyle{\frac{i(p-i)}{2}} h_i:X &\rightarrow \textstyle{\frac{1}{2}\Z},
\end{align}
so
$\deg(\delta) = p$,
$\deg(\varpi_i) = {\textstyle\frac{i(p-i)}{2}}$ and, most importantly,
$\deg(\alpha_i) = 1$.
Then we modify the definition of the function $\wt:\B \rightarrow X$ 
from \cref{wtdef} by
redefining the weight 
$\wt(b) \in X$ of $b \in \B$ 
so that $\deg(\wt(b))$ is the degree of $L(b)$, that is, $n$ if $b \in \B_n$, and $h_i(\wt(b)) = \phi_i(b)-\eps_i(b)$
for each $i \in I$ as before.
Finally, we use this new definition of $\wt$
to obtain the weight decomposition 
\cref{blockdec} with weight subcategories parametrized now by 
the $\delta$-extended weight lattice.
Since $E_i$ takes $\kk S_n$-modules to $\kk S_{n+1}$-module
and $\deg(\alpha_i)=1$,
the functors $E_i$ and $F_i$ still satisfy (KM0). The remaining checks in the proof of \cref{maintheorem} go through as before, making $\catR$ into a Kac-Moody categorification with this extended weight lattice.

We would next like to understand the $\dot\U$-module
$\C \otimes_\Z G_0(\mod{\kk S})$. The general theory 
gives that it is an integrable $\dot \U$-module
with a perfect basis given by the isomorphism classes $[L(b)]\:(b \in \B)$. 
Since $m_\varnothing(x)=x$ and $n_\varnothing(x)=1$,
we have that 
\begin{align}\label{wts}
\eps_i(\varnothing) &= \delta_{i,0}, 
&\phi_i(\varnothing) &= 0
\end{align}
for $i \in I$.
Applying \cref{wtdef}, we deduce that
$\wt(\varnothing) = -\varpi_0$. 
%For $b \in \B_n$, $L(b)$ is a quotient of $E^n L(\varnothing)$, so $\wt(b)$ is of the form $-\varpi_0+\alpha_{i_1}+\cdots+\alpha_{i_n}$ for some $i_1,\dots,i_n \in I$.
Let $V(-\varpi_0)$ be the irreducible $\dot\U$-module
of lowest weight $-\omega_0$, that is, the quotient of 
$\dot\U 1_{-\omega_0}$ by the left ideal generated by
$e_i^{(1+\delta_{i,0})} 1_{-\omega_0}$ and 
$f_i 1_{-\omega_0}$ for all $i \in I$.
By \cref{wts}, 
there is a injective $\dot\U$-module homomorphism
$V(-\varpi_0) \rightarrow \C\otimes_\Q G_0(\mod{\kk S})$
mapping the generator of $V(-\varpi_0)$ to $[L(\varnothing)]$.
In fact, this homomorphism is an isomorphism, so that
\begin{equation}\label{chickenfood}
\C\otimes_\Z G_0(\mod{\kk S}) \cong V(-\varpi_0).
\end{equation}
There are several ways to prove this, depending on how much one is willing to assume about integrable representations of $\dot\U$.
The simplest approach is to compare the dimensions 
of the sum of the weight spaces of these two modules
for all weights $\lambda \in X$ with $\deg(\lambda)=n$,
that is, the weights $\lambda$ of the form $-\varpi_0+\alpha_{i_1}+\cdots+\alpha_{i_n}$ for $i_1,\dots,i_n \in I$.
For the module $\C\otimes_\Z G_0(\mod{\kk S})$, this is
the dimension of the subspace $\C\otimes_\Z G_0(\mod{\kk S_n})$, which 
by some basic representation theory of finite groups
is equal to the number of conjugacy classes of $S_n$ if $p=0$, or the number of $p$-regular conjugacy classes of $S_n$ if $p > 0$. This is the number of partitions of $n$ if $p=0$ or the number of $p$-regular partitions of $n$ if $p > 0$.
It remains to observe that the same number computes the dimension of the sum of the degree 
$n$ weight spaces of $V(-\varpi_0)$, hence, our homomorphism is an isomorphism.
This argument uses some knowledge of the 
basic representation $V(-\varpi_0)$, which we discuss a little further
below.

We focus in this paragraph on the $p=0$ case.
Then $V(-\varpi_0)$ is 
a minuscule representation of $\mathfrak{sl}_\infty$ with 
the following completely explicit construction: it is the $\C$-vector space with basis $v_\lambda$ indexed by the set of partitions, the vector $v_\varnothing$ is of weight 
$-\varpi_0$, and the Chevalley generators $e_i$ (resp.,  $f_i$)
act by mapping $v_\lambda$ to $v_\mu$ for $\mu$ obtained adding (resp.,  removing)
a node of content $i$ to (resp., from) 
the Young diagram\footnote{We use the English convention for Young diagrams, in which the content of a node is its column number minus its row number.} of $\lambda$ if such a node exists, or to 0 otherwise.
The same combinatorics gives a realization of the crystal
$\mathcal B(-\varpi_0)$ of $V(-\varpi_0)$ as Young's partition lattice.
Since all perfect bases for an
integrable lowest or highest weight module have canonically isomorphic underlying crystals by \cite[Th.~5.37]{BeK}, the crystal
$\mathcal B(-\varpi_0)$ is also the associated crystal of the Kac-Moody categorification
$\mod{\kk S}$.
It follows that the labelling 
set $\B$ is identified with the set of partitions, and we have in our hands explicit 
branching rules for the $i$-induction and $i$-restriction functors $E_i$ and $F_i$. The fact that the labelling of irreducible $\kk S_n$-modules thus obtained coincides with the standard labelling follows by comparing these branching rules with the classically known ones. The Okounkov-Vershik approach from \cite{OV} is particularly well suited for this comparison.

When $p > 0$, the crystal $\mathcal B(-\varpi_0)$ of the basic representation
$V(-\varpi_0)$ was described explicitly in
\cite{MM}. 
Its vertex set is naturally identified with the set of $p$-regular partitions; using this, the counting argument used to prove \cref{chickenfood} is clear.
Now, the result from \cite[Th.~5.37]{BeK} about uniqueness of crystals implies that the crystal associated to the Kac-Moody categorification
$\mod{\kk S}$ is identified with $\mathcal B(-\varpi_0)$.
With these arguments, which used some basic facts about integrable
representations of $\widehat{\mathfrak{sl}}_p$ but almost no knowledge of the modular representation theory of symmetric groups, we have constructed an explicit
parametrization of the isomorphism classes of irreducible representations of the symmetric
groups over the field $\kk$ of characteristic $p > 0$ 
by the set of $p$-regular partitions, and modular
branching rules for the $i$-induction and $i$-restriction functors.

Comparing with the modular branching rules proved in
\cite{Klesh}, it follows that
the labelling of irreducibles that we have constructed coincides with the usual labelling from \cite{James}. The results in \cite{Klesh}
are themselves not so easy to prove (see also \cref{losevify}),
but we have only used them in this final step to identify the
labelling arising from the general theory of categorification with the
classical labelling. 
Many other fundamental results about modular
representations of symmetric groups follow immediately from the general theory. 
The classical Nakayama Conjecture can be interpreted as the assertion that
the non-zero weight subcategories are exactly the
blocks of $\catR$.
Then property \ref{noods9} from \cref{s2-KM} recovers the Morita equivalences
between blocks of symmetric groups
constructed originally by Scopes \cite{Scopes}, and \ref{noods10} gives the
derived equivalences between blocks from \cite{CR}.

%% file: s6-example.tex
\setcounter{section}{5}

%=====================================
\section{Example: rational representations of \texorpdfstring{$\mathrm{GL}_n$}{} and its Frobenius kernels}\label{s6-example}
%=====================================

Let $\kk$ be an algebraically closed field of characteristic $p > 0$.
Let $G$ be the general linear group scheme $\mathrm{GL}_n$ over $\kk$,
fixing also the standard choices of maximal torus $T$ (diagonal matrices) and Borel subgroup $B$ (upper triangular matrices).
We begin with some reminders about the symmetric tensor category $\Rep(G)$
of finite-dimensional representations of $G$.
Let $V$ be the natural $G$-module with standard basis $v_1,\dots,v_n$;
the weight of $v_i$ is $\eps_i \in X(T)$.
We let
\begin{equation}
\rho := -  \eps_2 - 2 \eps_3 - \cdots - (n-1) \eps_n.
\end{equation}
For $\lambda \in X(T)$ and $1 \leq i \leq n$, we 
use the notation $\lambda_i$ to denote
$(\lambda+\rho, \eps_i) \in \Z$, where $(\cdot,\cdot)$ is the usual symmetric bilinear form on $X(T)$ defined so that $\eps_1,\dots,\eps_n$  are orthonormal.
The irreducible $G$-modules are parametrized by highest weight theory by the set 
\begin{equation}
X^+(T) = \{\lambda \in X(T)\:|\:\lambda_1 > \cdots > \lambda_n\}
\end{equation}
of {\em dominant weights}. We denote the irreducible $G$-module of highest weight $\lambda \in X^+(T)$ by $L(\lambda)$.

The {\em algebra of distributions} $\Dist(G)$ is identified with $\kk \otimes_\Z U_\Z$, where $U_\Z$ is the Kostant $\Z$-form for the universal enveloping algebra of the Lie algebra $\mathfrak{gl}_n(\C)$
generated by the divided power monomials 
$e_{i,j}^{(m)} := e_{i,j}^m / m!$
for $m \geq 1$ in the matrix units and 
$1 \leq i,j \leq n$ with $i \neq j$, plus the binomials
$\binom{e_{i,i}}{m}$ for $m \geq 1$ and $1 \leq i \leq n$.
It is useful because the canonical functor
from $\Rep(G)$ to $\mod{\Dist(G)}$ is an isomorphism of categories.
Note also that the Lie algebra $\mathfrak{g}$ of $G$ is identified with the Lie subalgebra of $\Dist(G)$
with basis $e_{i,j}\:(1 \leq i,j \leq n)$.
Let
\begin{equation}
\Omega := \sum_{i,j=1}^n e_{i,j} \otimes e_{j,i} \in \mathfrak{g} \otimes \mathfrak{g}
\end{equation}
be the Casimir tensor.
Acting by $\Omega$ defines a $G$-equivariant endomorphism of any tensor product $M_1 \otimes M_2$ of $G$-modules.

Now we can make $\Rep(G)$ into a degenerate Heisenberg categorification of central charge $\kappa = 0$. Here is the required data:
\begin{itemize}
\item
The endofunctors $E$ and $F$ are $V \otimes-$ and $V^* \otimes -$, respectively. 
The adjunction $(E,F)$ is the canonical one.
\item 
The natural transformation $\begin{tikzpicture}[H,centerzero]
\draw[-to] (0,-0.2) -- (0,0.2);
\opendot{0,0};
\end{tikzpicture}:E \Rightarrow E$ is defined 
on $M \in \ob\Rep(G)$ by the endomorphism $V \otimes M \rightarrow V \otimes M$ arising from
the action of the Casimir tensor $\Omega$.
\item 
The natural transformation
$\begin{tikzpicture}[H,centerzero,scale=.9]
\draw[-to] (-0.2,-0.2) -- (0.2,0.2);
\draw[-to] (0.2,-0.2) -- (-0.2,0.2);
\end{tikzpicture}:E\otimes E \Rightarrow E \otimes E$
is defined on $M \in \ob\Rep(G)$ by the endomorphism
$V \otimes V \otimes M \rightarrow V \otimes V \otimes M,
u \otimes v \otimes m \mapsto v \otimes u \otimes m$.
\end{itemize}
The axioms of Heisenberg categorification are easy to check. The
isomorphism $V \otimes V^* \otimes M \stackrel{\sim}{\rightarrow} V^* \otimes V \otimes M$
defined by the rightward crossing is simply the tensor flip
$v \otimes f \otimes m \mapsto f \otimes v \otimes m$.

Now we apply the construction from \cref{s4-HtoKM}, making
$\Rep(G)$ into a Kac-Moody categorification.
As with any example, the next steps are to determine the spectrum $I$, the integrable representation that is the Grothendieck group of $\Rep(G)$, and the associated crystal. This requires a bit more knowledge about $\Rep(G)$. For each $\lambda \in X^+(T)$, the {\em Weyl module}
$\Delta(\lambda)$ is the universal highest weight module of highest weight $\lambda$ in $\Rep(G)$. In fact, $\Rep(G)$ is a {\em highest weight category}, and the Weyl modules are its standard modules.

\begin{lem}\label{dumb}
For $\Rep(G)$, the spectrum $I$ is $\Z / p \Z$.
For $i \in I$ and $\lambda \in X^+(T)$,
$E_i \Delta(\lambda)$ (resp., $F_i \Delta(\lambda)$) has a multiplicity-free filtration
with sections that are isomorphic to the Weyl modules
$\Delta(\lambda+\eps_j)$ (resp., $\Delta(\lambda-\eps_j)$) for $1 \leq j \leq n$ such that $\lambda+\eps_j \in X^+(T)$ and $\lambda_j \equiv i\pmod{p}$ (resp., $\lambda-\eps_j \in X^+(T)$ and $\lambda_j - 1 \equiv -n-i\pmod{p}$).
\end{lem}

\begin{proof}
We prove this for $E_i$.
The first step is to observe that $E \Delta(\lambda) = V \otimes \Delta(\lambda)$ has a multiplicity-free filtration with sections that are isomorphic to the Weyl modules
$\Delta(\lambda+\eps_j)$ for all $1 \leq j \leq n$ such that $\lambda+\eps_j \in X^+(T)$.
To see this, we let $M_0 := \{0\}$, then for $j=1,\dots,n$ we recursively define
$M_j$ to be the submodule of $V \otimes \Delta(\lambda)$ generated by 
$M_{j-1}$ and $v_j \otimes v_+$, where $v_+$ is a highest weight vector. The image of $v_j \otimes v_+$ in $M_j / M_{j-1}$ is easily seen to be a (possibly zero) highest weight vector of weight $\lambda+\eps_j$. Hence, $M_j / M_{j-1}$ is a quotient of $\Delta(\lambda+\eps_j)$ if $\lambda+\eps_j$ is a dominant weight, and it is zero 
if $\lambda + \eps_j$ is not a dominant weight.
It remains to observe that the sections $M_j / M_{j-1}$
are actually isomorphic to $\Delta(\lambda+\eps_j)$ (rather than being proper quotients) for all $j$ with $\lambda+\eps_j \in X^+(T)$. This follows by a calculation using Weyl's dimension formula.

Now we take $i \in \kk$ and consider $E_i \Delta(\lambda)$, i.e., 
the projection of $V \otimes \Delta(\lambda)$ onto the generalized $i$-eigenspace of the endomorphism $\Omega$. Let $z_r$ be the unique element of the
Harish-Chandra center of $\Dist(G)$ which acts on $\Delta(\lambda)$ by ${\mathrm e}_r(\lambda)$, the $r$th elementary symmetric polynomial evaluated at the scalars $\lambda_1,\dots,\lambda_n$.
By \cite[Lem.~5.1]{cyclo}, $z_2$ acts on $V \otimes \Delta(\lambda)$ in the same way
as $1 \otimes z_2 + 1 \otimes z_1 - \Omega$.
Consequently, $\Omega$ leaves each of the submodules $M_j$ invariant and, on an $\Omega$-eigenvector of eigenvalue $i$
in $M_j / M_{j-1} \cong \Delta(\lambda+\eps_j)$,
we have that 
$$
{\mathrm e}_2(\lambda+\eps_j) = {\mathrm e}_2(\lambda) + {\mathrm e}_1(\lambda)-i.
$$
This identity implies that $i = \lambda_j\text{ mod }p$. We deduce that $E_i \Delta(\lambda)$ is zero unless $i \in \Z / p\Z$, hence, $I = \Z / p \Z$. Moreover, for $i \in \Z / p\Z$, we have shown that $E_i \Delta(\lambda)$ has the filtration claimed.

The proof for $F_i$ is very similar, using instead that $z_2$ acts
on $V^* \otimes \Delta(\lambda)$ in the same way as
$1 \otimes z_2 + 1 \otimes (1-n-z_1) - \Omega$.
\end{proof}

Henceforth, $I := \Z / p\Z \subset \kk$.
Taking the vanilla flavor of Kac-Moody categorification with weight lattice $X = \bigoplus_{i \in I} \varpi_i$, the 
algebra $\dot \U$ is the modified form of the universal enveloping algebra of $\widehat{\mathfrak{sl}}_p'$. 
It is a tolerable 
nuisance that there are two weight lattices in play now, $X(T)$ for $G$ and $X$ for $\dot \U_\Z$,
so that $\lambda$ could either mean an element of $X(T)$ or an element of the weight lattice $X$ depending on the context.
We use $\M_\Z$ to denote the level zero representation of $\dot\U_\Z$ which is the restricted dual of the natural representation of $\widehat{\mathfrak{sl}}_p'$. It is a free $\Z$-module with basis $m_j\:(j \in \Z)$, on which the Chevalley generators act by
\begin{align}\label{thursday}
e_i m_j &= \begin{cases}
m_{j+1}&\text{if $j \equiv i \pmod{p}$}\\
0&\text{otherwise,}
\end{cases}&
f_i m_j &= \begin{cases}
m_{j-1}&\text{if $j \equiv i+1 \pmod{p}$}\\
0&\text{otherwise.}
\end{cases}
\end{align}
The weight of $m_j$ is $-\epsilon_j$ where $\epsilon_j := \varpi_{j\text{ mod }p} - \varpi_{(j-1)\text{ mod }p}$.
As $\Rep(G)$ is a highest weight category, the isomorphism classes of the Weyl modules give a basis for $G_0(\Rep(G))$. 
From \cref{dumb}, it follows that there is an isomorphism of $\dot\U_\Z$-modules
\begin{align}\label{giardia}
G_0(\Rep(G)) &\stackrel{\sim}{\rightarrow} {\bigwedge}^n \M_\Z,
&
[\Delta(\lambda)] &\mapsto m_{\lambda_1} \wedge\cdots\wedge m_{\lambda_n}.
\end{align}

To complete this first part of the story, we describe the associated crystal $(X^+(T), \tilde e_i, \tilde f_i, \eps_i, \phi_i, \wt)$.
Recall that Kashiwara introduced an associative (but not commutative) tensor product operation on normal crystals; see \cite{Kashcrystal}. Given normal crystals $\B_1$ and $\B_2$, their tensor product
$\B_1 \otimes \B_2$ is the normal crystal with underlying set 
$\{b_1\otimes b_2\:|\:b_1 \in \B_1, b_2 \in \B_2\}$
with $\wt(b_1 \otimes b_2) := \wt(b_1)+\wt(b_2)$.
The crystal operators $\tilde e_i, \tilde f_i$ are defined by
\begin{align}
\tilde e_i(b_1\otimes b_2) &= 
\begin{cases}
b_1 \otimes \tilde e_i(b_2)&\text{if $\eps_i(b_2) > \phi_i(b_1)$}\\
\tilde e_i(b_1) \otimes b_2&\text{if $\eps_i(b_2) \leq \phi_i(b_1)$,}
\end{cases}\\
\tilde f_i(b_1\otimes b_2) &= 
\begin{cases}
\tilde f_i(b_1) \otimes b_2&\text{if $\phi_i(b_1) > \eps_i(b_2)$}\\
b_1 \otimes \tilde f_i(b_2)&\text{if $\phi_i(b_1) \leq \eps_i(b_2)$,}
\end{cases}
\end{align}
assuming this makes sense (the crystal operators are undefined otherwise).
Associated to the based 
module $\M_\Z$, there is a normal crystal with vertex set
$\Z$, 
crystal operators defined by $\tilde e_i(j) = j+1 \Leftrightarrow j \equiv i \pmod{p}$, $\tilde f_i(j) = j-1\Leftrightarrow j \equiv i+1\pmod{p}$,
and $\wt:\Z \rightarrow X$ defined by
$\wt(j) := -\epsilon_j$.
Taking the $n$th tensor power of this gives a crystal
whose underlying set is $\Z^n$. We transport this to $X(T)$ using the bijection
$X(T) \stackrel{\sim}{\rightarrow} \Z^n, 
\lambda \mapsto (\lambda_1,\dots,\lambda_n)$.

\begin{rem}\label{exp}
Here is a more user-friendly 
description of the crystal operators $\tilde e_i$ and $\tilde f_i$
on $\lambda \in X(T)$.
Let $(\sigma_1,\dots,\sigma_n) \in \{0,+,-\}^n$ be the {\em reduced $i$-signature} of $\lambda$.
This is obtained by starting from the sequence with
$\sigma_j := +$ if $\lambda_j = i$,
$\sigma_j := -$ if $\lambda_j-1 = i$, and $\sigma_j := 0$ otherwise,
then repeatedly replacing subsequences of the form $-,0,\dots,0,+$
with $0,0,\dots,0,0$ so that at the end there is no $-$ appearing to the left of a $+$.
Then $\eps_i(\lambda)$ is the number of entries equal to $+$ and, if this is non-zero, $\tilde e_i(\lambda) = \lambda+\eps_j$ where $j$ is the index of the rightmost $+$. Similarly, 
$\phi_i(\lambda)$ is the number of entries equal to $-$ and, if this is non-zero, $\tilde f_i(\lambda) = \lambda-\eps_j$ where $j$ is the index of the leftmost $-$.
\end{rem}

\begin{theo}\label{itwas}
The crystal associated to the Heisenberg categorification
$\Rep(G)$ is the connected component of the 
$\widehat{\mathfrak{sl}}_p'$-crystal 
$(X(T), \tilde e_i, \tilde f_i, \eps_i, \phi_i, \wt)$ 
just described with vertex set $X^+(T) \subset X(T)$.
\end{theo}

\begin{proof}
This follows from \cite{Klesh} using the explicit description in \cref{exp}.
\end{proof}

Everything so far is well known. To say something not quite so standard, we look instead at Frobenius kernels. The basic reference  is \cite[Ch.~II.3]{Jantzen}.
The {\em Frobenius morphism} $F:G \rightarrow G$ is the morphism of group schemes defined by raising matrix entries to the $p$th power.
Let 
$$
G_r := \ker F^r
$$ 
be the $r$th {\em Frobenius kernel}, which is a closed normal subgroup scheme.
Also set $B_r := B \cap G_r$ and $T_r := T \cap G_r$.
The category $\Rep(G_r)$ of finite-dimensional rational representations of $G_r$ can be viewed equivalently as the category
$\mod{\Dist(G_r)}$, where $\Dist(G_r)$ is the subalgebra of $\Dist(G)$
generated by the divided powers $e_{i,j}^{(m)}\:(i \neq j)$ and the binomials $\binom{e_{i,i}}{m}$ for all $m$ with $1 \leq m \leq p^r-1$. 
It is a finite-dimensional algebra of dimension $p^{r n^2}$.

Let $X(T_r) := X(T) / p^r X(T)$. This Abelian group naturally 
labels the one-dimensional representations of $T_r$.
For $\bar\lambda \in X(T_r)$ with pre-image $\lambda \in X(T)$, we use $\bar\lambda_i$
to denote the unique 
element of $\Z / p^r \Z$ such that $(\lambda +\rho,\eps_i) \equiv \bar\lambda_i \pmod{p^r}$.
Let $\kk_{\bar\lambda}$ be the $B_r$-module that is the inflation of the one-dimensional 
$T_r$-module corresponding to $\bar\lambda\in X(T_r)$.
The {\em baby Verma module} $Z_r(\bar\lambda)$ is
the $G_r$-module $\Dist(G_r) \otimes_{\Dist(B_r)} \kk_{\bar\lambda}$. It has simple head $L_r(\bar\lambda)$, and the modules $L_r(\bar\lambda)\:(\bar\lambda \in X(T_r))$ give representatives for the isomorphism classes of irreducible $G_r$-modules.

We make $\Rep(G_r)$ into a degenerate Heisenberg categorification of central charge 0 in exactly the same way as we did for $\Rep(G)$.
Hence, it is a Kac-Moody categorification by the construction from \cref{s4-HtoKM}. The counterpart of \cref{dumb} is as follows.

\begin{lem}\label{dumber}
For $\Rep(G_r)$, the spectrum $I$ is $\Z / p \Z$.
For $i \in I$ and $\bar\lambda \in X(T_r)$,
$E_i Z_r(\bar\lambda)$ (resp., $F_i Z_r(\bar\lambda)$) has a multiplicity-free filtration
with sections that are isomorphic to the baby Verma modules
$Z_r(\bar\lambda+\eps_j)$ (resp., $Z_r(\bar\lambda-\eps_j)$) for $1 \leq j \leq n$ such that $\bar\lambda_j \equiv i\pmod{p}$ (resp., $\bar\lambda_j - 1 \equiv -n-i\pmod{p}$).
\end{lem}

\begin{proof}
This follows from the same argument as used to prove \cref{dumb}. 
It is even a bit easier since there is no dominance constraint,
and $\dim Z_r(\bar\lambda) = p^{r n(n-1)/2}$ always so that the dimension calculation is quite trivial.
\end{proof}

\cref{dumber} allows us to describe the Grothendieck group as a $\dot\U_\Z$-module:
letting 
$\M^{(r)}_\Z$ be the level zero representation of $\dot\U_\Z$ 
which is free as a $\Z$-module with basis $m^{(r)}_j\:(j \in \Z / p^r \Z)$
and action of Chevalley generators defined analogously to \cref{thursday},
there is an isomorphism
\begin{align}
K_0(\Rep(G_r)) &\stackrel{\sim}{\rightarrow} 
\left(\M^{(r)}_\Z\right)^{\otimes n},
&
[Z_r(\bar\lambda)]\ &\mapsto m^{(r)}_{\bar\lambda_1} \otimes \cdots \otimes m^{(r)}_{\bar\lambda_n}.
\end{align}
Associated to the module $\M^{(r)}_\Z$, there is a normal $\widehat{\mathfrak{sl}}_p'$-crystal
with vertex set $\Z / p^r \Z$. The crystal operators and
weight function are defined analogously to the crystal of $\M_\Z$ described above.
We take the $n$th tensor power of this to obtain a crystal with underlying set $(\Z / p^r \Z)^n$, then transport this to $X(T_r)$
using the bijection $X(T_r) \stackrel{\sim}{\rightarrow}
(\Z/ p^r \Z)^n, \bar\lambda \mapsto (\bar\lambda_1,\dots,\bar\lambda_n)$.
It is easy to see from this definition that the quotient map
$X(T) \twoheadrightarrow X(T_r), \lambda \mapsto \bar\lambda$ is a morphism of $\widehat{\mathfrak{sl}}_p'$-crystals.

\begin{theo}\label{tpc}
The crystal associated to the Heisenberg categorification $\Rep(G_r)$ is the crystal with underlying set $X(T_r)$ just described,
that is, it is the $n$th tensor power of the $\widehat{\mathfrak{sl}}_p'$-crystal of $\M^{(r)}_\Z$.
\end{theo}

\begin{proof}
We will deduce this from \cref{itwas} using the following basic facts.
A weight $\lambda \in X^+(T)$ is {\em $p^r$-restricted}
if $\lambda_i - \lambda_{i+1} \leq p^r$ for $i=1,\dots,n-1$.
By \cite[Prop.~II.3.15]{Jantzen} and the Steinberg tensor product theorem, 
for $\lambda \in X^+(T)$, the restriction
$\operatorname{res}^G_{G_r} L(\lambda)$ is completely reducible, and it is irreducible if and only if $\lambda$ is $p^r$-restricted, in which case it is isomorphic to $L(\bar \lambda)$ where $\bar\lambda$ is the canonical 
image of $\lambda$ in $X(T_r)$.

Now take $\bar\lambda \in X(T_r)$ and let
$\lambda \in X^+(T)$ be a $p^r$-restricted lift. 
We know that the $G$-module $E_i L(\lambda)$ is non-zero if and only if
$\eps_i(\lambda) > 0$, in which case it has irreducible head $\cong L(\tilde e_i(\lambda))$.
Hence $E_i L_r(\bar \lambda)\cong E_i \left(\operatorname{res}^G_{G_r} L(\lambda)\right) \cong \operatorname{res}^G_{G_r}
\left(E_i L(\lambda)\right)$ is non-zero if and only if 
$\eps_i(\lambda) > 0$. 
Assuming it is non-zero, we also know that the $G_r$-module 
$E_i L_r(\bar \lambda)$ has irreducible head. 
This implies that $\tilde e_i(\lambda)$ must be $p^r$-restricted (this can also be seen directly using the combinatorics from
\cref{exp}) and the irreducible head of $E_i L_r(\bar\lambda)$ is isomorphic to $L_r(\overline{\tilde e_i(\lambda)}) \cong \operatorname{res}^{G}_{G_r} L(\tilde e_i(\lambda))$. The lemma follows since
the map $X(T) \twoheadrightarrow X(T_r)$ is a crystal morphism.
\end{proof}

The category $\Rep(G)$ has a natural degree decomposition defined by the action of the 
one-dimensional torus that is the center of $G$, e.g.,
the trivial module is in degree 0.
We find it more convenient to shift degrees so as to put the trivial module instead into degree $-\frac{n(n-1)}{2}$.
Applying $E$ adds 1 to the degree of a homogeneous 
object, i.e., the axiom (H0) formulated after \cref{degdef} holds.
Consequently,
we can modify the above construction to make $\Rep(G)$ into a
Kac-Moody categorification for
the extended weight lattice from \cref{enlarged}. 
This is done in the same way as explained in the previous section,
redefining simple roots and simple coroots as there, and 
modifying the weight function
$\wt:X^+(T) \rightarrow X$ so that $\deg(\wt(\lambda))$ is the degree of $L(\lambda)$, and 
$h_i(\wt(\lambda)) = \phi_i(\lambda)-\eps_i(\lambda)$
for each $i$ as usual.
We redefine the weights $\epsilon_j$ by setting
\begin{equation}\label{food}\textstyle
\epsilon_j := \varpi_{j\text{ mod }p} - \varpi_{(j-1)\text{ mod }p} +\left(\frac{p-1}{2p}-\left\lceil \frac{j}{p}\right\rceil\right)\delta
\in X.
\end{equation}
The term $-\left\lceil\frac{j}{p}\right\rceil\delta$ is needed here so that
$\epsilon_j - \epsilon_{j+1} = \alpha_{j\text{ mod }p}$ for every $j \in \Z$; cf. \cref{taps}.
The additional term $\frac{p-1}{2p} \delta$ in \cref{food} is somewhat arbitrary; we included it 
so that the function $\deg$ from \cref{sisters}
satisfies the memorable 
formula $\deg(\eps_j) = -j$ for any $j \in \Z$.
We then have for any
$\lambda \in X(T)$
that
\begin{equation}\label{extendedwt}
\wt(\lambda) =-\epsilon_{\lambda_1} - \cdots - \epsilon_{\lambda_n}.
\end{equation}
Now, $\dot\U$ is the modified form for the enveloping algebra of
$\widehat{\mathfrak{sl}}_p$ (rather than $\widehat{\mathfrak{sl}}_p'$).
The action of $\dot\U_\Z$ on $\M_\Z$ from before extends by declaring
that the weight of the basis vector $m_j$ is the new
$-\epsilon_j$.
It is again the case that
\begin{align}
K_0(\Rep(G)) &\stackrel{\sim}{\rightarrow} {\bigwedge}^n \M_\Z,
&
[\Delta(\lambda)] &\mapsto m_{\lambda_1}\wedge\cdots\wedge m_{\lambda_n}.
\end{align}
The associated crystal has the same crystal operators as in \cref{itwas}, but using 
\cref{extendedwt} as its weight function so that 
it is now an $\widehat{\mathfrak{sl}}_p$-crystal.

For Frobenius kernels, there is no degree decomposition, but this can be fixed by working instead with the {\em thickened Frobenius kernels}
$G_r T$ as in \cite[Sec.~II.9]{Jantzen}.
Irreducible modules in $\Rep(G_r T)$ are parametrized by their highest weights: for $\lambda \in X(T)$, we write $\hat L_r(\lambda)$ for the irreducible $G_r T$-module of highest weight $\lambda$, which may be constructed explicitly as the unique irreducible quotient of
$\hat Z_r(\lambda) := \Dist(G_r T) \otimes_{\Dist(B_r T)} \kk_\lambda$.
Then we make $\Rep(G_r T)$ into a degenerate Heisenberg categorification of central charge zero as above.
It becomes a Kac-Moody categorification using the extended weight lattice $X$ from \cref{enlarged} and the weight function \cref{extendedwt}. Letting $\dot\U_\Z$ and $\M_\Z$ be as in the previous paragraph, 
there is an isomorphism of $\dot\U_\Z$-modules
\begin{align}
K_0(\Rep(G_r T)) &\stackrel{\sim}{\rightarrow} \M_\Z^{\otimes n},
&
[\hat Z_r(\lambda)]\ &\mapsto m_{\lambda_1} \otimes \cdots \otimes m_{\lambda_n}.
\end{align}
This follows from the obvious extended analogue of \cref{dumber}.
Also the associated crystal is the crystal 
$(X(T), \tilde e_i, \tilde f_i, \eps_i, \phi_i, \wt)$ defined just before \cref{itwas} but with $\wt$ as in \cref{extendedwt}, 
that is, it is 
the $n$th tensor power 
of the $\widehat{\mathfrak{sl}}_p$-crystal of $\M_\Z$.
This may be deduced from \cref{tpc} using that $\operatorname{res}^{G_r T}_{G_r} \hat L_r(\lambda) \cong L_r(\bar \lambda)$ for $\lambda \in X(T)$ with image $\bar\lambda \in X(T_r)$.
We observe this crystal is {\em the same} for all $r \geq 1$, although the categories $\Rep(G_r T)$ themselves are quite different as $r$ varies.

\begin{rem}\label{losevify}
In \cite{losev} (see also \cite[Sec.~7]{LW}),
Losev has developed a more conceptual approach for computing the associated crystal of a {\em tensor product categorification} in the sense of \cite[Rem.~3.6]{LW} via Kashiwara's tensor product rule.
This can be applied to $\Rep(G_r T)$, which is an $n$-fold tensor product categorification, thereby 
giving a more direct way to identify the associated crystal in this case. The earlier descriptions of associated crystals in \cref{itwas,tpc}, and also Kleshchev's result for 
symmetric groups discussed in the previous section, can be deduced 
as consequences. This is desirable since the ad hoc argument in 
\cite{Klesh} is rather complicated.
\end{rem}

%% file: s7-bubbles.tex
\setcounter{section}{6}

%=====================================
\section{Explicit formulae for the second adjunction}\label{s7-bubbles}
%=====================================

We go back now to the setup of \cref{s4-HtoKM}.
Prior to \cref{maintheorem}, we defined the Kac-Moody 
natural transformations
$\ \begin{tikzpicture}[KM,centerzero,scale=.8]
\draw[-to] (0,-0.3) \botlabel{i} -- (0,0.3);
\opendot{0,0};
\end{tikzpicture}\ ,
\ \begin{tikzpicture}[KM, centerzero,scale=.8]
\draw[to-] (0,-0.3)\botlabel{i} -- (0,0.3);
\opendot{0,0};
\end{tikzpicture}\ ,\ \begin{tikzpicture}[KM,centerzero,scale=.8]
\draw[-to] (-0.3,-0.3) \botlabel{i} -- (0.3,0.3);
\draw[-to] (0.3,-0.3) \botlabel{j} -- (-0.3,0.3);
\end{tikzpicture}\ ,\  \begin{tikzpicture}[KM,centerzero,scale=.8]
\draw[-to] (-0.25,-0.15) \botlabel{i} to [out=90,in=90,looseness=3](0.25,-0.15);
\end{tikzpicture}\ $ and $\ \begin{tikzpicture}[KM,centerzero,scale=.8]
\draw[-to] (-0.25,0.15) \toplabel{i} to[out=-90,in=-90,looseness=3] (0.25,0.15);
\end{tikzpicture}\ $.
There are also natural transformations depicted by the
leftward Kac-Moody caps and cups 
$\ \begin{tikzpicture}[KM,centerzero,scale=.8]
\draw[to-] (-0.25,-0.15) \botlabel{i} to [out=90,in=90,looseness=3](0.25,-0.15);
\end{tikzpicture}\ $
and
$\ \begin{tikzpicture}[KM,centerzero,scale=.8]
\draw[to-] (-0.25,0.15) \toplabel{i} to[out=-90,in=-90,looseness=3] (0.25,0.15);
\end{tikzpicture}\ $,
and the rightward, leftward and downward Kac-Moody crossings
$\ \begin{tikzpicture}[KM,centerzero,scale=.8]
\draw[-to] (-0.3,-0.3) \botlabel{i} -- (0.3,0.3);
\draw[to-] (0.3,-0.3) \botlabel{j} -- (-0.3,0.3);
\end{tikzpicture}\ $,
$\ \begin{tikzpicture}[KM,centerzero,scale=.8]
\draw[to-] (-0.3,-0.3) \botlabel{i} -- (0.3,0.3);
\draw[-to] (0.3,-0.3) \botlabel{j} -- (-0.3,0.3);
\end{tikzpicture}\ $
and
$\ \begin{tikzpicture}[KM,centerzero,scale=.8]
\draw[to-] (-0.3,-0.3) \botlabel{i} -- (0.3,0.3);
\draw[to-] (0.3,-0.3) \botlabel{j} -- (-0.3,0.3);
\end{tikzpicture}\ $,
which are defined by \cref{rightpivot} and the properties itemized after \cref{leftpivots}.
Then we can introduce the bubble generating functions
$\ \begin{tikzpicture}[KM,baseline=-1mm,scale=.9]
\draw[to-] (-0.25,0) arc(180:-180:0.25);
\node at (0,-.4) {\strandlabel{i}};
\region{1,0}{\lambda};
\node at (.58,0) {$(u)$};
\end{tikzpicture}$
and
$\ \begin{tikzpicture}[KM,baseline=-1mm,scale=.9]
\draw[-to] (-0.25,0) arc(180:-180:0.25);
\node at (0,-.4) {\strandlabel{i}};
\region{1,0}{\lambda};
\node at (.58,0) {$(u)$};
\end{tikzpicture}$, which are the formal power series in 
$Z(\catR_\lambda)\lround u^{-1}\rround$
with leading terms $u^{h_i(\lambda)}$
and $u^{-h_i(\lambda)}$, respectively,
satisfying \cref{miles,infgrass}.

A shortcoming of \cref{maintheorem} is that it gives no clue as to exactly 
how to express the leftward Kac-Moody caps, cups and crossings 
in terms of the natural transformations arising from the initially given Heisenberg action on $\catR$.
We are going to fill this gap in this section.
We find it easier to go in the other direction, so we focus initially on the problem of expressing the natural transformations defined by the leftward Heisenberg caps, cups and crossings in terms of the ones defined by the leftward Kac-Moody ones.
Some preparation is needed before we can write this down.

The definition of the Kac-Moody bubble generating functions gives that
\begin{align}
\begin{tikzpicture}[KM,baseline=-1mm,scale=.9]
\draw[to-] (-0.25,0) arc(180:-180:0.25);
\node at (0,-.4) {\strandlabel{i}};
\region{0.95,0}{\lambda};
\node at (.55,0) {$(u)$};
\end{tikzpicture}
&=
\begin{tikzpicture}[KM,baseline=-1mm]
\filledanticlockwisebubble{0,0}{u};
\bubblelabel{0,0}{u};
\node at (0,-.35) {\strandlabel{i}};
\region{0.4,0}{\lambda};
\end{tikzpicture}
+
\begin{tikzpicture}[KM,baseline=-1mm,scale=.9]
\draw[to-] (-0.25,0) arc(180:-180:0.25);
\node at (0,-.4) {\strandlabel{i}};
\circled{.25,0}{u};
\region{0.65,0}{\lambda};
\end{tikzpicture},&
\begin{tikzpicture}[KM,baseline=-1mm,scale=.9]
\draw[to-] (-0.25,0) arc(-180:180:0.25);
\node at (0,-.4) {\strandlabel{i}};
\region{0.95,0}{\lambda};
\node at (.55,0) {$(u)$};
\end{tikzpicture}
&=
\begin{tikzpicture}[KM,baseline=-1mm]
\filledclockwisebubble{0,0}{u};
\node at (0,-.35) {\strandlabel{i}};
\region{0.4,0}{\lambda};
\bubblelabel{0,0}{u};
\end{tikzpicture}
+
\begin{tikzpicture}[KM,baseline=-1mm,scale=.9]
\draw[to-] (0.25,0) arc(0:360:0.25);
\node at (0,-.4) {\strandlabel{i}};
\circled{-.25,0}{u};
\region{0.55,0}{\lambda};
\end{tikzpicture},\label{bgf}
\end{align}
with
$\begin{tikzpicture}[KM,baseline=-1mm]
\filledanticlockwisebubble{0,0}{u};
\bubblelabel{0,0}{u};
\node at (0,-.35) {\strandlabel{i}};
\region{0.4,0}{\lambda};
\end{tikzpicture}$ (resp., 
$\begin{tikzpicture}[KM,baseline=-1mm]
\filledclockwisebubble{0,0};
\bubblelabel{0,0}{u};
\node at (0,-.35) {\strandlabel{i}};
\region{0.4,0}{\lambda};
\end{tikzpicture}$) being a monic polynomial of degree $h_i(\lambda)$ 
(resp., $-h_i(\lambda)$) if this is a non-negative integer, or
$0$ otherwise.
We refer to these as the {\em fake bubble polynomials}.
Replacing $u$ by $u-i$ gives the equivalent identities
\begin{align}\label{sbgf}
\begin{tikzpicture}[KM,baseline=-1mm,scale=.9]
\draw[to-] (-0.25,0) arc(180:-180:0.25);
\node at (0,-.4) {\strandlabel{i}};
\region{1.65,0}{\lambda};
\node at (.9,0) {$(u-i)$};
\end{tikzpicture}
&=
\begin{tikzpicture}[KM,baseline=-1mm]
\stretchedanticlockwisebubble{0,0}{u};
\bubblelabel{0,0}{u-i};
\node at (0,-.35) {\strandlabel{i}};
\region{0.5,0}{\lambda};
\end{tikzpicture}
+
\begin{tikzpicture}[KM,baseline=-1mm,scale=.9]
\draw[to-] (-0.25,0) arc(180:-180:0.25);
\node at (0,-.4) {\strandlabel{i}};
\region{0.65,0}{\lambda};
\pin{.25,0}{1.5,0}{(u-x-i)^{-1}};
\region{2.65,0}{\lambda};
\end{tikzpicture},&
\begin{tikzpicture}[KM,baseline=-1mm,scale=.9]
\draw[to-] (-0.25,0) arc(-180:180:0.25);
\node at (0,-.4) {\strandlabel{i}};
\region{1.65,0}{\lambda};
\node at (.9,0) {$(u-i)$};
\end{tikzpicture}
&=
\begin{tikzpicture}[KM,baseline=-1mm]
\stretchedclockwisebubble{0,0}{u};
\bubblelabel{0,0}{u-i};
\node at (0,-.35) {\strandlabel{i}};
\region{0.5,0}{\lambda};
\end{tikzpicture}
+
\begin{tikzpicture}[KM,baseline=-1mm,scale=.9]
\draw[to-] (0.25,0) arc(0:360:0.25);
\node at (0,-.4) {\strandlabel{i}};
\pin{-.25,0}{-1.5,0}{(u-x-i)^{-1}};
\region{.55,0}{\lambda};
\end{tikzpicture}.
\end{align}
These are decompositions of formal Laurent series in $u^{-1}$ as the sum of its polynomial and non-polynomial parts.
It is still the case that
the {\em shifted fake bubble polynomial}
$\begin{tikzpicture}[KM,baseline=-1mm]
\stretchedanticlockwisebubble{0,0}{u};
\bubblelabel{0,0}{u-i};
\node at (0,-.35) {\strandlabel{i}};
\region{0.5,0}{\lambda};
\end{tikzpicture}$ (resp., 
$\begin{tikzpicture}[KM,baseline=-1mm]
\stretchedclockwisebubble{0,0};
\bubblelabel{0,0}{u-i};
\node at (0,-.35) {\strandlabel{i}};
\region{0.5,0}{\lambda};
\end{tikzpicture}$) is monic of degree $h_i(\lambda)$ 
(resp., $-h_i(\lambda)$) if this is a non-negative integer, or
$0$ otherwise.

Now comes the main new definition.
For $i \neq j$ in $I$, we define the counterclockwise and clockwise {\em internal bubbles of color $j$} to be the natural transformations
\begin{align}
\begin{tikzpicture}[KM,centerzero]
\draw[-to] (0,-.5)\botlabel{i} to (0,.5);
\anticlockwiseinternalbubbleR[j]{0,0};
\region{0.4,0}{\lambda};
\end{tikzpicture}
&:=
\begin{tikzpicture}[KM,centerzero]
\draw[-to] (0,-.5)\botlabel{i} to (0,.5);
\draw[to-] (0.7,0) arc(360:0:0.2);
\pinpin{0.3,0}{0,0}{-1.3,0}{(i-j+x-y)^{-1}};
\strand{.5,-.38}{j};
\region{0.95,0}{\lambda};
\end{tikzpicture}
+
\left[\ 
\begin{tikzpicture}[KM,centerzero]
\draw[-to] (0,-.5)\botlabel{i} to (0,.5);
\pin{0,0}{-1.1,0}{(u-x-i)^{-1}};
\strand{.6,-.38}{j};
\stretchedanticlockwisebubble{.6,0};
\bubblelabel{.6,0}{u-j};
\region{1.2,0}{\lambda};
\end{tikzpicture} 
\right]_{u:-1}\!\!\!,\label{int1}\\
\begin{tikzpicture}[KM,centerzero]
\draw[-to] (0,-.5)\botlabel{i} to (0,.5);
\clockwiseinternalbubbleL[j]{0,0};
\region{-.4,0}{\lambda};
\end{tikzpicture}
&:=
\begin{tikzpicture}[KM,centerzero]
\draw[-to] (0,-.5)\botlabel{i} to (0,.5);
\draw[to-] (-0.7,0) arc(-180:180:0.2);
\pinpin{-0.3,0}{0,0}{1.3,0}{(i-j+y-x)^{-1}};
\strand{-.5,-.38}{j};
\region{-.95,0}{\lambda};
\end{tikzpicture}
+
\left[
\begin{tikzpicture}[KM,centerzero]
\draw[-to] (-1,-.5)\botlabel{i} to (-1,.5);
\pin{-1,0}{.1,0}{(u-x-i)^{-1}};
\stretchedclockwisebubble{-1.6,0};
\bubblelabel{-1.6,0}{u-j};
\strand{-1.6,-.38}{j};
\region{-2.15,0}{\lambda};
\end{tikzpicture}\  
\right]_{u:-1}\!\!\!.\label{int2}
\end{align}
The first terms on the right hand sides of \cref{int1,int2} make sense because $i-j \in \kk^\times$ and $x$ and $y$ are locally nilpotent. 
If $h_j(\lambda) = 0$, as is the case for all but finitely many $j \in I$, then the shifted fake bubble polynomials 
appearing in the second terms on the right hand sides 
are identities, hence, these terms are simply equal to
$\begin{tikzpicture}[KM,anchorbase]
\draw[-to] (0,-0.2)\botlabel{i} -- (0,0.2);
\region{0.2,0}{\lambda};
\end{tikzpicture}$
 or $\begin{tikzpicture}[KM,anchorbase]
\draw[-to] (0,-0.2)\botlabel{i} -- (0,0.2);
\region{-0.2,0}{\lambda};
\end{tikzpicture}$, respectively.
Moreover, still with $h_j(\lambda) = 0$, all of the bubbles with a non-negative number of dots appearing in the expansions of the first terms 
on the right hand sides of \cref{int1,int2} are locally nilpotent.
This follows because all
such bubbles act as zero on irreducible objects of $\catR_\lambda$
according to property \ref{noods5} from \cref{s2-KM}.
Consequently, for any indecomposable object $V \in \ob\catR_\lambda$
or $\ob\catR_{\lambda-\alpha_i}$, respectively, and all but finitely many $i \neq j \in I$, the internal bubbles of color $j$ from \cref{int1,int2}
define elements of $\End_{\catR}(E_i V)$ of the form $1+z$
for $z$ in the Jacobson radical of this finite-dimensional algebra.
Thus, it makes sense to take the (possibly infinite) product over all $i \neq j \in I$ of these commuting natural transformations, defining the counterclockwise and clockwise {\em internal bubbles} to be the vertical compositions
\begin{align}\label{int3}
\begin{tikzpicture}[KM,centerzero]
\draw[-to] (0,-.5)\botlabel{i} to (0,.5);
\anticlockwiseinternalbubbleR{0,0};
\region{0.4,0}{\lambda};
\end{tikzpicture}
&:= {\textstyle\prod_{i \neq j \in I}}\ 
\begin{tikzpicture}[KM,centerzero]
\draw[-to] (0,-.5)\botlabel{i} to (0,.5);
\anticlockwiseinternalbubbleR[j]{0,0};
\region{0.4,0}{\lambda};
\end{tikzpicture}\ ,&
\begin{tikzpicture}[KM,centerzero]
\draw[-to] (0,-.5)\botlabel{i} to (0,.5);
\clockwiseinternalbubbleL{0,0};
\region{-.4,0}{\lambda};
\end{tikzpicture}
&:=-
{\textstyle\prod_{i \neq j \in I}}
\begin{tikzpicture}[KM,centerzero]
\draw[-to] (0,-.5)\botlabel{i} to (0,.5);
\clockwiseinternalbubbleL[j]{0,0};
\region{-.4,0}{\lambda};
\end{tikzpicture}\ .
\end{align}
Opening some parentheses and using \cref{trickconsequence} 
gives the following more explicit formulae:
\begin{align}\label{daunting1}
\begin{tikzpicture}[KM,centerzero]
\draw[-to] (0,-.5)\botlabel{i} to (0,.5);
\anticlockwiseinternalbubbleR{0,0};
\region{0.4,0}{\lambda};
\end{tikzpicture}
&=
\sum_{\substack{J, K\text{ with }|J|<\infty\\J \sqcup K=I-\{i\}}}
\left[\ 
\begin{tikzpicture}[KM,centerzero]
\draw[-to] (0,-.6)\botlabel{i} to (0,.5);
\node at (-2.8,.2) {$\prod_{j \in J}$};
\draw[to-] (.7,.2) arc(360:0:0.2);
\strand{.5,-.18}{j};
\pinpin{0.3,.2}{0,.2}{-1.3,.2}{(i-j+x-y)^{-1}};
\node at (1.55,.2) {$\prod_{k \in K}$};
\stretchedanticlockwisebubble{2.45,.2};
\bubblelabel{2.45,.2}{u-k};
\strand{2.45,-.15}{k};
\pin{0,-.4}{-1.1,-.4}{(u-x-i)^{-1}};
\region{1.2,-.35}{\lambda};
\end{tikzpicture} 
\right]_{u:-1}\!\!\!,\\\label{daunting2}
\begin{tikzpicture}[KM,centerzero]
\draw[-to] (0,-.5)\botlabel{i} to (0,.5);
\clockwiseinternalbubbleL{0,0};
\region{-0.4,0}{\lambda};
\end{tikzpicture}\ \ 
&=
-\sum_{\substack{J, K\text{ with }|J|<\infty\\J \sqcup K=I-\{i\}}}
\left[\ 
\begin{tikzpicture}[KM,centerzero]
\draw[-to] (0,-.6)\botlabel{i} to (0,.5);
\node at (-1.25,.2) {$\prod_{j \in J}$};
\draw[to-] (-.7,.2) arc(-180:180:0.2);
\strand{-.5,-.18}{j};
\pinpin{-0.3,.2}{0,.2}{1.3,.2}{(i-j+y-x)^{-1}};
\node at (-3.3,.2) {$\prod_{k \in K}$};
\stretchedclockwisebubble{-2.4,.2};
\bubblelabel{-2.4,.2}{u-k};
\strand{-2.4,-.15}{k};
\pin{0,-.4}{1.1,-.4}{(u-x-i)^{-1}};
\region{-1.55,-.35}{\lambda};
\end{tikzpicture} 
\right]_{u:-1}\!\!\!.
\end{align}
In these expressions, the product over $J$ applies to the double pin at the top (the pin labelled $(u-x-i)^{-1}$ is only taken once), and the product over $K$ applies to the shifted fake bubble polynomial.
The (possibly infinite) product over $K$ makes sense because all but finitely many of the shifted fake bubble polynomials are 1.
The (possibly infinite) sum at the front makes sense because, on any given indecomposable $V$, the expression inside 
is 0 for all but finitely many $J$.

We also introduce counterclockwise and clockwise internal bubbles 
on downward strings
\begin{align}\label{int5}
\begin{tikzpicture}[KM,centerzero]
\draw[to-] (0,-.5)\botlabel{i} to (0,.5);
\anticlockwiseinternalbubbleL{0,0};
\region{-0.4,0}{\lambda};
\end{tikzpicture}
&= {\textstyle\prod_{i \neq j \in I}}
\begin{tikzpicture}[KM,centerzero]
\draw[to-] (0,-.5)\botlabel{i} to (0,.5);
\anticlockwiseinternalbubbleL[j]{0,0};
\region{-0.4,0}{\lambda};
\end{tikzpicture}\ ,&
\begin{tikzpicture}[KM,centerzero]
\draw[to-] (0,-.5)\botlabel{i} to (0,.5);
\clockwiseinternalbubbleR{0,0};
\region{.4,0}{\lambda};
\end{tikzpicture}
&=
-{\textstyle\prod_{i \neq j \in I}}\ 
\begin{tikzpicture}[KM,centerzero]
\draw[to-] (0,-.5)\botlabel{i} to (0,.5);
\clockwiseinternalbubbleR[j]{0,0};
\region{.4,0}{\lambda};
\end{tikzpicture}\ ,
\end{align}
whose 
definitions are obtained from the ones in the previous paragraph by rotating all diagrams through 180$^\circ$. This ensures that all sorts of internal bubbles slide across both rightward and leftward caps and cups.

\begin{theo}\label{average}
The leftward Heisenberg caps and cups are related to the leftward Kac-Moody caps and cups by 
\begin{align}\label{leftwardbuggers}
\begin{tikzpicture}[KM,anchorbase,scale=1.5]
\draw (0.3,-0.4)\botlabel{j} to (0.3,-.2);
\draw[-to] (-0.1,-.2) to (-0.1,-0.4)\botlabel{i};
\draw[H] (.3,-.2) to[out=90,in=90,looseness=1.5] (-.1,-.2);
\notch{-.1,-.2};
\notch{.3,-.2};
\region{.5,-.3}{\lambda};
\end{tikzpicture}
&=
\delta_{i,j}\ 
\begin{tikzpicture}[KM,anchorbase,scale=1.5]
\draw (0.3,-0.4) to (0.3,-.2);
\draw[-to] (-0.1,-.2) to (-0.1,-0.4)\botlabel{i};
\draw (.3,-.2) to[out=90,in=90,looseness=2] (-.1,-.2);
\anticlockwiseinternalbubbleR{.27,-0.18};
\region{.6,-.2}{\lambda};
\end{tikzpicture}\ ,&
\begin{tikzpicture}[KM,anchorbase,scale=1.5]
\draw[-] (0.3,0.4)\toplabel{j} to (0.3,.2);
\draw[-to] (-0.1,.2) to (-0.1,0.4)\toplabel{i};
\draw[H] (.3,.2) to[out=-90,in=-90,looseness=1.5] (-.1,.2);
\notch{-.1,.2};
\notch{.3,.2};
\region{.5,.3}{\lambda};
\end{tikzpicture}
&= \delta_{i,j}\ 
\begin{tikzpicture}[KM,anchorbase,scale=1.5]
\draw[-] (0.3,0.4) to (0.3,.2);
\draw[-to] (-0.1,.2) to (-0.1,0.4)\toplabel{i};
\draw (.3,.2) to[out=-90,in=-90,looseness=2] (-.1,.2);
\clockwiseinternalbubbleL{-.07,0.18};
\region{.5,.2}{\lambda};
\end{tikzpicture}\end{align}
for $i,j \in I$ and $\lambda \in X$.
Hence, the Heisenberg bubble generating functions are related to the Kac-Moody bubble generating functions by
\begin{align}
\label{newbub}
\begin{tikzpicture}[H,baseline=-1mm]
\draw[to-] (-0.25,0) arc(180:-180:0.25);
\node at (.53,0) {\color{black}$(u)$};
\region{1,0}{\lambda};
\end{tikzpicture}
&=
{\textstyle\prod_{i \in I}}\ 
\begin{tikzpicture}[KM,baseline=-1mm]
\draw[to-] (-0.25,0) arc(180:-180:0.25);
\node at (.82,0) {\color{black}$(u-i)$};
\region{1.55,0}{\lambda};
\strand{0,-.4}{i};
\end{tikzpicture},&
\begin{tikzpicture}[H,baseline=-1mm]
\draw[-to] (-0.25,0) arc(180:-180:0.25);
\node at (.53,0) {$\color{black}(u)$};
\region{1,0}{\lambda};
\end{tikzpicture}
&=-
{\textstyle\prod_{i \in I}}\ 
\begin{tikzpicture}[KM,baseline=-1mm]
\draw[-to] (-0.25,0) arc(180:-180:0.25);
\node at (.82,0) {$\color{black}(u-i)$};
\strand{0,-.4}{i};
\region{1.55,0}{\lambda};
\end{tikzpicture}.
\end{align}
\end{theo}

\begin{proof}
We prove this just in the case that $\kappa \leq 0$; an analogous
argument treats $\kappa > 0$. 
All of the equations to be proved are 
trivial when $\catR_\lambda$ is zero, so we assume from
now on that we are given some fixed $\lambda \in X$ such that $\catR_\lambda$ is non-zero.
We will use the shorthand $\beta_{i,r}$ to denote coefficients of the shifted fake bubble polynomials on the Kac-Moody side:
\begin{equation}\label{how}
\begin{tikzpicture}[KM,baseline=-1mm]
\stretchedanticlockwisebubble{0,0}{u};
\bubblelabel{0,0}{u-i};
\node at (0,-.35) {\strandlabel{i}};
\region{0.5,0}{\lambda};
\end{tikzpicture}=
\sum_{r=0}^{h_i(\lambda)}
\beta_{i,r} u^r.
\end{equation}
Assuming that $h_i(\lambda) \geq 0$ so that it is defined, 
we have that
$\beta_{i,h_i(\lambda)} = 1$.
We then have that
\begin{equation}\label{abitmore}
{\textstyle\prod_{k \in K}}
\begin{tikzpicture}[KM,baseline=-1mm]
\stretchedanticlockwisebubble{0,0}{u};
\bubblelabel{0,0}{u-k};
\node at (0,-.35) {\strandlabel{k}};
\region{0.5,0}{\lambda};
\end{tikzpicture}=
\sum_{\substack{(r_k)_{k \in K}\\0 \leq r_k \leq h_k(\lambda)}}
u^{\sum_k r_k}
{\textstyle\prod_{k \in K}} \beta_{k,r_k}.
\end{equation}
We will use this several times below to rewrite the fake bubble polynomials appearing in \cref{daunting1}.

Now we define two new natural transformations denoted
$\ \begin{tikzpicture}[H,centerzero,scale=.8]
\draw[-,shadow] (-0.25,-0.2) to [out=90,in=90,looseness=3](0.25,-0.2);
\draw[-to] (-0.25,-0.22) to [out=-90,in=90,looseness=3](-0.25,-0.23);
\end{tikzpicture}\ $
and
$\ \begin{tikzpicture}[H,centerzero,scale=.8]
\draw[-,shadow] (-0.25,0.2)  to [out=-90,in=-90,looseness=3](0.25,0.2);
\draw[to-] (-0.25,0.22) to [out=-90,in=90,looseness=3](-0.25,0.23);
\end{tikzpicture}\ $ on the Heisenberg side so that they satisfy \cref{leftwardbuggers}, that is,
\begin{align}\label{again}
\begin{tikzpicture}[KM,anchorbase,scale=1.5]
\draw (0.3,-0.4)\botlabel{j} to (0.3,-.3);
\draw[-to] (-0.1,-.3) to (-0.1,-0.4)\botlabel{i};
\draw[H,shadow] (.3,-.2) to[out=90,in=90,looseness=2] (-.1,-.2);
%\draw[H,-to] (-0.1,-0.22) to [out=-90,in=90](-0.1,-0.23);
\draw[H](-.1,-.2) to (-.1,-.3);\draw[H](.3,-.2) to (.3,-.3);
\notch{-.1,-.3};
\notch{.3,-.3};
\region{.5,-.3}{\lambda};
\end{tikzpicture}
&:=
\delta_{i,j}\ 
\begin{tikzpicture}[KM,anchorbase,scale=1.5]
\draw (0.3,-0.4) to (0.3,-.2);
\draw[-to] (-0.1,-.2) to (-0.1,-0.4)\botlabel{i};
\draw (.3,-.2) to[out=90,in=90,looseness=2] (-.1,-.2);
\anticlockwiseinternalbubbleR{.27,-0.18};
\region{.6,-.2}{\lambda};
\end{tikzpicture}\ ,&
\begin{tikzpicture}[KM,anchorbase,scale=1.5]
\draw (0.3,0.4)\toplabel{j} to (0.3,.3);
\draw[-to] (-0.1,.3) to (-0.1,0.4)\toplabel{i};
\draw[H,shadow] (.3,.2) to[out=-90,in=-90,looseness=2] (-.1,.2);
%\draw[H,-to] (-0.1,0.22) to [out=90,in=-90](-0.1,0.23);
\draw[H](-.1,.2) to (-.1,.3);\draw[H](.3,.2) to (.3,.3);
\notch{-.1,.3};
\notch{.3,.3};
\region{.5,.3}{\lambda};
\end{tikzpicture}
&= \delta_{i,j}\ 
\begin{tikzpicture}[KM,anchorbase,scale=1.5]
\draw[-] (0.3,0.4) to (0.3,.2);
\draw[-to] (-0.1,.2) to (-0.1,0.4)\toplabel{i};
\draw (.3,.2) to[out=-90,in=-90,looseness=2] (-.1,.2);
\clockwiseinternalbubbleL{-.07,0.18};
\region{.5,.2}{\lambda};
\end{tikzpicture}
\end{align}
for $i,j \in I$.
Using this new version of the leftward Heisenberg cap, 
we obtain a new version of the counterclockwise 
Heisenberg bubble generating function, denoted
$\ \begin{tikzpicture}[H,baseline=-1mm,scale=.9]
\draw[-,shadow] (-0.25,0) arc(180:0:0.25);
\draw[-to] (-0.25,-0.02) to [out=-90,in=90,looseness=3](-0.25,-0.03);
\draw[-] (-0.25,0) arc(-180:0:0.25);
\node at (.55,0) {\color{black}$(u)$};
\region{1,0}{\lambda};
\end{tikzpicture}$. It is defined like in \cref{Hmiles} 
so that it is a 
formal Laurent series in 
$Z(\catR)\llbracket u^{-1}\rrbracket$ with leading term $u^\kappa$ 
and
$\left[\ \begin{tikzpicture}[H,baseline=-1mm,scale=.9]
\draw[-,shadow] (-0.25,0) arc(180:0:0.25);
\draw[-to] (-0.25,-0.02) to [out=-90,in=90,looseness=3](-0.25,-0.03);
\draw[-] (0.25,0) arc(0:-180:0.25);
\node at (.55,0) {\color{black}$(u)$};
\region{1,0}{\lambda};
\end{tikzpicture}\right]_{u:<0} = \begin{tikzpicture}[H,baseline=-1mm]
\draw[-,shadow] (-0.25,0) arc(180:0:0.25);
\draw[-to] (-0.25,-0.02) to [out=-90,in=90,looseness=3](-0.25,-0.03);
\draw[-] (0.25,0) arc(0:-180:0.25);
\circled{.18,-.18}{u};
\end{tikzpicture}$.
Now we proceed to prove a series of claims, which together prove the theorem.

\vspace{2mm}
\noindent
\underline{Claim 1}. {\em 
$\begin{tikzpicture}[H,baseline=-1mm,scale=.9]
\draw[-,shadow] (-0.25,0) arc(180:0:0.25);
\draw[-to] (-0.25,-0.02) to [out=-90,in=90,looseness=3](-0.25,-0.03);
\draw[-] (0.25,0) arc(0:-180:0.25);
\node at (.55,0) {\color{black}$(u)$};
\region{1,0}{\lambda};
\end{tikzpicture}=
{\textstyle\prod_{i \in I}}\ 
\begin{tikzpicture}[KM,baseline=-1mm]
\draw[to-] (-0.25,0) arc(180:-180:0.25);
\node at (.82,0) {\color{black}$(u-i)$};
\region{1.55,0}{\lambda};
\strand{0,-.4}{i};
\end{tikzpicture}$.}

\noindent
\underline{Proof}.
To prove this, we will use the following fancy notation.
Suppose we are given a finite subset $J \subseteq I$.
Let $A_J$ be the polynomial algebra $\kk[x_j\:|\:j \in J]$.
There is a linear map (usually {\em not} an algebra homomorphism!)
$\theta_J:A_J \rightarrow Z(\catR_\lambda)$ 
mapping the polynomial $f(x_j\:|\:j \in J)$ to the natural transformation obtained by pinning
$f(x_j+j\:|\:j \in J)$ to 
the string diagram $\prod_{j \in J} 
\begin{tikzpicture}[KM,baseline=-1mm,scale=.8]
\draw[to-] (-0.25,0) arc(180:-180:0.25);
\node at (0,-.4) {\strandlabel{j}};
\opendot{.25,0};
\region{.6,0}{\lambda};
\end{tikzpicture}$ 
with the variable $x_j$ corresponding to the dot on the $j$th bubble.
Let $\widehat{A}_J$ be
the completion of $A_J$ at the maximal
ideal $(x_j - j\:|\:j \in J)$.
Since the Kac-Moody dots are locally nilpotent, 
$\theta_J$ naturally extends 
to a linear map $\widehat\theta_J:\widehat{A}_J \rightarrow
Z(\catR_\lambda)$. For example, we have that
\begin{equation}\label{why} 
\widehat\theta_{\{i,j\}}
\left( \frac{1}{x_i-x_j}\right) = \begin{tikzpicture}[KM,baseline=-1mm,scale=.9]
\draw[to-] (-0.25,0) arc(180:-180:0.25);
\node at (0,-.4) {\strandlabel{i}};
\draw[to-] (1.25,0) arc(0:-360:0.25);
\node at (1,-.4) {\strandlabel{j}};
\pinpin{.75,0}{.25,0}{.25,.75}{(i-j+x-y)^{-1}};
\region{1.6,0}{\lambda};
\end{tikzpicture}
\end{equation}
for $i \neq j$ in $I$.
To prove the claim, as $\kappa \leq 0$, 
the left hand side of the equation to be proved is equal to
$\delta_{\kappa,0}+
\begin{tikzpicture}[H,baseline=-1mm,scale=.9]
\draw[-,shadow] (-0.25,0) arc(180:0:0.25);
\draw[-to] (-0.25,-0.02) to [out=-90,in=90,looseness=3](-0.25,-0.03);
\draw[-] (0.25,0) arc(0:-180:0.25);
\circled{.18,-.18}{u};
\region{.6,0}{\lambda};
\end{tikzpicture}$. As 
$\sum_{i \in I} h_i(\lambda) = \kappa$ by \cref{dunking}
and each
$\begin{tikzpicture}[KM,baseline=-1mm,scale=.8]
\draw[to-] (-0.25,0) arc(180:-180:0.25);
\node at (1,0) {\color{black}$(u-i)$};
\region{1.85,0}{\lambda};
\strand{0,-.4}{i};
\end{tikzpicture}$ has leading term $u^{h_i(\lambda)}$,  
the right hand side is also equal to $\delta_{\kappa,0}$ plus a term in
$u^{-1} Z(\catR_\lambda)\llbracket u^{-1}\rrbracket$.
Thus, we are reduced to proving that
\begin{equation}\label{todo}
\begin{tikzpicture}[H,baseline=-1mm]
\draw[-,shadow] (-0.25,0) arc(180:0:0.25);
\draw[-to] (-0.25,-0.02) to [out=-90,in=90,looseness=3](-0.25,-0.03);
\draw[-] (0.25,0) arc(0:-180:0.25);
\circled{.18,-.18}{u};
\region{.6,0}{\lambda};
\end{tikzpicture}
=
\left[\ {\textstyle\prod_{i \in I}}\ 
\begin{tikzpicture}[KM,baseline=-1mm]
\draw[to-] (-0.25,0) arc(180:-180:0.25);
\node at (.83,0) {\color{black}$(u-i)$};
\region{1.55,0}{\lambda};
\strand{0,-.4}{i};
\end{tikzpicture}\right]_{u:<0}.
\end{equation}
Using \cref{sbgf} and opening some parentheses, the right hand side 
of \cref{todo} is equal to
$$
\sum_{\substack{J, K\text{ with }|J|<\infty\\
J \sqcup K=I}}
\left[\left({\textstyle\prod_{j \in J}}\ 
\begin{tikzpicture}[KM,baseline=-1mm,scale=.9]
\draw[to-] (-0.25,0) arc(180:-180:0.25);
\node at (0,-.4) {\strandlabel{j}};
\region{0.65,0}{\lambda};
\pin{.25,0}{1.5,0}{(u-x-j)^{-1}};
\region{2.65,0}{\lambda};
\end{tikzpicture}\right)\left( {\textstyle\prod_{k \in K}}
\begin{tikzpicture}[KM,baseline=-1mm]
\stretchedanticlockwisebubble{0,0}{u};
\bubblelabel{0,0}{u-k};
\node at (0,-.35) {\strandlabel{k}};
\region{0.5,0}{\lambda};
\end{tikzpicture}\right)
\right]_{u:<0}.
$$
Using the fancy notation and \cref{abitmore}, this is equal to
\begin{equation}
\sum_{\substack{J, K\text{ with }|J|<\infty\\
J \sqcup K=I}}
\sum_{\substack{(r_k)_{k \in K}\\0 \leq r_k \leq h_k(\lambda)}}
\widehat\theta_J \left(
\left[u^{\sum_k r_k} {\textstyle\prod_{j \in J}}
\frac{1}{u-x_j}
\right]_{u:<0}
\right)
{\textstyle\prod_{k \in K}} \beta_{k,r_k}.\label{couples}
\end{equation}
By partial fractions, we have in $\widehat A_J\llbracket u^{-1}\rrbracket$ that
$$
\left[u^r {\textstyle\prod_{j \in J}}
\frac{1}{u-x_j}
\right]_{u:<0} =
\left[
\sum_{i \in J} \frac{u^r}{u-x_i} {\textstyle\prod_{i \neq j \in J}}
\frac{1}{x_i-x_j}
\right]_{u:<0}
\stackrel{\cref{trick}}{=}\ 
\sum_{i \in J} \frac{x_i^{r}}{u-x_i} \prod_{i \neq j \in J}
\frac{1}{x_i-x_j}.
$$
Using this identity, we can simplify \cref{couples}, that is, the right hand side of \cref{todo}, to obtain
\begin{equation}\label{todo2}
\sum_{\substack{J, K\text{ with }0 < |J|<\infty\\
J \sqcup K=I}}
\sum_{\substack{(r_k)_{k \in K}\\0 \leq r_k \leq h_k(\lambda)}}
\widehat\theta_J
\left(\sum_{i \in J} 
\frac{x_i^{\sum_k r_k}}{(u-x_i) \prod_{i \neq j \in J}(x_i-x_j)}
\right) {\textstyle\prod_{k \in K}} \beta_{k,r_k}.
\end{equation}
Now we turn our attention to the left hand side of \cref{todo}.
Using \cref{gofigure,prodots,sliders}
then the definition \cref{leftwardbuggers}, we get that
$$
\begin{tikzpicture}[H,baseline=-1mm]
\draw[-,shadow] (-0.25,0) arc(180:0:0.25);
\draw[-to] (-0.25,-0.02) to [out=-90,in=90,looseness=3](-0.25,-0.03);
\draw[-] (0.25,0) arc(0:-180:0.25);
\circled{.18,-.18}{u};
\region{.6,0}{\lambda};
\end{tikzpicture}
=
\sum_{i \in I}
\begin{tikzpicture}[KM,baseline=-1mm]
\draw[H,-,shadow] (-0.25,0.2) arc(180:0:0.25);
\draw[H,-to] (-0.25,0.18) to [out=-90,in=90,looseness=3](-0.25,0.17);
\draw[H,-] (-.25,.2) to (-.25,.1);
\draw[H] (.25,.2) to (.25,.1);
\draw[-] (0.25,0) arc(0:-180:0.25);
\draw[-] (-.25,.1) to (-.25,0); 
\draw (.25,.1) to (.25,0);
\notch{-.25,.1};
\notch{.25,.1};
\strand{0,-.4}{i};
\pin{.25,0}{1.4,0}{(u-x-i)^{-1}};
\region{2.5,0}{\lambda};
\end{tikzpicture}
=\sum_{i \in I}
\begin{tikzpicture}[KM,baseline=-1mm]
\draw[-] (-0.25,0.2) arc(180:0:0.25);
\draw[-] (0.25,0) arc(0:-180:0.25);
\draw[-to] (-.25,.2) to (-.25,0); 
\draw (.25,.2) to (.25,0);
\anticlockwiseinternalbubbleL{.25,.25};
\strand{0,-.4}{i};
\pin{.25,-.1}{1.4,-.1}{(u-x-i)^{-1}};
\region{2.5,0}{\lambda};
\end{tikzpicture}\ .
$$
Applying \cref{daunting1} (replacing $u$ there with $v$), this equals
$$
\sum_{i \in I}
\sum_{\substack{J, K\text{ with }|J|<\infty\\J \sqcup K=I-\{i\}}}
\left[\ 
\begin{tikzpicture}[KM,centerzero]
\draw[-to] (-3.5,0) to (-3.5,-.6) to [out=-90,in=-90,looseness=.9] (0,-.6)
to (0,.4) to [out=90,in=90,looseness=.9] (-3.5,.4)
to (-3.5,0);
\node at (-2.8,.25) {$\prod_{j \in J}$};
\draw[to-] (.7,.25) arc(360:0:0.2);
\strand{.5,-.08}{j};
\pinpin{0.3,.25}{0,.25}{-1.3,.25}{(i-j+x-y)^{-1}};
\node at (1.55,.25) {$\prod_{k \in K}$};
\stretchedanticlockwisebubble{2.45,.25};
\bubblelabel{2.45,.25}{v-k};
\strand{2.45,-.1}{k};
\pin{0,-.3}{-1.1,-.3}{(v-x-i)^{-1}};
\pin{0,-.6}{1.2,-.6}{(u-x-i)^{-1}};
\region{2.4,-.6}{\lambda};
\strand{-.5,-1.4}{i};
\end{tikzpicture} 
\right]_{v:-1}\!\!\!.
$$
Switching to the fancy notation using a variation on \cref{why} and \cref{how}, this is
$$
\sum_{i \in I}
\sum_{\substack{J, K\text{ with }|J|<\infty\\J \sqcup K=I-\{i\}}}
\sum_{\substack{(r_k)_{k \in K}\\0 \leq r_k \leq h_k(\lambda)}}
\widehat\theta_{J \cup \{i\}}
\left(\left[
    \frac{v^{\sum_k r_k}}{(u-x_i)(v-x_i)\prod_{j \in J} (x_i-x_j)} 
 \right]_{v:-1}\right) {\textstyle\prod_{k \in K}} \beta_{k,r_k}.
$$
Applying \cref{trick}, this is
$$
\sum_{i \in I}
\sum_{\substack{J, K\text{ with }|J|<\infty\\J \sqcup K=I-\{i\}}}
\sum_{\substack{(r_k)_{k \in K}\\0 \leq r_k \leq h_k(\lambda)}}
\widehat\theta_{J \cup \{i\}}
\left(
    \frac{x_i^{\sum_k r_k}}{(u-x_i)\prod_{j \in J} (x_i-x_j)} 
\right) {\textstyle\prod_{k \in K}} \beta_{k,r_k}.
$$
Finally, to see that this is equal to \cref{todo2}, one just needs to reindex the second summation, replacing $J$ by $J \sqcup\{i\}$,
then move the first summation over $i \in I$ to the inside.

\vspace{2mm}
\noindent
\underline{Claim 2}.
{\em 
$\begin{tikzpicture}[H,anchorbase,scale=1]
\draw[to-] (0,-0.5) to[out=90,in=-120] (0.08,0.05);
\draw[shadow] (0.08,.05) to [out=60,in=up,looseness=2] (0.6,0);
%\draw[-to] (.08,.03) to [out=-120,in=60] (.07,.028);
\draw (.6,0) to[looseness=2,out=down,in=-60] (.08,-.05) to[out=120,in=down] (0,0.5);
\end{tikzpicture}
=
\left[
\begin{tikzpicture}[H,anchorbase,scale=1]
\draw[to-] (-0.8,-0.5) -- (-0.8,0.5);
\circled{-0.8,0}{u};
\draw[-,shadow] (0,0) arc(0:180:0.2);
\draw[-to,H] (-.4,0) to (-.4,-.02);
\draw[-] (-.4,0) arc(-180:0:0.2);
\node at (0.29,0) {$\color{black}(u)$};
\end{tikzpicture}
\!\right]_{u:-1}\!\!\!.$}

\noindent
\underline{Proof}.
We need to extend slightly
the fancy notation introduced in the proof of Claim 1.
Suppose we are given a finite subset $J \subseteq I$ and $i \in I$.
There is a linear map (usually {\em not} an algebra homomorphism!)
$\theta_{J,i}:A_J[x] \rightarrow \End(F_i|_{\catR_\lambda})$ 
mapping the polynomial $f(x,x_j\:|\:j \in J)$ to the natural transformation obtained by pinning
$f(x+i,x_j+j\:|\:j \in J)$ to $
\begin{tikzpicture}[KM,baseline=-1mm,scale=.8]
\draw[-to] (-.6,.3) to (-.6,-.3)\botlabel{i};
\opendot{-.6,0};
\region{-.2,0}{\lambda};
\end{tikzpicture}\prod_{j \in J} 
\begin{tikzpicture}[KM,baseline=-1mm,scale=.8]
\draw[to-] (-0.25,0) arc(180:-180:0.25);
\node at (0,-.4) {\strandlabel{j}};
\opendot{.25,0};
\end{tikzpicture}$ 
with $x$ corresponding to the dot on the propagating string and $x_j$ corresponding to the dot on the $j$th bubble.
Letting $\widehat{A_J[x]}$ be
the completion of $A_J[x]$ at the maximal
ideal $(x-i, x_j-j\:|\:j \in J)$, the map $\theta_{J,i}$ extends to
$\widehat\theta_{J,i}:\widehat{A_J[x]} \rightarrow
\End(F_i|_{\catR_\lambda})$.
Claim 2
follows if we can show that
\begin{equation}\label{beer}
\begin{tikzpicture}[KM,anchorbase,scale=1.1]
\draw[KM,to-](0,-.5) \botlabel{i} to (0,-.35);
\draw[KM](0,.5) \toplabel{i} to (0,.35);
\notch{0,-.35};\notch{0,.35};
\draw[H,-] (0,-0.35) to[out=90,in=-120] (0.08,0.05);
\draw[H,shadow] (0.08,.05) to [out=60,in=up,looseness=2] (0.6,0);
%\draw[H,-to] (.08,.03) to [out=-120,in=60] (.07,.028);
\draw[H] (.6,0) to[looseness=2,out=down,in=-60] (.08,-.05) to[out=120,in=down] (0,0.35);
\region{.85,0}{\lambda};
\end{tikzpicture}
=
\left[\ 
\begin{tikzpicture}[H,anchorbase,scale=1.1]
\draw[to-,KM](-.8,-.5) \botlabel{i} to (-.8,-.3);
\draw[KM](-.8,.5) \toplabel{i} to (-.8,.3);
\notch{-.8,-.3};\notch{-.8,.3};
\draw[-] (-0.8,-0.3) -- (-0.8,0.3);
\circled{-0.8,0}{u};
\draw[-,shadow] (0,0) arc(0:180:0.2);
\draw[-to,H] (-.4,0) to (-.4,-.02);
\draw[-] (-.4,0) arc(-180:0:0.2);
\node at (0.26,0) {$\color{black}(u)$};
\region{.7,0}{\lambda};
\end{tikzpicture}
\right]_{u:-1}\!\!\!.
\end{equation}
for all $i \in I$.
By the definition of Kac-Moody dots and Claim 1, the right hand side is equal to 
$$
\left[\ 
\begin{tikzpicture}[KM,anchorbase,scale=1.1]
\draw[to-] (-0.8,-0.5)\botlabel{i} -- (-0.8,0.5);
\pin{-0.8,0}{-1.8,0}{(u-x-i)^{-1}};
\node at (-.13,0) {$\prod_{j \in I}$};
\draw[-to] (.35,0) arc(-180:180:0.2);
\strand{.55,-.35}{j};
\node at (1.3,0) {$\color{black}(u-j)$};
\region{2,0}{\lambda};
\end{tikzpicture}
\right]_{u:-1}\!\!\!.
$$
Opening some parentheses,
using also the fancy notation and \cref{abitmore}, this expands as
\begin{equation*}
\sum_{\substack{J, K\text{ with }|J|<\infty\\
J \sqcup K=I}}
\sum_{\substack{(r_k)_{k \in K}\\0 \leq r_k \leq h_k(\lambda)}}
\theta_{J,i}\left(\left[
 \frac{u^{\sum_k r_k}}{(u-x)\prod_{j \in J}(u-x_j)}\right]_{u:-1}\right)
{\textstyle\prod_{k \in K}} \beta_{k,r_k}.
\end{equation*}
Using partial fractions, this decomposes as the sum of the following two expressions:
\begin{align}\label{rhs1}
\sum_{\substack{J, K\text{ with }|J|<\infty\\
J \sqcup K=I}}
\sum_{\substack{(r_k)_{k \in K}\\0 \leq r_k \leq h_k(\lambda)}}
\theta_{J,i}\left(
\frac{x^{\sum_k r_k}}{\prod_{j \in J}(x-x_j)}
\right)
{\textstyle\prod_{k \in K}} \beta_{k,r_k},\\\label{rhs2}
\sum_{\substack{J, K\text{ with }|J|<\infty\\
J \sqcup K=I}}
\sum_{\substack{(r_k)_{k \in K}\\0 \leq r_k \leq h_k(\lambda)}}
\theta_{J,i}\left(
\sum_{j \in J}
\frac{x_j^{\sum_k r_k}}{(x_j-x)\prod_{j\neq h \in J}(x_j-x_h)}
\right)
{\textstyle\prod_{k \in K}} \beta_{k,r_k}.
\end{align}
Now consider the left hand side of \cref{beer}. It is
\begin{align}\notag
\sum_{j \in I}
\begin{tikzpicture}[KM,anchorbase,scale=1.1]
\draw[KM,to-](0,-.5) \botlabel{i} to (0,-.35);
\draw[KM](0,.5) \toplabel{i} to (0,.35);
\notch{0,-.35};\notch{0,.35};
\draw[H,-] (0,-0.35) to[out=90,in=-120] (0.08,0.05);
\draw[H,shadow] (0.08,.05) to [out=60,in=up,looseness=2] (0.7,.05);
%\draw[H,-to] (.08,.05) to [out=-120,in=60] (.07,.029);
\draw[H] (.7,-0.05) to (.7,.05);
\draw (.7,-0.05) to (.7,-.25);
\draw[H] (.7,-.25) to[looseness=2,out=down,in=-60] (.08,-.05) to[out=120,in=down] (0,0.35);
\notch{.7,-0.05};\notch{.7,-.25};
\region{1.1,0}{\lambda};
\strand{.84,-0.13}{j};
\end{tikzpicture}
&=
\begin{tikzpicture}[KM,anchorbase,scale=1.1]
\region{.8,0}{\lambda};
\projcr{-0.06,0};
\draw[to-] (-0.2,-0.5)\botlabel{i} to[out=90,in=180] (0.3,0.23) to[out=0,in=up] (0.5,0) to[out=down,in=0] (0.3,-0.23) to[out=180,in=-90] (-0.2,0.5) \toplabel{i};
\anticlockwiseinternalbubbleR{.45,0};
\strand{.3,-0.34}{i};
\end{tikzpicture}
+\sum_{i \neq j \in I}
\begin{tikzpicture}[KM,anchorbase,scale=1.1]
\region{.8,0}{\lambda};
\projcr{-0.06,0};
\draw[to-] (-0.2,-0.5)\botlabel{i} to[out=90,in=180] (0.3,0.23) to[out=0,in=up] (0.5,0) to[out=down,in=0] (0.3,-0.23) to[out=180,in=-90] (-0.2,0.5) \toplabel{i};
\anticlockwiseinternalbubbleR{.45,0};
\strand{.3,-0.34}{j};
\end{tikzpicture}\\
&
\substack{\cref{billy}\\\textstyle=\\\cref{eggs}}\ 
\begin{tikzpicture}[KM,anchorbase,scale=1.1]
\draw[to-] (-0.2,-0.5)\botlabel{i} to[out=90,in=180] (0.3,0.23) to[out=0,in=up] (0.5,0) to[out=down,in=0] (0.3,-0.23) to[out=180,in=-90] (-0.2,0.5);
\pinpin{-0.15,-.2}{-0.15,.2}{-.9,.2}{1+x-y};
\anticlockwiseinternalbubbleR{.45,0};
\region{.8,0}{\lambda};
\end{tikzpicture}
-\sum_{i \neq j \in I}
\begin{tikzpicture}[KM,anchorbase,scale=1.1]
\draw[to-] (-0.25,-0.5)\botlabel{i} to (-0.25,0.5);
\draw[-to] (0,0) arc(-180:180:0.27);
\anticlockwiseinternalbubbleR{.5,0};
\strand{.3,-0.38}{j};
\region{.9,0}{\lambda};
\pinpin{0.02,-.1}{-.25,-.1}{-1.5,-.1}{(i-j+x-y)^{-1}};
\end{tikzpicture}\ .\label{julia}
\end{align}
We look at the two terms on the right hand side of this separately.
Expanding the definition of the internal bubble like we did in the proof of Claim 1, the first term is
$$
\sum_{\substack{J, K\text{ with }|J|<\infty\\J \sqcup K=I-\{i\}}}
\sum_{\substack{(r_k)_{k \in K}\\0 \leq r_k \leq h_k(\lambda)}}
\begin{tikzpicture}[KM,centerzero]
\draw[-to] (-3.8,1.2)[out=-90,in=120] to (-3.5,-.2) to [out=-60,in=-90] (0,-.4) to[out=90,in=-90] (0,.4) to [out=90,in=60,looseness=.9] (-3.5,.2) to [out=-120,in=90] (-3.8,-1.2)\botlabel{i};
\node at (-2.8,.25) {$\prod_{j \in J}$};
\draw[to-] (.7,.25) arc(360:0:0.2);
\strand{.5,-.08}{j};
\pinpin{0.3,.25}{0,.25}{-1.3,.25}{(i-j+x-y)^{-1}};
\pin{0,-.3}{-1.2,-.3}{(x+i)^{\sum_k r_k}};
\region{1,0}{\lambda};
\pinpin{-3.7,-.23}{-3.7,.23}{-4.6,.23}{1+x-y};
\end{tikzpicture}\ {\textstyle\prod_{k \in K}} \beta_{k,r_k}.
$$
For $j \neq i$, the expression $(i-j+x-y)^{-1}$ here is a formal power series in $x$ and $y$, which makes sense because they are locally nilpotent. To ensure that that they are actually nilpotent,
from now on, we consider 
this natural transformation evaluated on some fixed choice of
indecomposable object $V \in \ob\catR_\lambda$ (this is sufficient for the proof).
Then, 
$(i-j+x-y)^{-1}$ can be replaced by a polynomial $f_j(x,y) \in \kk[x,y]$
of some large degree in $x$ and $y$ 
(depending implicitly on the choice of $V$).
After this reduction we can
``open'' the curl using \cref{Hcurlrels,trick} as in \cite[(3.21)]{HKM}
to obtain
$$
\sum_{\substack{J, K\text{ with }|J|<\infty\\J \sqcup K=I-\{i\}}}
\sum_{\substack{(r_k)_{k \in K}\\0 \leq r_k \leq h_k(\lambda)}}
\left[
(u+i)^{\sum_k r_k}\ 
\begin{tikzpicture}[KM,baseline=3pt]
\draw[-to] (-2.5,.7)to (-2.5,-.2)\botlabel{i};
\node at (-.5,.25) {$\prod_{j \in J}$};
\draw[to-] (1.75,.25) arc(360:0:0.2);
\strand{1.55,-.08}{j};
\pin{1.35,.25}{.55,.25}{f_j(u,x)};
\draw[to-] (-2.1,.25) arc(180:-180:0.2);
\strand{-1.9,-.1}{i};
\node at (-1.4,.25) {$\color{black}(u)$};
\region{2,.25}{\lambda};
\circled{-2.5,0.25}{u};
\end{tikzpicture} 
\right]_{u:-1}{\textstyle\prod_{k \in K}} \beta_{k,r_k}.
$$
Replacing $u$ by $u-i$,
this is
$$
\sum_{\substack{J, K\text{ with }|J|<\infty\\J \sqcup K=I-\{i\}}}
\sum_{\substack{(r_k)_{k \in K}\\0 \leq r_k \leq h_k(\lambda)}}
\!\!\!\!\!\left[
u^{\sum_k r_k}\ 
\begin{tikzpicture}[KM,baseline=3pt]
\draw[-to] (-2.5,.7)to (-2.5,-.2)\botlabel{i};
\node at (-.1,.25) {$\prod_{j \in J}$};
\draw[to-] (2.55,.25) arc(360:0:0.2);
\strand{2.35,-.08}{j};
\pin{2.15,.25}{1.15,.25}{f_j(u-i,x)};
\draw[to-] (-2.2,.25) arc(180:-180:0.2);
\strand{-2,-.1}{i};
\node at (-1.23,.25) {$\color{black}(u-i)$};
\region{2.8,.25}{\lambda};
\pin{-2.5,0.25}{-3.6,.25}{(u-x-i)^{-1}};
\end{tikzpicture}
\right]_{u:-1}{\textstyle\prod_{k \in K}} \beta_{k,r_k}.
$$
Then we expand the bubble generating function using \cref{sbgf}
to split this as the sum of the following two terms, writing them using fancy notation:
\begin{align}
\sum_{\substack{J, K\text{ with }|J|<\infty\\i \in K, J \sqcup K=I}}
\sum_{\substack{(r_k)_{k \in K}\\0 \leq r_k \leq h_k(\lambda)}}
\widehat\theta_{J,i}\left(\left[
\frac{u^{\sum_k r_k}\prod_{j \in J} f_j(u-i,x_j-j)}{u-x}
\right]_{u:-1}\right)
{\textstyle\prod_{k \in K}} \beta_{k,r_k},\label{lhs1}\\
\sum_{\substack{J, K\text{ with }|J|<\infty\\i \in J, J \sqcup K=I}}
\sum_{\substack{(r_k)_{k \in K}\\0 \leq r_k \leq h_k(\lambda)}}
\widehat\theta_{J,i}\left(\left[\frac{u^{\sum_k r_k}\prod_{i \neq j \in J} f_j(u-i,x_j-j)}{(u-x)(u-x_i)}\right]_{u:-1}\right)
{\textstyle\prod_{k \in K}} \beta_{k,r_k},\label{lhs2}
\end{align}
We simplify \cref{lhs1} by first using \cref{trick} to obtain
\begin{equation*}
\sum_{\substack{J, K\text{ with }|J|<\infty\\i \in K, J \sqcup K=I}}
\sum_{\substack{(r_k)_{k \in K}\\0 \leq r_k \leq h_k(\lambda)}}
\widehat\theta_{J,i}\left(
x^{\sum_k r_k}\textstyle\prod_{j \in J} f_j(x-i,x_j-j)
\right){\textstyle\prod_{k \in K}} \beta_{k,r_k},
\end{equation*}
then noting that each $f_j(x-i,x_j-j)$
can be replaced by $(x-x_j)^{-1}$ as they have the same image 
under $\widehat\theta_{J,i}$, so that \cref{lhs1} equals
\begin{equation}\label{LHS1}
\sum_{\substack{J, K\text{ with }|J|<\infty\\i \in K, J \sqcup K=I}}
\sum_{\substack{(r_k)_{k \in K}\\0 \leq r_k \leq h_k(\lambda)}}
\widehat\theta_{J,i}\left(
\frac{x^{\sum_k r_k}}{\prod_{j \in J} (x-x_j)}
\right){\textstyle\prod_{k \in K}} \beta_{k,r_k}.
\end{equation}
For \cref{lhs2}, we use partial fractions to expand it as the sum of the following two terms:
\begin{align*}
\sum_{\substack{J, K\text{ with }|J|<\infty\\i \in J, J \sqcup K=I}}
\sum_{\substack{(r_k)_{k \in K}\\0 \leq r_k \leq h_k(\lambda)}}
\widehat\theta_{J,i}\left(\frac{x^{\sum_k r_k}\prod_{i \neq j \in J} f_j(x-i,x_j-j)}{x-x_i}\right)
{\textstyle\prod_{k \in K}} \beta_{k,r_k},\\
\sum_{\substack{J, K\text{ with }|J|<\infty\\i \in J, J \sqcup K=I}}
\sum_{\substack{(r_k)_{k \in K}\\0 \leq r_k \leq h_k(\lambda)}}
\widehat\theta_{J,i}\left(\frac{x_i^{\sum_k r_k}\prod_{i \neq j \in J} f_j(x_i-i,x_j-j)}{x_i-x}\right)
{\textstyle\prod_{k \in K}} \beta_{k,r_k}.
\end{align*}
Then we move $f_j(x-i,x_j-j)$ and
$f_j(x_i-i,x_j-j)$ in the numerator
to $(x-x_j)$ and $(x_i-x_j)$ in the denominator, respectively, to obtain finally the sum of the terms
\begin{align}\label{LHS2a}
\sum_{\substack{J, K\text{ with }|J|<\infty\\i \in J, J \sqcup K=I}}
\sum_{\substack{(r_k)_{k \in K}\\0 \leq r_k \leq h_k(\lambda)}}
\widehat\theta_{J,i}\left(\frac{x^{\sum_k r_k}}{\prod_{j \in J} (x-x_j)}\right)
{\textstyle\prod_{k \in K}} \beta_{k,r_k},\\\label{LHS2b}
\sum_{\substack{J, K\text{ with }|J|<\infty\\i \in J, J \sqcup K=I}}
\sum_{\substack{(r_k)_{k \in K}\\0 \leq r_k \leq h_k(\lambda)}}
\widehat\theta_{J,i}\left(\frac{x_i^{\sum_k r_k}}{(x_i-x)\prod_{i \neq j \in J} (x_i-x_j)}\right)
{\textstyle\prod_{k \in K}} \beta_{k,r_k}.
\end{align}
Finally we look at the second term from 
\cref{julia}, which expands as
$$
-\sum_{i \neq j\in I}
\sum_{\substack{J, K\text{ with }|J|<\infty\\J \sqcup K=I-\{j\}}}
\sum_{\substack{(r_k)_{k \in K}\\0 \leq r_k \leq h_k(\lambda)}}
\begin{tikzpicture}[KM,centerzero]
\draw[-to] (-3.8,1.4) to (-3.8,-1.4) \botlabel{i};
\draw[-to] (-3.5,0) to (-3.5,-.6) to [out=-90,in=-90,looseness=.9] (0,-.6)
to (0,.4) to [out=90,in=90,looseness=.9] (-3.5,.4)
to (-3.5,0);
\pinpin{-3.5,-.2}{-3.8,-.2}{-5.1,-.2}{(i-j+x-y)^{-1}};
\node at (-2.84,.25) {$\prod_{h \in J}$};
\draw[to-] (.7,.25) arc(360:0:0.2);
\strand{.5,-.08}{h};
\pinpin{0.3,.25}{0,.25}{-1.3,.25}{(j-h+x-y)^{-1}};
\pin{0,-.3}{-1.1,-.3}{(x+j)^{\sum_k r_k}};
\region{1,0}{\lambda};
\strand{-.5,-1.4}{j};
\end{tikzpicture} \ 
{\textstyle\prod_{k \in K}} \beta_{k,r_k}.
$$
Reindexing the summation and using fancy notation, this is
\begin{equation}\label{lhs3}
\sum_{\substack{J, K\text{ with }|J|<\infty\\J \sqcup K=I}}
\sum_{\substack{(r_k)_{k \in K}\\0 \leq r_k \leq h_k(\lambda)}}
\sum_{i \neq j \in J}
\widehat\theta_{J,i}\left(
\frac{x_j^{\sum_k r_k}}{(x_j-x)\prod_{j \neq h \in J}(x_j-x_h)}\right)
{\textstyle\prod_{k \in K}} \beta_{k,r_k}.
\end{equation}
Now we have expanded the right hand side of \cref{beer}
as the sum of \cref{rhs1,rhs2}, and the left hand side as the sum of
\cref{LHS1,LHS2a,LHS2b,lhs3}.
The two sides are equal: the sum of the terms from \cref{rhs1} with $i \in J$ 
is \cref{LHS2a}, the sum of the terms from \cref{rhs1} with $i \in K$
is \cref{LHS1}, the term from \cref{rhs2} with $j=i$ is
\cref{LHS2b}, and the sum of the 
remaining terms from \cref{rhs2} is \cref{lhs3}.
Thus, \cref{beer} is proved.

\vspace{2mm}
\noindent
\underline{Claim 3}.
{\em We have that
$\ 
\begin{tikzpicture}[H,centerzero,scale=1.6]
\draw[-,shadow] (0.2,-0.1) to [looseness=2.2,out=90,in=90] (-0.2,-0.1);
\draw[-to] (-.2,-.1) to (-.2,-.12);
\region{.35,0}{\lambda};
\end{tikzpicture}
=\ 
\begin{tikzpicture}[H,centerzero,scale=1.6]
\draw[-to] (0.2,-0.1) to [looseness=2.2,out=90,in=90] (-0.2,-0.1);
\region{.35,0}{\lambda};
\end{tikzpicture},$ proving the first equality in \cref{leftwardbuggers}.}

\noindent
\underline{Proof}.
We first prove this assuming that $\kappa < 0$.
By the definition of $\ \begin{tikzpicture}[H,centerzero,scale=1.2]
\draw[-to] (0.2,-0.1) to [looseness=2.2,out=90,in=90] (-0.2,-0.1);
\region{.35,0}{\lambda};
\end{tikzpicture}$
as the last entry of the inverse of the matrix $M_\kappa$ from \cref{rybnikov}, this reduces to checking the following identities:
\begin{align}\label{thai}
\begin{tikzpicture}[H,centerzero,scale=1.2]
\draw[-] (-.2,-.3) to [out=90,in=-90,looseness=1] (.18,.15);
\draw[shadow] (.18,.15) to [out=90,in=90,looseness=1.7] (-.18,.15);
%\draw[-to] (-.18,.15) to (-.18,.13);
\draw[-to] (-.18,.15) to [out=-90,in=90,looseness=1] (.2,-.3);
\region{.35,0}{\lambda};
\end{tikzpicture}
&=0,
&
\begin{tikzpicture}[H,baseline=-1mm]
\draw[-,shadow] (-0.25,0) arc(180:0:0.25);
\draw[-to] (-.25,0) to (-.25,-.02);
\draw[-] (0.25,0) arc(0:-180:0.25);
\multopendot{0.25,0}{west}{n};
\region{.8,0}{\lambda};
\end{tikzpicture}
&=\delta_{n,-\kappa-1}\text{ for $0 \leq n \leq -\kappa-1$.}
\end{align}
To check these, Claim 1 plus 
\cref{bubblegeneratingfunction1,dunking}
implies that
$\ \begin{tikzpicture}[H,anchorbase,scale=1]
\draw[-,shadow] (-0.2,0) arc(180:0:0.2);
\draw[-to] (-.2,0) to (-.2,-.02);
\draw[-] (0.2,0) arc(0:-180:0.2);
\node at (0.5,0) {$\color{black}(u)$};
\region{.95,0}{\lambda};
\end{tikzpicture}\in u^\kappa + u^{\kappa-1} 
Z(\catR_\lambda)\llbracket u^{-1}\rrbracket$.
The second of the equalities \cref{thai} follows immediately from this.
For the first one, we attach 
 $\ \begin{tikzpicture}[H,centerzero,scale=1.2]
\draw[to-] (0.2,-0.1) to [looseness=2.2,out=90,in=90] (-0.2,-0.1);
\region{.35,0}{\lambda};
\end{tikzpicture}$ to the top of the identity from Claim 2, noting that the right hand side is zero when $\kappa < 0$. 

It remains to treat the case that $\kappa = 0$.
We have that
$\ \begin{tikzpicture}[H,centerzero,scale=.8]
\draw[to-] (-0.2,-0.4) to[out=45,in=down] (0.15,0) to[out=up,in=-45] (-0.2,0.4);
\draw[-to] (0.2,-0.4) to[out=135,in=down] (-0.15,0) to[out=up,in=225] (0.2,0.4);
\end{tikzpicture}
=\begin{tikzpicture}[H,centerzero,scale=1]
\draw[to-] (-0.14,-0.3) -- (-0.14,0.3);
\draw[-to] (0.14,-0.3) -- (0.14,0.3);
\end{tikzpicture}\ 
$
by a special case of \cref{Haltquadratic}.
If we attach 
$\ \begin{tikzpicture}[H,centerzero,scale=1.2]
\draw[-,shadow] (0.2,-0.1) to [looseness=2.2,out=90,in=90] (-0.2,-0.1);
\draw[-to] (-.2,-.1) to (-.2,-.12);
\region{.35,0}{\lambda};
\end{tikzpicture}$
to the top of this identity, keeping the definition of 
$\ \begin{tikzpicture}[H,centerzero,scale=1.2]
\draw[-to] (0.2,-0.1) to [looseness=2.2,out=90,in=90] (-0.2,-0.1);
\region{.35,0}{\lambda};
\end{tikzpicture}$ in mind, the problem in hand reduces
to showing that
$$
\begin{tikzpicture}[H,centerzero,scale=1.2]
\draw[-] (-.2,-.3) to [out=90,in=-90,looseness=1] (.18,.15);
\draw[shadow] (.18,.15) to [out=90,in=90,looseness=1.7] (-.18,.15);
%\draw[->] (-.18,.15) to (-.18,.13);
\draw[-to] (-.18,.15) to [out=-90,in=90,looseness=1] (.2,-.3);
\region{.35,0}{\lambda};
\end{tikzpicture} = \begin{tikzpicture}[H,baseline=1mm]
\draw[-to] (-0.15,-0.05) to [out=90,in=90,looseness=3](0.35,-0.05);
\region{.6,0.2}{\lambda};
\end{tikzpicture}\ .
$$
This is easily checked in a similar way to the proof of the first equality in \cref{thai}.

\vspace{2mm}
\noindent
\underline{Claim 4}.
{\em 
Both of the equalities in \cref{newbub} hold.}

\noindent
\underline{Proof}.
Claim 3 implies that $\ \begin{tikzpicture}[H,baseline=-1mm,scale=.8]
\draw[-,shadow] (-0.25,0) arc(180:0:0.25);
\draw[-to] (-.25,0) to (-.25,-.02);
\draw[-] (0.25,0) arc(0:-180:0.25);
\node at (.58,0) {\color{black}$(u)$};
\region{1.07,0}{\lambda};
\end{tikzpicture}
=\ \begin{tikzpicture}[H,baseline=-1mm,scale=.8]
\draw[to-] (-0.25,0) arc(180:-180:0.25);
\node at (.58,0) {\color{black}$(u)$};
\region{1.07,0}{\lambda};
\end{tikzpicture}$. Consequently, the first of the identities in \cref{newbub} follows from Claim 1.
The second of these identities then follows by taking inverses on both sides, using \cref{infgrass,Hmiles}.

\noindent
\underline{Claim 5}.
{\em 
$\begin{tikzpicture}[H,anchorbase,scale=1]
\draw[to-] (0,0.5) to[out=-90,in=120] (0.08,-0.05);
\draw[shadow] (0.08,-.05) to [out=-60,in=down,looseness=2] (0.6,0);
%\draw[-to] (.08,-.03) to [out=120,in=-60] (.07,-.028);
\draw (.6,0) to[looseness=2,out=up,in=60] (.08,.05) to[out=-120,in=up] (0,-0.5);
\end{tikzpicture}
=-
\left[
\begin{tikzpicture}[H,anchorbase,scale=1.1]
\draw[-to] (-0.8,-0.5) -- (-0.8,0.5);
\circled{-0.8,0}{u};
\draw[-to] (-.4,0) arc(180:-180:0.2);
\node at (0.26,0) {$\color{black}(u)$};
\end{tikzpicture}
\!\right]_{u:-1}\!\!\!$.}

\noindent
\underline{Proof}.
The follows by a similar argument to the proof of Claim 2, using 
the second of the identities \cref{newbub} 
established in
Claim 4 in place of the identity from Claim 1. 

\vspace{2mm}
\noindent
\underline{Claim 6}.
{\em $\ \begin{tikzpicture}[H,centerzero,scale=1.6]
\draw[-,shadow] (0.2,0.1) to [looseness=2.2,out=-90,in=-90] (-0.2,0.1);
\draw[-to] (-.2,.1) to (-.2,.12);
\region{.35,0}{\lambda};
\end{tikzpicture}
=
\begin{tikzpicture}[H,centerzero,scale=1.6]
\draw[-to] (0.2,0.1) to [looseness=2.2,out=-90,in=-90] (-0.2,0.1);
\region{.35,0}{\lambda};
\end{tikzpicture}$, proving the second equality in \cref{leftwardbuggers}.}

\noindent
\underline{Proof}.
Using Claim 5 for the unlabelled equality, we have that
\begin{align*}
\begin{tikzpicture}[H,centerzero,scale=1.6]
\draw[-,shadow] (0.2,0.1) to [looseness=2.2,out=-90,in=-90] (-0.2,0.1);
\draw[-to] (-.2,.1) to (-.2,.12);
\region{.35,0}{\lambda};
\end{tikzpicture}
&\stackrel{\cref{Haltquadratic}}{=}\ 
\begin{tikzpicture}[H,centerzero,scale=1.2]
\draw[-to] (-0.15,-0.4) to[out=90,in=down] (0.15,0) to[out=up,in=-45] (-0.15,0.4);
\draw [shadow](-.15,-.4) to [out=-90,in=-90,looseness=1.5] (0.15,-0.4);
%\draw[-to] (-.15,-.4) to (-.15,-.38);
\draw (.15,-.4) to[out=90,in=down] (-0.15,0) to[out=up,in=225] (0.2,0.4);
\region{.35,0}{\lambda};
\end{tikzpicture}-
\left[\ \begin{tikzpicture}[H,centerzero,scale=1.6]
\draw[-to] (0.2,0.3) to [looseness=2.2,out=-90,in=-90] (-0.2,0.3);
\circled{-.15,.12}{u};
\draw[to-] (0.27,0) arc(180:-180:0.142);
\draw[-] (0.083,-.2) arc(0:180:0.142);
\draw[shadow] (-.2,-.2) arc(180:360:0.142);
\draw[-to] (-.2,-.2) to (-.2,-.18);
\circled{-.12,-.1}{u};
\node at (0.73,0) {$\color{black}(u)$};
\region{1,0}{\lambda};
\end{tikzpicture}\ \right]_{u:-1}
\stackrel{\cref{Hmiles}}{=}\ 
\begin{tikzpicture}[H,centerzero,scale=1.2]
\draw[-to] (-0.15,-0.4) to[out=90,in=down] (0.15,0) to[out=up,in=-45] (-0.15,0.4);
\draw [shadow](-.15,-.4) to [out=-90,in=-90,looseness=1.5] (0.15,-0.4);
%\draw[-to] (-.15,-.4) to (-.15,-.38);
\draw (.15,-.4) to[out=90,in=down] (-0.15,0) to[out=up,in=225] (0.2,0.4);
\region{.35,0}{\lambda};
\end{tikzpicture}\\
&
=
-\left[\ 
\begin{tikzpicture}[H,centerzero,scale=1.6]
\draw[to-] (-.2,.3) to [out=-90,in=90,looseness=1] (.15,-.15) to [out=-90,in=-90,looseness=1.7] (-.15,-.15) to [out=90,in=-90,looseness=1] (.2,.3);
\circled{0.14,-0.15}{u};
\draw[to-] (.5,0) arc(-180:180:0.142);
\node at (.95,0) {$\color{black}(u)$};
\region{1.2,0}{\lambda};
\end{tikzpicture}
\right]_{u:-1}
\stackrel{\cref{Hcurlrels}}{=}
-\left[\ 
\begin{tikzpicture}[H,centerzero,scale=1.6]
\draw[-to] (0.35,0.15) to [looseness=2.2,out=-90,in=-90] (-.1,0.15);
\circled{-.05,-.03}{u};
\draw[to-] (1.2,0) arc(-180:180:0.142);
\node at (1.65,0) {$\color{black}(u)$};
\draw[to-] (.5,0) arc(180:-180:0.142);
\node at (.95,0) {$\color{black}(u)$};
\region{1.9,0}{\lambda};
\end{tikzpicture}\right]_{u:-1}
\stackrel{\cref{Hmiles}}{=}
\begin{tikzpicture}[H,centerzero,scale=1.6]
\draw[-to] (0.2,0.1) to [looseness=2.2,out=-90,in=-90] (-0.2,0.1);
\region{.35,0}{\lambda};
\end{tikzpicture}\ .
\end{align*}
\end{proof}

\begin{cor}\label{averagecor}
The leftward Heisenberg crossing is related to the leftward Kac-Moody crossings by
\begin{align}
\begin{tikzpicture}[centerzero,KM,scale=1.4]
\draw[-to] (0.28,-.28) \botlabel{i} to (-0.28,.28)\toplabel{i};
\draw[to-] (-0.28,-.28) \botlabel{j} to (0.28,.28)\toplabel{j};
\draw[H] (0.1,-.1) to (-0.1,.1);
\draw[H] (-0.1,-.1) to (0.1,.1);
\notch[45]{.1,-.1};
\notch[45]{-.1,.1};
\notch[-45]{.1,.1};
\notch[-45]{-.1,-.1};
\region{.4,0}{\lambda};
\end{tikzpicture}
=
\begin{tikzpicture}[centerzero,KM,scale=1.4]
\draw[-to] (0.28,-.28) \botlabel{i} to (-0.28,.28)\toplabel{i};
\draw[to-] (-0.28,-.28) \botlabel{j} to (0.28,.28)\toplabel{j};
\projcr{0,0};
\region{.4,0}{\lambda};
\end{tikzpicture}
&=
\begin{dcases}
\begin{tikzpicture}[centerzero,KM,scale=1.4]
\draw[-to] (0.28,-.28) \botlabel{i} to (-0.28,.28)\toplabel{i};
\draw[to-] (-0.28,-.28) \botlabel{i} to (0.28,.28)\toplabel{i};
\anticlockwiseinternalbubbleNW{-.09,-.06};
\clockwiseinternalbubbleSE{.07,.03};
\pinpin{-.21,.21}{.21,.21}{.8,.21}{1+x-y};
\region{.7,-.1}{\lambda};
\end{tikzpicture}
&\text{if $i=j$}\\
\begin{tikzpicture}[centerzero,KM,scale=1.4]
\draw[-to] (0.28,-.28) \botlabel{i} to (-0.28,.28)\toplabel{i};
\draw[to-] (-0.28,-.28) \botlabel{i-1} to (0.28,.28)\toplabel{i-1};
\anticlockwiseinternalbubbleNW{-.09,-.06};
\clockwiseinternalbubbleSE{.07,.03};
\pinpin{-.21,.21}{.21,.21}{1.05,.21}{(1+x-y)^{-1}};
\region{.9,-.15}{\lambda};
\end{tikzpicture}
&\text{if $i=j+1$}\\
(-1)^{\delta_{i,j-1}}
\begin{tikzpicture}[centerzero,KM,scale=1.4]
\draw[-to] (0.28,-.28) \botlabel{i} to (-0.28,.28)\toplabel{i};
\draw[to-] (-0.28,-.28) \botlabel{j} to (0.28,.28)\toplabel{j};
\anticlockwiseinternalbubbleNW{-.09,-.06};
\clockwiseinternalbubbleSE{.07,.03};
\pinpin{-.21,.21}{.21,.21}{1.9,.21}{(i-j+x-y)^{-1}(i-j-1+x-y)};
\region{1.2,-.15}{\lambda};
\end{tikzpicture}
&\text{if $i \neq j, j+1$}
\end{dcases}
\label{lastgasp1}\\\intertext{for any $i,j \in I$.  For $i \neq j$ in $I$, we also have that}
\begin{tikzpicture}[centerzero,KM,scale=1.4]
\draw[-to] (0.28,-.28) \botlabel{i} to (-0.28,.28)\toplabel{j};
\draw[to-] (-0.28,-.28) \botlabel{i} to (0.28,.28)\toplabel{j};
\draw[H] (0.1,-.1) to (-0.1,.1);
\draw[H] (-0.1,-.1) to (0.1,.1);
\notch[45]{.1,-.1};
\notch[45]{-.1,.1};
\notch[-45]{.1,.1};
\notch[-45]{-.1,-.1};
\region{.4,0}{\lambda};
\end{tikzpicture}
=
\begin{tikzpicture}[centerzero,KM,scale=1.4]
\draw[-to] (0.28,-.28) \botlabel{i} to (-0.28,.28)\toplabel{j};
\draw[to-] (-0.28,-.28) \botlabel{i} to (0.28,.28)\toplabel{j};
\projcr{0,0};
\region{.4,0}{\lambda};
\end{tikzpicture}
&=
\begin{tikzpicture}[KM,anchorbase,scale=1.2]
\draw (0.4,-0.55) to (0.4,-.35);
\draw[-to] (-0.1,-.35) to (-0.1,-0.55)\botlabel{i};
\draw (.4,-.35) to[out=90,in=90,looseness=2] (-.1,-.35);
\draw[-to,internalbubble] (.5,-.32)  arc(0:360:.13);
\region{1.3,-.3}{\lambda};
\draw[-] (0.4,0.55) to (0.4,.35);
\draw[-to] (-0.1,.35) to (-0.1,0.55)\toplabel{j};
\draw (.4,.35) to[out=-90,in=-90,looseness=2] (-.1,.35);
\draw[-to,internalbubble] (-.2,.32)  arc(180:-180:.13);
\pinpin{.3,-.11}{.3,.11}{1.4,.11}{(i-j+x-y)^{-1}};
\end{tikzpicture}\ .
\label{lastgasp2}
\end{align}
We know already that all other projected leftward Heisenberg crossings are zero.
\end{cor}

\begin{proof}
The identity \cref{lastgasp1} follows by attaching 
$\ \begin{tikzpicture}[KM,anchorbase,scale=1]
\draw (0.3,-0.4) to (0.3,-.2);
\draw[-to] (-0.1,-.2) to (-0.1,-0.4)\botlabel{j};
\draw (.3,-.2) to[out=90,in=90,looseness=2] (-.1,-.2);
\draw[-to,internalbubble] (.38,-.15) arc(0:360:.12);
\end{tikzpicture}\ $
to the top left and 
$\ \begin{tikzpicture}[KM,anchorbase,scale=1]
\draw[-] (0.3,0.4) to (0.3,.2);
\draw[-to] (-0.1,.2) to (-0.1,0.4)\toplabel{j};
\draw (.3,.2) to[out=-90,in=-90,looseness=2] (-.1,.2);
\draw[-to,internalbubble] (-.18,.15) arc(180:-180:.12);
\end{tikzpicture}\ $
to the bottom right of
the identity \cref{upwardmagicalt}, then simplifying using other
relations including \cref{leftwardbuggers,wax}.
The proof of 
\cref{lastgasp2} is similar, starting from \cref{eggs} instead of \cref{upwardmagicalt}.
\end{proof}

\begin{cor}\label{other}
$\begin{tikzpicture}[KM,centerzero,scale=.8]
\draw[-to] (0,-.5)\botlabel{i} to (0,.5);
\anticlockwiseinternalbubbleR{0,0};
\region{0.5,0}{\lambda};
\end{tikzpicture}
=
\left(\ \, \begin{tikzpicture}[KM,centerzero,scale=.8]
\draw[-to] (0,-.5)\botlabel{i} to (0,.5);
\clockwiseinternalbubbleL{0,0};
\region{0.5,0}{\lambda};
\end{tikzpicture}\right)^{-1}$.
\end{cor}

\begin{proof}
The following checks that the counterclockwise and clockwise internal bubbles are inverses on one side:
\begin{align*}
\begin{tikzpicture}[KM,centerzero,scale=.8]
\draw[-to] (0,-.8)\botlabel{i} to (0,.8);
\clockwiseinternalbubbleL{0,.35};
\anticlockwiseinternalbubbleR{0,-.35};
\region{0.55,0}{\lambda};
\end{tikzpicture}
&=
\begin{tikzpicture}[KM,centerzero,scale=.8]
\draw[-to] (1,-.8)\botlabel{i} to (1,0.17) to [out=90,in=90,looseness=1.5] (0,.17) to (0,-.17) to[out=-90,in=-90,looseness=1.5] (-1,-.17) to (-1,.8);
\clockwiseinternalbubbleL{-1,.2};
\anticlockwiseinternalbubbleR{1,-.2};
\region{1.5,0}{\lambda};
\end{tikzpicture}=
\sum_{j \in I}\ 
\begin{tikzpicture}[KM,centerzero,scale=.8]
\draw[-] (1,-.8)\botlabel{i} to (1,0.17);
\draw[H] (1,.17) to [out=90,in=90,looseness=1.5] (0,.17);
\draw (0,.17) to (0,-.17);
\draw[H] (0,-.17) to[out=-90,in=-90,looseness=1.5] (-1,-.17);
\draw[-to] (-1,-.17) to (-1,.8)\toplabel{i};
\strand{0.2,0}{j};
\notch{0,.17};\notch{0,-.17};
\notch{1,.17};\notch{-1,-.17};
\region{1.3,0}{\lambda};
\end{tikzpicture}
=
\begin{tikzpicture}[KM,centerzero,scale=.8]
\draw[-] (1,-.8)\botlabel{i} to (1,0.17);
\draw[H] (1,.17) to [out=90,in=90,looseness=1.5] (0,.17) to
(0,-.17) to[out=-90,in=-90,looseness=1.5] (-1,-.17);
\draw[-to] (-1,-.17) to (-1,.8)\toplabel{i};
\notch{1,.17};\notch{-1,-.17};
\region{1.3,0}{\lambda};
\end{tikzpicture}
=
\begin{tikzpicture}[KM,centerzero,scale=.8]
\draw[-to] (0,-.8)\botlabel{i} to (0,.8)\toplabel{i};
\draw[H] (0,-.4) to (0,.4);
\notch{0,-.4};\notch{0,.4};
\region{0.3,0}{\lambda};
\end{tikzpicture}
=
\begin{tikzpicture}[KM,centerzero,scale=.8]
\draw[-to] (0,-.8)\botlabel{i} to (0,.8)\toplabel{i};
\region{0.3,0}{\lambda};
\end{tikzpicture}\ .
\end{align*}
A similar calculation with the other leftward zig-zag identity 
shows that they are inverses on the other side too.
\end{proof}

\begin{rem}
The formulae \cref{leftwardbuggers,lastgasp1,lastgasp2} are new, but
\cref{newbub} was proved already 
in \cite[(5.37)]{HKM}
by a more indirect method. The proof just explained depends on \cref{maintheorem}, hence, 
on \cite[Lem.~4.9, Lem.~4.10]{HKM}, but is independent of all of the results of \cite[Sec.~5]{HKM}.
\end{rem}

\cref{average,averagecor} explain how to write
the leftward Heisenberg cap, cup and crossing in terms of the Kac-Moody ones and internal bubbles. 
When applying \cref{maintheorem}, one really wants to be able to go in the other direction. Since the internal bubbles are invertible
by \cref{other}, it is 
quite trivial to rearrange the formulae obtained so far 
to write the leftward Kac-Moody caps, cups and crossings in terms of the Heisenberg ones and internal bubbles.
However, one still needs to be able to determine the internal bubbles, which were
defined explicitly in terms of the Kac-Mooody dotted bubbles.
The final theorem in this section completes the picture by explaining
how the Kac-Moody bubble generating functions can be obtained 
from the Heisenberg bubble generating functions. It is not very explicit, but this is to be expected in this place.

Given an indecomposable $V \in \ob\catR$,
we denote its endomorphism algebra by $Z_V$. This is a finite-dimensional commutative local algebra. We denote its unique maximal ideal by $J_V$.

\begin{lem}\label{badder}
Let $Z$ be any finite-dimensional commutative local algebra
and $J$ be its unique maximal ideal.
Let $\bar Z := Z / J \cong \kk$, denoting the canonical homomorphism
$Z[x] \rightarrow \bar Z[x]$ by $f(x) \mapsto \bar f(x)$.
Suppose that we are given $k \in \Z$ and a rational function
$\frac{f(u)}{g(u)}$ in $Z(u)
\cap \left(u^k + u^{k-1} Z\llbracket u^{-1}\rrbracket\right)$.
Let $h_i \in \Z$ be the multiplicity of $i \in \kk$ as a zero or pole
of $\frac{\bar f(u)}{\bar g(u)} \in \bar Z(u)$.
There are unique rational functions
$\frac{f_i(u)}{g_i(u)}$ in $Z(u) \cap \left(u^{h_i}+u^{h_i-1} J\llbracket u^{-1}\rrbracket\right)$ such that 
$\frac{f_i(u)}{g_i(u)} = 1$ for all but finitely many $i$ and
$$
\frac{f(u)}{g(u)} = \prod_{i \in \kk} \frac{f_i(u-i)}{g_i(u-i)}.
$$
\end{lem}

\begin{proof}
This follows from \cite[Cor.~2.4]{HKM}, which is the analogous statement for polynomials rather than rational functions.
\end{proof}

\begin{lem}\label{bad}
Let $V \in \ob\catR_\lambda$ be an indecomposable object.
The Heisenberg bubble generating functions
$\ \begin{tikzpicture}[H,baseline=-1mm,scale=.9]
\draw[to-] (-0.25,0) arc(180:-180:0.25);
\node at (.58,0) {$\color{black}(u)$};
\end{tikzpicture}$ and
$\ \begin{tikzpicture}[H,baseline=-1mm,scale=.9]
\draw[to-] (-0.25,0) arc(-180:180:0.25);
\node at (.58,0) {$\color{black}(u)$};
\end{tikzpicture}$
act on $V$ as rational functions belonging to
$Z_V(u) \cap \left(u^{\kappa} + u^{\kappa-1} Z_V\llbracket u^{-1}\rrbracket\right)$
and
$Z_V(u) \cap \left(-u^{-\kappa} + u^{-\kappa-1} Z_V\llbracket u^{-1}\rrbracket\right)$, respectively.
\end{lem}

\begin{proof}
It suffices to prove this for
the counterclockwise bubble generating function, since the other one is its inverse (up to a sign) by \cref{Hmiles}.
Let $f(x) \in \kk[x]$ be the (monic) minimal 
polynomial of the endomorphism of $E V$ defined by $\begin{tikzpicture}[H,centerzero]
\draw[-to] (0,-0.2) -- (0,0.2);
\opendot{0,-0.02};
\end{tikzpicture}$. By \cite[Lem.~4.3(1)]{HKM},
there is a monic polynomial $g(u) \in Z_V[u]$
such that $\ \begin{tikzpicture}[H,baseline=-1mm,scale=.9]
\draw[to-] (-0.25,0) arc(180:-180:0.25);
\node at (.58,0) {$\color{black}(u)$};
\end{tikzpicture}$ acts on $V$ in the same way as the rational function
$\frac{g(u)}{f(u)}$.
%and
%$\begin{tikzpicture}[H,centerzero]
%\draw[to-] (0,-0.2) -- (0,0.2);
%\pin{0,0.02}{.6,.02}{f(x)};
%\end{tikzpicture}$ acts on $V$ as $0$.
The lemma follows using also that we know the leading term from 
\cref{Hbubblegeneratingfunction1}.
\end{proof}

\begin{theo}\label{worse}
The following hold for all $\lambda \in X$:
\begin{enumerate}
\item
The bubble generating functions
$\ \begin{tikzpicture}[KM,baseline=-1mm,scale=.9]
\draw[to-] (-0.25,0) arc(180:-180:0.25);
\node at (.59,0) {\color{black}$(u)$};
\region{1.1,0}{\lambda};
\strand{0,-.4}{i};
\end{tikzpicture}$
are the unique elements of $Z(\catR_\lambda)\lround u^{-1}\rround$
such that
$\ \begin{tikzpicture}[H,baseline=-1mm,scale=.9]
\draw[to-] (-0.25,0) arc(180:-180:0.25);
\node at (.59,0) {\color{black}$(u)$};
\region{1.1,0}{\lambda};
\end{tikzpicture}
=
{\textstyle\prod_{i \in I}}\ 
\begin{tikzpicture}[KM,baseline=-1mm,scale=.9]
\draw[to-] (-0.25,0) arc(180:-180:0.25);
\node at (1,0) {\color{black}$(u-i)$};
\region{1.9,0}{\lambda};
\strand{0,-.4}{i};
\end{tikzpicture}
$
and
$\ \begin{tikzpicture}[KM,baseline=-1mm,scale=.9]
\draw[to-] (-0.25,0) arc(180:-180:0.25);
\node at (.59,0) {\color{black}$(u)$};
\region{1.1,0}{\lambda};
\strand{0,-.4}{i};
\end{tikzpicture}$
acts on $V$ as
a rational function in $Z_V(u)\cap \left(u^{h_i(\lambda)} 
+ u^{h_i(\lambda)-1} J_V\llbracket u^{-1}\rrbracket\right)$
for each indecomposable $V \in \ob\catR_\lambda$.
\item
The bubble generating functions
$\ \begin{tikzpicture}[KM,baseline=-1mm,scale=.9]
\draw[to-] (-0.25,0) arc(-180:180:0.25);
\node at (.59,0) {\color{black}$(u)$};
\region{1.1,0}{\lambda};
\strand{0,-.4}{i};
\end{tikzpicture}$
are the unique elements of $Z(\catR_\lambda)\lround u^{-1}\rround$
such that $\ \begin{tikzpicture}[H,baseline=-1mm,scale=.9]
\draw[-to] (-0.25,0) arc(180:-180:0.25);
\node at (.59,0) {$\color{black}(u)$};
\region{1.1,0}{\lambda};
\end{tikzpicture}
=-
{\textstyle\prod_{i \in I}}\ 
\begin{tikzpicture}[KM,baseline=-1mm,scale=.9]
\draw[-to] (-0.25,0) arc(180:-180:0.25);
\node at (1,0) {$\color{black}(u-i)$};
\strand{0,-.4}{i};
\region{1.9,0}{\lambda};
\end{tikzpicture}
$
and $\ \begin{tikzpicture}[KM,baseline=-1mm,scale=.9]
\draw[to-] (-0.25,0) arc(180:-180:0.25);
\node at (.59,0) {\color{black}$(u)$};
\region{1.1,0}{\lambda};
\strand{0,-.4}{i};
\end{tikzpicture}$
acts on $V$ as a rational function in
$Z_V(u) \cap \left(u^{-h_i(\lambda)} + u^{-h_i(\lambda)-1} J_V\llbracket u^{-1}\rrbracket\right)$
for each indecomposable $V \in \ob\catR_\lambda$.
\end{enumerate}
\end{theo}

\begin{proof}
This follows from \cref{badder,bad} using the known formulae 
\cref{newbub}.
\end{proof}

%% file: s8-quantum.tex
\setcounter{section}{7}

%=====================================
\section{The quantum case}\label{s8-quantum}
%=====================================

There is also a quantum version of the Heisenberg category and of \cref{maintheorem}, which allows further examples
of Kac-Moody categorifications to be obtained from various categories of representations of
quantum linear groups and associated Hecke algebras (perhaps at roots
of unity). We briefly recall the general setup in this section, then write down
the analogues of the results of \cref{s7-bubbles}, omitting the detailed
proofs since they are similar in spirit to the degenerate case.

Fix once again an algebraically closed field $\kk$.
As well as the central charge $\kappa \in \Z$, 
the definition of the quantum Heisenberg category
depends on an additional parameter $q \in \kk-\{0,\pm 1\}$.
It occurs often so we use the shorthand
\begin{equation}
z := q-q^{-1}.
\end{equation}
Then we define $\qHeis_\kappa$ to be the strict $\kk[t,t^{-1}]$-linear monoidal
category with
generating objects $E$ and $F$, whose identity endomorphisms are represented by the oriented strings $\ \begin{tikzpicture}[Q,anchorbase]
\draw[-to] (0,-0.2) -- (0,0.2);
\end{tikzpicture}\ $ and
$\ \begin{tikzpicture}[Q,anchorbase]
\draw[to-] (0,-0.2) -- (0,0.2);
\end{tikzpicture}\ $,
and the following four generating morphisms:
\begin{align*}
\begin{tikzpicture}[Q,centerzero]
\draw[-to] (0,-0.3) -- (0,0.3);
\opendot{0,0};
\end{tikzpicture}&
\colon E \rightarrow E,&
\begin{tikzpicture}[Q,centerzero,scale=.9]
\draw[-to] (0.3,-0.3) -- (-0.3,0.3);
\draw[wipe] (-0.3,-0.3) -- (0.3,0.3);
\draw[-to] (-0.3,-0.3) -- (0.3,0.3);
\end{tikzpicture}
&\colon E\otimes E  \rightarrow E \otimes E,&
\begin{tikzpicture}[Q,baseline=0mm]
\draw[-to] (-0.25,-0.15) to [out=90,in=90,looseness=3](0.25,-0.15);
\end{tikzpicture}\ \ 
&\colon E\otimes F \rightarrow \one,&
\begin{tikzpicture}[Q,baseline=-1.3mm]
\draw[-to] (-0.25,0.15) to[out=-90,in=-90,looseness=3] (0.25,0.15);
\end{tikzpicture}\ \ 
&\colon \one \rightarrow F\otimes E.
\end{align*}
The dot and the crossing, which from now on we will call the {\em
  positive} crossing, are required to be invertible.
The invertibility of the dot
means that now it makes sense to label dots by an arbitrary integer
rather than just by $n \in \N$, and it makes sense to pin Laurent polynomials
$f(x) \in \kk[x,x^{-1}]$ rather than mere polynomials to them. 
We denote the inverse of the positive crossing by
\begin{align*}
\begin{tikzpicture}[Q,centerzero,scale=.9]
\draw[-to] (-0.3,-.3) to (0.3,.3);
\draw[wipe] (0.3,-.3) to (-0.3,.3);
\draw[-to] (0.3,-.3) to (-0.3,.3);
\end{tikzpicture}
&:E\otimes E\rightarrow E \otimes E,
\end{align*}
and call this the {\em negative} crossing.
\iffalse
Thus, we have that
\begin{align}
\begin{tikzpicture}[Q,centerzero,scale=.8]
\draw[-] (0.28,-.6) to[out=90,in=-90] (-0.28,0);
\draw[-to] (-0.28,0) to[out=90,in=-90] (0.28,.6);
\draw[wipe] (-0.28,-.6) to[out=90,in=-90] (0.28,0);
\draw[-] (-0.28,-.6) to[out=90,in=-90] (0.28,0);
\draw[wipe] (0.28,0) to[out=90,in=-90] (-0.28,.6);
\draw[-to] (0.28,0) to[out=90,in=-90] (-0.28,.6);
\end{tikzpicture}
&=
\begin{tikzpicture}[Q,centerzero,scale=.8]
\draw[-to] (0.18,-.6) to (0.18,.6);
\draw[-to] (-0.18,-.6) to (-0.18,.6);
\end{tikzpicture}
=
\begin{tikzpicture}[Q,centerzero,scale=.8]
\draw[-to] (0.28,0) to[out=90,in=-90] (-0.28,.6);
\draw[wipe] (-0.28,0) to[out=90,in=-90] (0.28,.6);
\draw[-to] (-0.28,0) to[out=90,in=-90] (0.28,.6);
\draw[-] (-0.28,-.6) to[out=90,in=-90] (0.28,0);
\draw[wipe] (0.28,-.6) to[out=90,in=-90] (-0.28,0);
\draw[-] (0.28,-.6) to[out=90,in=-90] (-0.28,0);
\end{tikzpicture}
\ .\label{AHA1}
\end{align}
\fi
We also introduce the negative rightward crossing:
\begin{align}\label{chocice}
%\begin{tikzpicture}[Q,centerzero,scale=.9]
%\draw[-to] (-0.3,-.3) to (0.3,.3);
%\draw[wipe] (0.3,-.3) to (-0.3,.3);
%\draw[to-] (0.3,-.3) to (-0.3,.3);
%\end{tikzpicture}
%&:=
%\begin{tikzpicture}[Q,centerzero,scale=.7]
%\draw[-to] (0.3,-.5) to (-0.3,.5);
%\draw[wipe] (-0.2,-.2) to (0.2,.3);
%\draw[-] (-0.2,-.2) to (0.2,.3);
%\draw[-] (0.2,.3) to[out=50,in=180] (0.5,.5);
%\draw[-to] (0.5,.5) to[out=0,in=90] (0.8,-.5);
%\draw[-] (-0.2,-.2) to[out=230,in=0] (-0.6,-.5);
%\draw[-] (-0.6,-.5) to[out=180,in=-90] (-0.85,.5);
%\end{tikzpicture}
%\ ,&
\begin{tikzpicture}[Q,centerzero,scale=.9]
\draw[to-] (0.3,-.3) to (-0.3,.3);
\draw[wipe] (-0.3,-.3) to (0.3,.3);
\draw[-to] (-0.3,-.3) to (0.3,.3);
\end{tikzpicture}
&:=
\begin{tikzpicture}[Q,centerzero,scale=.7]
\draw[-] (-0.2,-.2) to (0.2,.3);
\draw[wipe] (0.3,-.5) to (-0.3,.5);
\draw[-to] (0.3,-.5) to (-0.3,.5);
\draw[-] (0.2,.3) to[out=50,in=180] (0.5,.5);
\draw[-to] (0.5,.5) to[out=0,in=90] (0.8,-.5);
\draw[-] (-0.2,-.2) to[out=230,in=0] (-0.6,-.5);
\draw[-] (-0.6,-.5) to[out=180,in=-90] (-0.85,.5);
\end{tikzpicture}
\ .
\end{align}
The other relations are as follows.
First, the affine Hecke algebra relations:
\begin{align}\label{AHA2}
\begin{tikzpicture}[Q,centerzero,scale=.9]
\draw[-to] (-0.3,-.3) to (0.3,.3);
\opendot{-0.17,-0.17};
\draw[wipe] (0.3,-.3) to (-0.3,.3);
\draw[-to] (0.3,-.3) to (-0.3,.3);
\end{tikzpicture}
&=
\begin{tikzpicture}[Q,centerzero,scale=.9]
\draw[-to] (0.3,-.3) to (-0.3,.3);
\draw[wipe] (-0.3,-.3) to (0.3,.3);
\draw[-to] (-0.3,-.3) to (0.3,.3);
\opendot{0.17,0.17};
\end{tikzpicture}
\ ,&
\begin{tikzpicture}[Q,centerzero,scale=.9]
\draw[-to] (0.3,-.3) to (-0.3,.3);
\draw[wipe] (-0.3,-.3) to (0.3,.3);
\draw[-to] (-0.3,-.3) to (0.3,.3);
\opendot{0.17,-0.17};
\end{tikzpicture}&=
\begin{tikzpicture}[Q,centerzero,scale=.9]
\draw[-to] (-0.3,-.3) to (0.3,.3);
\draw[wipe] (0.3,-.3) to (-0.3,.3);
\draw[-to] (0.3,-.3) to (-0.3,.3);
\opendot{-0.17,0.17};
\end{tikzpicture}
\ ,\\\label{AHA3}
\begin{tikzpicture}[Q,centerzero,scale=.9]
\draw[-to] (0.3,-.3) to (-0.3,.3);
\draw[wipe] (-0.3,-.3) to (0.3,.3);
\draw[-to] (-0.3,-.3) to (0.3,.3);
\end{tikzpicture}
-
\begin{tikzpicture}[Q,centerzero,scale=.9]
\draw[-to] (-0.3,-.3) to (0.3,.3);
\draw[wipe] (0.3,-.3) to (-0.3,.3);
\draw[-to] (0.3,-.3) to (-0.3,.3);
\end{tikzpicture}
&=
z\:
\begin{tikzpicture}[Q,centerzero,scale=.9]
\draw[-to] (0.18,-.3) to (0.18,.3);
\draw[-to] (-0.18,-.3) to (-0.18,.3);
\end{tikzpicture}
\ ,&
\begin{tikzpicture}[Q,centerzero,scale=.8]
\draw[-to] (0.45,-.6) to (-0.45,.6);
\draw[-] (0,-.6) to[out=90,in=-90] (-.45,0);
\draw[wipe] (-0.45,0) to[out=90,in=-90] (0,0.6);
\draw[-to] (-0.45,0) to[out=90,in=-90] (0,0.6);
\draw[wipe] (0.45,.6) to (-0.45,-.6);
\draw[to-] (0.45,.6) to (-0.45,-.6);
\end{tikzpicture}
&=
\begin{tikzpicture}[Q,centerzero,scale=.8]
\draw[-to] (0.45,-.6) to (-0.45,.6);
\draw[wipe] (0,-.6) to[out=90,in=-90] (.45,0);
\draw[-] (0,-.6) to[out=90,in=-90] (.45,0);
\draw[-to] (0.45,0) to[out=90,in=-90] (0,0.6);
\draw[wipe] (0.45,.6) to (-0.45,-.6);
\draw[to-] (0.45,.6) to (-0.45,-.6);
\end{tikzpicture}
\ .
\end{align}
Next, the \emph{right adjunction relations}
implying that $F$ is right dual to $E$:
\begin{align}\label{Qrightadj}
\begin{tikzpicture}[Q,centerzero,scale=1.2]
\draw[-to] (-0.3,0.4) -- (-0.3,0) arc(180:360:0.15) arc(180:0:0.15) -- (0.3,-0.4);
\end{tikzpicture}
&=
\begin{tikzpicture}[Q,centerzero,scale=1.2]
\draw[to-] (0,-0.4) -- (0,0.4);
\end{tikzpicture}
\ ,&
\begin{tikzpicture}[Q,centerzero,scale=1.2]
\draw[-to] (-0.3,-0.4) -- (-0.3,0) arc(180:0:0.15) arc(180:360:0.15) -- (0.3,0.4);
\end{tikzpicture}
&=
\begin{tikzpicture}[Q,centerzero,scale=1.2]
\draw[-to] (0,-0.4) -- (0,0.4);
\end{tikzpicture}\ .
\end{align}
There is an
\emph{inversion relation} asserting that the 
morphism
\begin{equation}
\begin{dcases}
\begin{pmatrix}
\begin{tikzpicture}[Q,anchorbase]
\draw[to-] (0.25,-0.25) -- (-0.25,0.25);
\draw[wipe] (-0.25,-0.25) -- (0.25,0.25);
\draw[-to] (-0.25,-0.25) -- (0.25,0.25);
\end{tikzpicture} &
\begin{tikzpicture}[Q,anchorbase]
\draw[-to] (-0.25,0.15) to[out=-90,in=-90,looseness=3] (0.25,0.15);
\node at (0,.2) {$\phantom.$};\node at (0,-.3) {$\phantom.$};
\end{tikzpicture}
&
\begin{tikzpicture}[Q,anchorbase]
\draw[-to] (-0.25,0.15) to[out=-90,in=-90,looseness=3] (0.25,0.15);
\node at (0,.2) {$\phantom.$};\node at (0,-.3) {$\phantom.$};
\opendot{0.23,-0.03};
\end{tikzpicture}
&\!\!\cdots\!\!
&
\begin{tikzpicture}[Q,anchorbase]
\draw[-to] (-0.25,0.15) to[out=-90,in=-90,looseness=3] (0.25,0.15);
\node at (0,.2) {$\phantom.$};\node at (0,-.3) {$\phantom.$};
\multopendot{0.23,-0.03}{west}{-\kappa-1};
\end{tikzpicture}
\end{pmatrix}\phantom{_T}:E \otimes F \oplus \one^{\oplus (-\kappa)}
\rightarrow F \otimes E
&\text{if } \kappa \leq 0\\
\begin{pmatrix}   
\begin{tikzpicture}[Q,centerzero]
\draw[to-] (0.25,-0.25)  -- (-0.25,0.25);
\draw[wipe] (-0.25,-0.25) -- (0.25,0.25);
\draw[-to] (-0.25,-0.25) -- (0.25,0.25);
\end{tikzpicture} &
\begin{tikzpicture}[Q,centerzero]
\draw[-to] (-0.25,-0.15) to [out=90,in=90,looseness=3](0.25,-0.15);
\node at (0,.3) {$\phantom.$};
\node at (0,-.4) {$\phantom.$};
\end{tikzpicture}
&
\begin{tikzpicture}[Q,centerzero]
\draw[-to] (-0.25,-0.15) to [out=90,in=90,looseness=3](0.25,-0.15);
\node at (0,.3) {$\phantom.$};
\node at (0,-.4) {$\phantom.$};
\opendot{-0.23,.03};
\end{tikzpicture}
&
\!\!\!\cdots\!\!\!
&
\begin{tikzpicture}[Q,centerzero]
\draw[-to] (-0.25,-0.15) to [out=90,in=90,looseness=3](0.25,-0.15);
\node at (0,.3) {$\phantom.$};
\node at (0,-.4) {$\phantom.$};
\multopendot{-0.23,.03}{east}{\kappa-1};
\end{tikzpicture}\ 
\end{pmatrix}^\transpose\phantom{_T}:E \otimes F \rightarrow 
F \otimes E \oplus \one^{\oplus \kappa}
&\text{if } \kappa > 0
\end{dcases}\label{timmy}
\end{equation}
is invertible.
Then there is one more {\em bubble relation} which had no counterpart in the degenerate case:
\begin{itemize}
\item
If $\kappa < 0$, we let
$\begin{tikzpicture}[Q,centerzero,scale=1.2]
\draw[to-] (-0.2,-.2) to (0.2,.2);
\draw[wipe] (0.2,-.2) to (-0.2,.2);
\draw[-to] (0.2,-.2) to (-0.2,.2);
\end{tikzpicture}$
and $t^{-1} z\; \begin{tikzpicture}[Q,centerzero,scale=1]
\draw[to-] (-0.25,-0.15)  to [out=90,in=90,looseness=3](0.25,-0.15);
\opendot{.23,0.05};
\end{tikzpicture}\:$
be the first and last entries of 
the inverse of the matrix \cref{timmy}, respectively,
define $\ \begin{tikzpicture}[Q,centerzero,scale=1]
\draw[to-] (-0.25,0.15)  to [out=-90,in=-90,looseness=3](0.25,0.15);
\end{tikzpicture}\ := t^{-1}
\begin{tikzpicture}[Q,centerzero,scale=1.2]
\draw[-] (.15,-.15) to [out=-90,in=-90,looseness=1.7] (-.15,-.15);
\draw[-] (-.15,-.15) to [out=90,in=-90,looseness=1] (.2,.3);
\draw[wipe] (-.2,.3) to [out=-90,in=90,looseness=1] (.15,-.15);
\draw[to-] (-.2,.3) to [out=-90,in=90,looseness=1] (.15,-.15);
\multopendot{0.14,-0.15}{west}{-\kappa};
\end{tikzpicture}$,
and require that
$\ \begin{tikzpicture}[Q,anchorbase,scale=1.2]
\draw[-to] (0,0.4) to[out=180,in=90] (-.2,0.2);
\draw[-] (0.2,0.2) to[out=90,in=0] (0,.4);
\draw[-] (-.2,0.2) to[out=-90,in=180] (0,0);
\draw[-] (0,0) to[out=0,in=-90] (0.2,0.2);
\end{tikzpicture}
=-t^{-1} z^{-1} \id_\one$.
\item If $\kappa = 0$, 
we let
$\begin{tikzpicture}[Q,centerzero,scale=1.2]
\draw[to-] (-0.2,-.2) to (0.2,.2);
\draw[wipe] (0.2,-.2) to (-0.2,.2);
\draw[-to] (0.2,-.2) to (-0.2,.2);
\end{tikzpicture}
 := 
\Big(\begin{tikzpicture}[Q,centerzero,scale=1.2]
\draw[to-] (0.2,-.2) to (-0.2,.2);
\draw[wipe] (-0.2,-.2) to (0.2,.2);
\draw[-to] (-0.2,-.2) to (0.2,.2);
\end{tikzpicture}
\Big)^{-1}$,
$\ \begin{tikzpicture}[Q,centerzero,scale=1]
\draw[to-] (-0.25,-0.15)  to [out=90,in=90,looseness=3](0.25,-0.15);
\end{tikzpicture}\ := t^{-1}
\begin{tikzpicture}[Q,centerzero,scale=1.2]
\draw[to-] (-.2,-.3) to [out=90,in=-90,looseness=1] (.15,.15);
\draw[-] (.15,.15) to [out=90,in=90,looseness=1.7] (-.15,.15);
\draw[wipe] (-.15,.15) to [out=-90,in=90,looseness=1] (.2,-.3);
\draw[-] (-.15,.15) to [out=-90,in=90,looseness=1] (.2,-.3);
\end{tikzpicture}$,
$\ \begin{tikzpicture}[Q,centerzero,scale=1]
\draw[to-] (-0.25,0.15)  to [out=-90,in=-90,looseness=3](0.25,0.15);
\end{tikzpicture}\ :=t^{-1}
\begin{tikzpicture}[Q,centerzero,scale=1.2]
\draw[-] (.15,-.15) to [out=-90,in=-90,looseness=1.7] (-.15,-.15);
\draw[-] (-.15,-.15) to [out=90,in=-90,looseness=1] (.2,.3);
\draw[wipe] (-.2,.3) to [out=-90,in=90,looseness=1] (.15,-.15);
\draw[to-] (-.2,.3) to [out=-90,in=90,looseness=1] (.15,-.15);
\end{tikzpicture}$, and require that
$\begin{tikzpicture}[Q,anchorbase,scale=1.2]
\draw[-] (0,0.4) to[out=180,in=90] (-.2,0.2);
\draw[-] (0.2,0.2) to[out=90,in=0] (0,.4);
\draw[to-] (-.2,0.2) to[out=-90,in=180] (0,0);
\draw[-] (0,0) to[out=0,in=-90] (0.2,0.2);
\end{tikzpicture}
=(tz^{-1}-t^{-1}z^{-1}) \id_\one$.
\item If $\kappa > 0$, we let 
$\begin{tikzpicture}[Q,centerzero,scale=1.2]
\draw[-to] (0.2,-.2) to (-0.2,.2);
\draw[wipe] (-0.2,-.2) to (0.2,.2);
\draw[to-] (-0.2,-.2) to (0.2,.2);
\end{tikzpicture}$
and
$t^{-1} z\ \begin{tikzpicture}[Q,centerzero,scale=1]
\draw[to-] (-0.25,0.15)  to [out=-90,in=-90,looseness=3](0.25,0.15);
\end{tikzpicture}\ $
be the first and {\em second} entries of the inverse of the matrix \cref{timmy}, respectively,
define
$\ \begin{tikzpicture}[Q,centerzero,scale=1]
\draw[to-] (-0.25,-0.15)  to [out=90,in=90,looseness=3](0.25,-0.15);
\end{tikzpicture}\ := t
\begin{tikzpicture}[Q,centerzero,scale=1.2]
\draw[-] (.15,.15) to [out=90,in=90,looseness=1.7] (-.15,.15);
\draw[-] (-.15,.15) to [out=-90,in=90,looseness=1] (.2,-.3);
\draw[wipe] (-.2,-.3) to [out=90,in=-90,looseness=1] (.15,.15);
\draw[to-] (-.2,-.3) to [out=90,in=-90,looseness=1] (.15,.15);
\multopendot{-0.14,0.15}{east}{\kappa};
\end{tikzpicture}$,
and require that
$\begin{tikzpicture}[Q,anchorbase,scale=1.2]
\draw[to-] (0,0.4) to[out=180,in=90] (-.2,0.2);
\draw[-] (0.2,0.2) to[out=90,in=0] (0,.4);
\draw[-] (-.2,0.2) to[out=-90,in=180] (0,0);
\draw[-] (0,0) to[out=0,in=-90] (0.2,0.2);
\multopendot{-.2,.2}{east}{\kappa};
\end{tikzpicture}
=-t^{-1}z^{-1} \id_\one$.
\end{itemize}
The leftward cup and cap just introduced automatically satisfy
\begin{align}\label{Qleftpivots}
\begin{tikzpicture}[Q,centerzero,scale=1.2]
\draw[to-] (-0.3,0.4) -- (-0.3,0) arc(180:360:0.15) arc(180:0:0.15) -- (0.3,-0.4);
\end{tikzpicture}
&=
\begin{tikzpicture}[Q,centerzero,scale=1.2]
\draw[-to] (0,-0.4) -- (0,0.4);
\end{tikzpicture}
\ ,&
\begin{tikzpicture}[Q,centerzero,scale=1.2]
\draw[to-] (-0.3,-0.4) -- (-0.3,0) arc(180:0:0.15) arc(180:360:0.15) -- (0.3,0.4);
\end{tikzpicture}
&=
\begin{tikzpicture}[Q,centerzero,scale=1.2]
\draw[to-] (0,-0.4) -- (0,0.4);
\end{tikzpicture}\ ,
\end{align}
so they define a second adjunction making $F$ into a left dual of $E$.
There is also a downward dot, a positive rightward crossing, another leftward crossing, and positive and negative downward crossings, all defined so that the resulting string diagrams are invariant under planar isotopy. Hence, $\qHeis_\kappa$ is strictly pivotal.
Proofs of these results, plus some other equivalent (and more symmetric) formulations of the definition, can be found in \cite{BSWquantum}.

Now that dots are invertible, we find it better to use a different
generating function for positive powers of dots, setting
\begin{align}\label{Qdgf}
\begin{tikzpicture}[Q,centerzero,scale=1.1]
\draw[-] (0,-0.3) -- (0,0.3);
\squared{0,0}{u};
\end{tikzpicture}
&:=
\begin{tikzpicture}[Q,centerzero,scale=1.1]
\draw[-] (0,-0.3) -- (0,0.3);
\pin{0,0}{.7,0}{\frac{x}{u-x}};
\end{tikzpicture}\
=
\sum_{r > 0}
\begin{tikzpicture}[Q,centerzero,scale=1.1]
\draw[-] (0,-0.3) -- (0,0.3);
\multopendot{0,0}{east}{r};
\end{tikzpicture}\ u^{-r}.
\end{align}
The bubble generating functions in this setting are the unique formal Laurent series
\begin{align}
\label{Qbubblegeneratingfunction1}
\begin{tikzpicture}[Q,baseline=-1mm]
\draw[to-] (-0.25,0) arc(180:-180:0.25);
\node at (.55,0) {\color{black}$(u)$};
\end{tikzpicture}
&\in 
tz^{-1} u^{\kappa}\id_{\one}
+ u^{\kappa-1} \kk\llbracket u^{-1}\rrbracket\End_{\qHeis_\kappa}(\one),\\
\begin{tikzpicture}[Q,baseline=-1mm]
\draw[-to] (-0.25,0) arc(180:-180:0.25);
\node at (.55,0) {$\color{black}(u)$};
\end{tikzpicture}
&\in -t^{-1} z^{-1} u^{-\kappa}\id_{\one}
+ u^{-\kappa-1} \kk\llbracket u^{-1}\rrbracket\End_{\qHeis_\kappa}(\one)
\label{Qbubblegeneratingfunction2}
\end{align}
such that
\begin{align}
\left[\ \begin{tikzpicture}[Q,baseline=-1mm]
\draw[to-] (-0.25,0) arc(180:-180:0.25);
\node at (.55,0) {$\color{black}(u)$};
\end{tikzpicture}\right]_{u:<0}\!\!
&=
\begin{tikzpicture}[Q,baseline=-1mm]
\draw[to-] (-0.25,0) arc(180:-180:0.25);
\squared{.25,0}{u};
\end{tikzpicture},
&
\left[\ \begin{tikzpicture}[Q,baseline=-1mm]
\draw[-to] (-0.25,0) arc(180:-180:0.25);
\node at (.55,0) {$\color{black}(u)$};
\end{tikzpicture}\right]_{u:<0}\!\!
&=
\begin{tikzpicture}[Q,baseline=-1mm]
\draw[-to] (0.25,0) arc(360:0:0.25);
\squared{-.25,0}{u};
\end{tikzpicture},&
\label{Qinfgrass}
\begin{tikzpicture}[Q,centerzero,scale=1]
\draw[to-] (-0.68,0) arc(180:-180:0.25);
\node[black] at (0.1,0) {$\color{black}(u)$};
\end{tikzpicture} 
\begin{tikzpicture}[Q,centerzero,scale=1]
\draw[-to] (-.25,0) arc(180:-180:0.25);
\node[black] at (0.54,0) {$\color{black}(u)$};
\end{tikzpicture} &= -z^{-2} \id_{\one}.
\end{align}
Using the same notation \cref{dgf} for the generating function for non-negative powers of dots, the counterparts of the relations \cref{Hbubslide,Hcurlrels,Haltquadratic} are
\begin{align}\label{Qbubslide}
\begin{tikzpicture}[centerzero,Q,scale=.9]
\draw[-to] (0.08,-.4) to (0.08,.4);
\draw[to-] (0.4,0) arc(180:-180:0.25);
\node at (1.2,0) {$\color{black}(u)$};
\end{tikzpicture}&=
\begin{tikzpicture}[centerzero,Q,scale=.9]
\draw[-to] (0.15,-.4) to (0.15,.4);
\draw[to-] (-1.1,0) arc(180:-180:0.25);
\pin{.15,0}{1.2,0}{R_q(u,x)};
\node at (-.3,0) {$\color{black}(u)$};
\end{tikzpicture}\ ,
&
\begin{tikzpicture}[centerzero,Q,scale=.9]
\draw[-to] (0.15,-.4) to (0.15,.4);
\draw[to-] (-1.1,0) arc(-180:180:0.25);
\node at (-.3,0) {$\color{black}(u)$};
\end{tikzpicture}
&=
\begin{tikzpicture}[centerzero,Q,scale=.9]
\draw[-to] (0.08,-.4) to (0.08,.4);
\draw[to-] (0.4,0) arc(-180:180:0.25);
\node at (1.2,0) {$\color{black}(u)$};
\pin{.08,0}{-1,0}{R_q(u,x)};
\end{tikzpicture}\\\intertext{where $R_q(x,y) := 1-z^2xy(x-y)^{-2} = \frac{(x-q^2y)(x-q^{-2}y)}{(x-y)^2}$,}
\begin{tikzpicture}[centerzero,Q]
\draw[to-] (0,0.6) to (0,0.3);
\draw[-] (0.3,-0.2) to [out=0,in=-90](.5,0);
\draw[-] (0.5,0) to [out=90,in=0](.3,0.2);
\draw[-] (0,0.3) to [out=-90,in=180] (.3,-0.2);
\draw[-] (0,-0.3) to (0,-0.6);
\draw[wipe] (0.3,.2) to [out=180,in=90](0,-0.3);
\draw[-] (0.3,.2) to [out=180,in=90](0,-0.3);
\squared{.5,0}{u};
\end{tikzpicture}
&=
-z\left[
\ \begin{tikzpicture}[centerzero,Q]
\draw[-to] (0.6,0.05) arc(180:-180:0.2);
\node at (1.28,0.05) {$\color{black}(u)$};
\draw[-to] (0.2,-.5) to (0.2,.6);
\squared{.2,.05}{u};
\end{tikzpicture}\ 
\right]_{u:<0}\!,&
\begin{tikzpicture}[centerzero,Q]
\draw[to-] (0,0.6) to (0,0.3);
\draw[-] (0,0.3) to [out=-90,in=0] (-.3,-0.2);
\draw[-] (0,-0.3) to (0,-0.6);
\draw[-] (-0.3,-0.2) to [out=180,in=-90](-.5,0);
\draw[-] (-0.5,0) to [out=90,in=180](-.3,0.2);
\draw[wipe] (-0.3,.2) to [out=0,in=90](0,-0.3);
\draw[-] (-0.3,.2) to [out=0,in=90](0,-0.3);
\squared{-.5,0}{u};
\end{tikzpicture}
&=
z\left[\ \begin{tikzpicture}[centerzero,Q]
\draw[-to] (-1,0.05) arc(-180:180:0.2);
\node at (-.33,0.05) {$\color{black}(u)$};
\draw[-to] (0.2,-.5) to (0.2,.6);
\squared{.2,.05}{u};
\end{tikzpicture}\ 
\right]_{u:<0}\!\!,
\label{Qcurlrels}\\
\begin{tikzpicture}[centerzero,Q,scale=.9]
\draw[-to] (0.25,0) to[out=90,in=-45] (-0.28,.6);
\draw[to-] (0.28,-.6) to[out=135,in=-90] (-0.25,0);
\draw[wipe] (-0.25,0) to[out=90,in=-135] (0.28,.6);
\draw[-] (-0.25,0) to[out=90,in=-135] (0.28,.6);
\draw[wipe] (-0.28,-.6) to[out=45,in=-90] (0.25,0);
\draw[-] (-0.28,-.6) to[out=45,in=-90] (0.25,0);
\end{tikzpicture}
&=
\begin{tikzpicture}[centerzero,Q,scale=.9]
\draw[to-] (0.08,-.6) to (0.08,.6);
\draw[-to] (-0.28,-.6) to (-0.28,.6);
\end{tikzpicture}
-t^{-1}z
\begin{tikzpicture}[centerzero,Q,scale=.9]
\draw[-] (0.3,0.6) to[out=-90, in=0] (0,0.1);
\draw[-to] (0,0.1) to[out = 180, in = -90] (-0.3,0.6);
\draw[to-] (0.3,-.6) to[out=90, in=0] (0,-0.1);
\draw[-] (0,-0.1) to[out = 180, in = 90] (-0.3,-.6);
\end{tikzpicture}
\!+z^2
\left[\ 
\begin{tikzpicture}[centerzero,Q,scale=.9]
\draw[-] (0.3,0.6) to[out=-90, in=0] (0,0.1);
\draw[-to] (0,0.1) to[out = 180, in = -90] (-0.3,0.6);
\squared{-.28,.3}{u};
\draw[to-] (0.3,-.6) to[out=90, in=0] (0,-0.1);
\draw[-] (0,-0.1) to[out = 180, in = 90] (-0.3,-.6);
\squared{-.28,-.3}{u};
\draw[-to] (.37,0.0) arc(-180:180:0.2);
\node at (1.07,0.0) {$\color{black}(u)$};
\end{tikzpicture}\right]_{u:0}\!\!,&
\begin{tikzpicture}[centerzero,Q,scale=.9]
\draw[to-] (-0.28,-.6) to [out=45,in=-90] (0.25,0);
\draw[wipe] (0.28,-.6) to[out=135,in=-90] (-0.25,0);
\draw[-] (0.28,-.6) to[out=135,in=-90] (-0.25,0);
\draw[-to] (-0.25,0) to[out=90,in=-135] (0.28,.6);
\draw[wipe] (0.25,0) to[out=90,in=-45] (-0.28,.6);
\draw[-] (0.25,0) to[out=90,in=-45] (-0.28,.6);
\end{tikzpicture}
&=
\begin{tikzpicture}[centerzero,Q,scale=.9]
\draw[-to] (0.08,-.6) to (0.08,.6);
\draw[to-] (-0.28,-.6) to (-0.28,.6);
\end{tikzpicture}
+tz
\begin{tikzpicture}[centerzero,Q,scale=.9]
\draw[to-] (0.3,0.6) to[out=-90, in=0] (0,.1);
\draw[-] (0,.1) to[out = 180, in = -90] (-0.3,0.6);
\draw[-] (0.3,-.6) to[out=90, in=0] (0,-0.1);
\draw[-to] (0,-0.1) to[out = 180, in = 90] (-0.3,-.6);
\end{tikzpicture}
\!+z^2
\left[\ \begin{tikzpicture}[centerzero,Q,scale=.9]
\draw[-] (-0.3,0.6) to[out=-90, in=180] (0,0.1);
\draw[-to] (0,0.1) to[out = 0, in = -90] (0.3,0.6);
\squared{.28,.3}{u};
\draw[to-] (-0.3,-.6) to[out=90, in=180] (0,-0.1);
\draw[-] (0,-0.1) to[out = 0, in = 90] (0.3,-.6);
\squared{.28,-.3}{u};
\draw[-to] (-1.26,0.0) arc(180:-180:0.2);
\node at (-.56,0.0) {$\color{black}(u)$};
\end{tikzpicture}\ \right]_{u:0}\!\!.
\label{Qaltquadratic}\end{align}

\begin{defin}\label{qhcatdef}
A {\em quantum Heisenberg categorification} of central charge $\kappa$ and $q \in \kk-\{0,\pm 1\}$
is a locally finite $\kk$-linear Abelian category $\catR$
plus an invertible element $t \in Z(\catR)^\times$ (usually, a scalar in $\kk^\times$) and
an adjoint pair $(E, F)$ of $\kk$-linear endofunctors
such that:
\begin{itemize}
\item[($q$-H1)]
The adjoint pair $(E,F)$ has a prescribed adjunction with unit 
and counit of adjunction denoted
$\;\begin{tikzpicture}[Q,baseline=-4pt,scale=.8]
\draw[-to] (-0.25,0.15) to[out=-90,in=-90,looseness=3] (0.25,0.15);
\end{tikzpicture}\;:\id_\catR \Rightarrow F \circ E$
and
$\;\begin{tikzpicture}[Q,baseline=-2pt,scale=.8]
\draw[-to] (-0.25,-0.15) to [out=90,in=90,looseness=3](0.25,-0.15);
\end{tikzpicture}\;:E\circ F \Rightarrow \id_\catR$.
\item[($q$-H2)]
There are given invertible natural transformations
$\begin{tikzpicture}[Q,centerzero]
\draw[-to] (0,-0.2) -- (0,0.2);
\opendot{0,0};
\end{tikzpicture}:E \Rightarrow E$ and $\begin{tikzpicture}[Q,centerzero,scale=.9]
\draw[-to] (0.2,-0.2) -- (-0.2,0.2);
\draw[wipe] (-0.2,-0.2) -- (0.2,0.2);
\draw[-to] (-0.2,-0.2) -- (0.2,0.2);\end{tikzpicture}:E \circ E \Rightarrow E \circ E$ satisfying the 
AHA relations \cref{AHA2,AHA3}, where
$\ \begin{tikzpicture}[Q,centerzero,scale=1]
\draw[-to] (-0.2,-0.2) -- (0.2,0.2);
\draw[wipe] (0.2,-0.2) -- (-0.2,0.2);
\draw[-to] (0.2,-0.2) -- (-0.2,0.2);
\end{tikzpicture}\  := \left(\ \begin{tikzpicture}[Q,centerzero,scale=1]
\draw[-to] (0.2,-0.2) -- (-0.2,0.2);
\draw[wipe] (-0.2,-0.2) -- (0.2,0.2);
\draw[-to] (-0.2,-0.2) -- (0.2,0.2);
\end{tikzpicture}\ 
\right)^{-1}$.
\item[($q$-H3)]
Defining the negative rightward crossing \begin{tikzpicture}[Q,centerzero,scale=1]
\draw[to-] (0.2,-0.2) -- (-0.2,0.2);
\draw[wipe] (-0.2,-0.2) -- (0.2,0.2);
\draw[-to] (-0.2,-0.2) -- (0.2,0.2);
\end{tikzpicture}
according to \cref{chocice},
the matrix
\begin{align*}
\begin{pmatrix}
\begin{tikzpicture}[Q,anchorbase]
\draw[to-] (0.25,-0.25) -- (-0.25,0.25);
\draw[wipe] (-0.25,-0.25) -- (0.25,0.25);
\draw[-to] (-0.25,-0.25) -- (0.25,0.25);
\end{tikzpicture} &
\begin{tikzpicture}[Q,anchorbase]
\draw[-to] (-0.25,0.15) to[out=-90,in=-90,looseness=3] (0.25,0.15);
\node at (0,.2) {$\phantom.$};\node at (0,-.3) {$\phantom.$};
\end{tikzpicture}
&
\begin{tikzpicture}[Q,anchorbase]
\draw[-to] (-0.25,0.15) to[out=-90,in=-90,looseness=3] (0.25,0.15);
\node at (0,.2) {$\phantom.$};\node at (0,-.3) {$\phantom.$};
\opendot{0.23,-0.03};
\end{tikzpicture}
&\!\!\cdots\!\!
&
\begin{tikzpicture}[Q,anchorbase]
\draw[-to] (-0.25,0.15) to[out=-90,in=-90,looseness=3] (0.25,0.15);
\node at (0,.2) {$\phantom.$};\node at (0,-.3) {$\phantom.$};
\multopendot{0.23,-0.03}{west}{-\kappa-1};
\end{tikzpicture}
\end{pmatrix}\phantom{_T}&:E \circ F \oplus \id_\catR^{\oplus (-\kappa)}
\Rightarrow F \circ E
&&\text{if } \kappa \leq 0,\text{ or}\\
\begin{pmatrix}   
\begin{tikzpicture}[Q,centerzero]
\draw[to-] (0.25,-0.25)  -- (-0.25,0.25);
\draw[wipe] (-0.25,-0.25) -- (0.25,0.25);
\draw[-to] (-0.25,-0.25) -- (0.25,0.25);
\end{tikzpicture} &
\begin{tikzpicture}[Q,centerzero]
\draw[-to] (-0.25,-0.15) to [out=90,in=90,looseness=3](0.25,-0.15);
\node at (0,.3) {$\phantom.$};
\node at (0,-.4) {$\phantom.$};
\end{tikzpicture}
&
\begin{tikzpicture}[Q,centerzero]
\draw[-to] (-0.25,-0.15) to [out=90,in=90,looseness=3](0.25,-0.15);
\node at (0,.3) {$\phantom.$};
\node at (0,-.4) {$\phantom.$};
\opendot{-0.23,.03};
\end{tikzpicture}
&
\!\!\!\cdots\!\!\!
&
\begin{tikzpicture}[Q,centerzero]
\draw[-to] (-0.25,-0.15) to [out=90,in=90,looseness=3](0.25,-0.15);
\node at (0,.3) {$\phantom.$};
\node at (0,-.4) {$\phantom.$};
\multopendot{-0.23,.03}{east}{\kappa-1};
\end{tikzpicture}\ 
\end{pmatrix}^\transpose\phantom{_T}&:E \circ F \Rightarrow 
F \circ E \oplus \id_{\catR}^{\oplus \kappa}
&&\text{if } \kappa > 0,
\end{align*}
is
invertible.
\item[($q$-H4)]
Defining
$\;\begin{tikzpicture}[Q,baseline=-4pt,scale=.8]
\draw[to-] (-0.25,0.15) to[out=-90,in=-90,looseness=3] (0.25,0.15);
\end{tikzpicture}\;:\id_\catR \Rightarrow E \circ F$
and
$\;\begin{tikzpicture}[Q,baseline=-2pt,scale=.8]
\draw[to-] (-0.25,-0.15) to [out=90,in=90,looseness=3](0.25,-0.15);
\end{tikzpicture}:F\circ E \Rightarrow \id_\catR$ as explained after \cref{timmy},
we have that
$\ \begin{tikzpicture}[Q,anchorbase,scale=1.2]
\draw[-to] (0,0.4) to[out=180,in=90] (-.2,0.2);
\draw[-] (0.2,0.2) to[out=90,in=0] (0,.4);
\draw[-] (-.2,0.2) to[out=-90,in=180] (0,0);
\draw[-] (0,0) to[out=0,in=-90] (0.2,0.2);
\end{tikzpicture}
=-t^{-1} z^{-1}$
if $\kappa < 0$,
$\ \begin{tikzpicture}[Q,anchorbase,scale=1.2]
\draw[-] (0,0.4) to[out=180,in=90] (-.2,0.2);
\draw[-] (0.2,0.2) to[out=90,in=0] (0,.4);
\draw[to-] (-.2,0.2) to[out=-90,in=180] (0,0);
\draw[-] (0,0) to[out=0,in=-90] (0.2,0.2);
\end{tikzpicture}
=tz^{-1}-t^{-1} z^{-1}$ if $\kappa = 0$,
and 
$\begin{tikzpicture}[Q,anchorbase,scale=1.2]
\draw[to-] (0,0.4) to[out=180,in=90] (-.2,0.2);
\draw[-] (0.2,0.2) to[out=90,in=0] (0,.4);
\draw[-] (-.2,0.2) to[out=-90,in=180] (0,0);
\draw[-] (0,0) to[out=0,in=-90] (0.2,0.2);
\multopendot{-.2,.2}{east}{\kappa};
\end{tikzpicture}
=t z^{-1}$ if $\kappa > 0$.
\end{itemize}
This data makes 
$\catR$ is a strict left $\qHeis_\kappa$-module category.
\end{defin}

Next we summarize the quantum Heisenberg to Kac-Moody construction from \cite[Sec.~4]{HKM}.
Let $\catR$ be a quantum Heisenberg categorification of central charge $\kappa$. We fix representatives $L(b)\:(b \in \B)$ for the isomorphism classes of irreducible objects in $\catR$,
and define the (monic) minimal polynomials $m_b(x), n_b(x) \in \kk[x]$ and the central character generating function $\chi_b(u) \in \kk\lround u^{-1}\rround$ for $b \in \B$ exactly as we did in the degenerate case in \cref{s4-HtoKM} (using the new bubble generating function \cref{Qbubblegeneratingfunction1} for $\chi_b(u)$). We also have the scalar $t_b \in \kk$ which is how the central endomorphism 
$t_L$ acts on $L$.
The analogue of \cref{mainfact} proved in \cite[Lem.~4.4]{HKM} is that
\begin{equation}\label{mainfuct}
\chi_b(u) = t_b z^{-1} n_b(u) / m_b(u).
\end{equation}
Moreover, 
$t_b^2 = m_b(0) / n_b(0)$ (this makes sense because $m_b(0)$ and $n_b(0)$ are non-zero in the quantum case as the dot is invertible).

The {\em spectrum} $I$ is 
the set of roots of all $m_b(x)$. It is a subset of $\kk^\times$ 
closed under the operations $i \mapsto q^{\pm 2} i$. 
The Cartan matrix $A = (a_{i,j})_{i,j \in I}$ is 
\begin{equation}\label{Qcartanmx}
a_{i,j} := \begin{cases}
2&\text{if $i = j$}\\
-1&\text{if $i = q^{\pm 2} j$ and $q^4 \neq 1$}\\
-2&\text{if $i = q^2 j = q^{-2} j$}\\
0&\text{otherwise.}
\end{cases}
\end{equation}
The connected components of $A$ are all of type $A_\infty$ if 
$q$ is not a root of unity
or $A_{p-1}^{(1)}$ if $q^2$ is a primitive $p$th root of 1 (note now that $p$ is not necessarily prime).
We also fix a root datum of this Cartan type in the same way as in \cref{s4-HtoKM}.
Then we define $\eps_i(b)$ and $\phi_i(b)$ as in \cref{epsphi},
and $\wt:\catB \rightarrow X$ as in \cref{wtdef}.
We obtain the weight decomposition \cref{blockdec}, letting
$\catR_\lambda$ be the Serre subcategory of $\catR$ generated by the irreducible objects $L(b)$ for $b \in \B_\lambda
:= \{b \in \B\:|\:\wt(b) = \lambda\}$ as before.
The eigenfunctors $E_i, F_i\:(i \in I)$ are defined as in \cref{eigenfunctors}.
Again, \cite[Lem.~4.7]{HKM}, which depends on
\cref{Qbubslide}, implies that
the restrictions of $E_i$ and $F_i$ are functors
\begin{align}\label{Qfirstaxiom}
E_i|_{\catR_\lambda}&:\catR_\lambda \rightarrow \catR_{\lambda+\alpha_i},
&F_i|_{\catR_{\lambda+\alpha_i}}&:\catR_{\lambda+\alpha_i} \rightarrow \catR_{\lambda}.
\end{align}

It remains to define the various Kac-Moody natural transformations.
We introduce the projectors as in \cref{projectors},
and use these to 
define projected rightward cups and caps, dots and crossings.
The only difference is that there are now {\em two} projected crossings, positive and negative, for each orientation. For example, the positive and negative upward projected crossings are 
\begin{align}\label{pully}
\begin{tikzpicture}[centerzero,KM]
\draw[-to] (0.28,-.28) \botlabel{j} to (-0.28,.28)\toplabel{i'};
\draw[line width=5.5pt,white] (-0.28,-.28) to (0.28,.28);
\draw[-to] (-0.28,-.28) \botlabel{i} to (0.28,.28)\toplabel{j'};
\projcr{0,0};
\end{tikzpicture}
&:=
\begin{tikzpicture}[centerzero,KM]
\draw[-to] (0.28,-.28) \botlabel{j} to (-0.28,.28)\toplabel{i'};
\draw[-to] (-0.28,-.28) \botlabel{i} to (0.28,.28)\toplabel{j'};
\draw[Q] (0.13,-.13) to (-0.13,.13);
\draw[wipe] (-0.13,-.13) to (0.13,.13);
\draw[Q] (-0.13,-.13) to (0.13,.13);
\notch[45]{.13,-.13};
\notch[45]{-.13,.13};
\notch[-45]{.13,.13};
\notch[-45]{-.13,-.13};
\end{tikzpicture}\ ,&
\begin{tikzpicture}[centerzero,KM]
\draw[-to] (-0.28,-.28) \botlabel{i} to (0.28,.28)\toplabel{j'};
\draw[line width=5.5pt,white] (0.28,-.28) to (-0.28,.28);
\draw[-to] (0.28,-.28) \botlabel{j} to (-0.28,.28)\toplabel{i'};
\projcr{0,0};
\end{tikzpicture}
&:=
\begin{tikzpicture}[centerzero,KM]
\draw[-to] (-0.28,-.28) \botlabel{i} to (0.28,.28)\toplabel{j'};
\draw[wipe] (0.28,-.28) to (-0.28,.28);
\draw[-to] (0.28,-.28) \botlabel{j} to (-0.28,.28)\toplabel{i'};
\draw[Q] (-0.13,-.13) to (0.13,.13);
\draw[wipe] (0.13,-.13) to (-0.13,.13);
\draw[Q] (0.13,-.13) to (-0.13,.13);
\notch[45]{.13,-.13};
\notch[45]{-.13,.13};
\notch[-45]{.13,.13};
\notch[-45]{-.13,-.13};
\end{tikzpicture}\ .
\end{align}
These have the same vanishing properties as before.
Then we define the Kac-Moody rightward cups and caps, the Kac-Moody dots (which are nilpotent), and the upward Kac-Moody crossings by
\begin{align}
\begin{tikzpicture}[KM,centerzero]
\draw[-to] (-0.25,0.15) \toplabel{i} to[out=-90,in=-90,looseness=2.5] (0.25,0.15);
\end{tikzpicture}\ &:=\begin{tikzpicture}[KM,anchorbase,scale=1.1]
\draw[to-] (0.3,0.4)\toplabel{i} to (0.3,.2);
\draw[-] (-0.1,.2) to (-0.1,0.4)\toplabel{i};
\draw[Q] (.3,.2) to[out=-90,in=-90,looseness=1.5] (-.1,.2);
\notch{-.1,.2};
\notch{.3,.2};
\end{tikzpicture}\ ,&
\begin{tikzpicture}[KM,centerzero]
\draw[-to] (-0.25,-0.15) \botlabel{i} to [out=90,in=90,looseness=2.5](0.25,-0.15);
\end{tikzpicture}\ 
&:=
\begin{tikzpicture}[KM,anchorbase,scale=1.1]
\draw[to-] (0.3,-0.4)\botlabel{i} to (0.3,-.2);
\draw[-] (-0.1,-.2) to (-0.1,-0.4)\botlabel{i};
\draw[Q] (.3,-.2) to[out=90,in=90,looseness=1.5] (-.1,-.2);
\notch{-.1,-.2};
\notch{.3,-.2};
\end{tikzpicture}\ ,\label{Qgofigure}\\
\label{Qgoprodots}
\begin{tikzpicture}[KM,centerzero]
\draw[-to] (0.08,-.35) \botlabel{i} to (.08,.35);
\opendot{0.08,0};
\end{tikzpicture}
&:=
\begin{tikzpicture}[Q,centerzero]
\draw[KM,-to] (0.08,-.35) \botlabel{i} to (.08,.35);
\draw (.08,-.15) to (0.08,.15);
\notch{0.08,.15};
\notch{0.08,-.15};
\pin{0.08,0}{.7,0}{\frac{x}{i}-1};
\end{tikzpicture}\ ,&
\begin{tikzpicture}[KM,centerzero]
\draw[to-] (0.08,-.35) \botlabel{i} to (.08,.35);
\opendot{0.08,0};
\end{tikzpicture}
&:=
\begin{tikzpicture}[Q,centerzero]
\draw[KM,to-] (0.08,-.35) \botlabel{i} to (.08,.35);
\draw (.08,-.15) to (0.08,.15);
\notch{0.08,.15};
\notch{0.08,-.15};
\pin{0.08,0}{.7,0}{\frac{x}{i}-1};
\end{tikzpicture}\ ,\end{align}\begin{align}
\begin{tikzpicture}[KM,centerzero,scale=.9]
\draw[-to] (-0.3,-0.3) \botlabel{i} -- (0.3,0.3);
\draw[-to] (0.3,-0.3) \botlabel{j} -- (-0.3,0.3);
\end{tikzpicture}
:=
\begin{dcases}
\begin{tikzpicture}[KM,centerzero,scale=.9]
\draw[-to] (0.4,-.4)\botlabel{i} to (-0.4,.4)\toplabel{i};
\draw[line width=5.5,white] (-0.4,-.4) to (0.4,.4);
\draw[-to] (-0.4,-.4)\botlabel{i} to (0.4,.4)\toplabel{i};
\projcr{0,0};
\pinpin{-.25,-.25}{.25,-.25}{2.3,-.25}{(q-q^{-1} + qx-q^{-1}y)^{-1}};
\end{tikzpicture}
+q^{-1}\begin{tikzpicture}[KM,centerzero,scale=.9]
\draw[-to] (0.5,-.4)\botlabel{i} to (0.5,.4);
\draw[-to] (0.1,-.4)\botlabel{i} to (0.1,.4);
\pinpin{.1,0}{.5,0}{2.5,0}{(q-q^{-1}+qx-q^{-1}y)^{-1}};
\end{tikzpicture}
&\text{if $i=j$}\\
\begin{tikzpicture}[KM,centerzero,scale=.9]
\draw[-to] (0.4,-.4)\botlabel{q^{-2} i} to (-0.4,.4)\toplabel{q^{-2} i};
\draw[line width=5.5,white] (-0.4,-.4) to (0.4,.4);
\draw[-to] (-0.4,-.4)\botlabel{i} to (0.4,.4)\toplabel{i};
\projcr{0,0};
\pinpin{-.25,-.25}{.25,-.25}{2.2,-.25}{q-q^{-1}+qx-q^{-1}y};
\end{tikzpicture}
&\text{if $i=q^2 j$}\\
-\ \begin{tikzpicture}[KM,centerzero,scale=.9]
\draw[-to] (0.4,-.4) \botlabel{j} to (-0.4,.4)\toplabel{j};
\draw[line width=5.5,white] (-0.4,-.4) to (0.4,.4);
\draw[-to] (-0.4,-.4) \botlabel{i} to (0.4,.4)\toplabel{i};
\projcr{0,0};
\pinpin{-.25,-.25}{.25,-.25}{3.6,-.25}{(i-j+ix-jy)(q^{-1}i-qj+q^{-1}ix-qjy)^{-1}};
\end{tikzpicture}
&\text{if $i \neq j, q^2 j$.}
\end{dcases}
\end{align}
when $i \neq j$.
We refer to \cite[Sec.~4]{HKM} for a more detailed development of the  properties of these natural transformations.

\begin{theo}[{\cite[Th.~4.11]{HKM}}]\label{Qmaintheorem}
Let $\catR$ be a quantum Heisenberg categorification.
The definitions summarized in the previous two paragraphs make
$\catR$ into a Kac-Moody categorification with the parameters defined from
$$
Q_{i,j}(x,y) := 
\begin{cases}
0&\text{if $i=j$}\\
y-x&\text{if $i = q^2 j \neq q^{-2} j$}\\
x-y&\text{if $i = q^{-2} j \neq q^2 j$}\\
(y-x)(x-y)&\text{if $i = q^2 j = q^{-2} j$}\\
1&\text{if $i \neq j, q^{\pm 2} j$.}
\end{cases}
$$
\end{theo}

We close by writing down the definition of internal bubbles in the quantum case and formulating the counterpart of \cref{average}. It is quite similar 
to the degenerate case from the previous section, so we will be brief.
\iffalse
Replacing $u$ by $\frac{u}{i}-1$ in \cref{bgf} gives
\begin{align}\label{sbgf-q1}
    \begin{tikzpicture}[KM,baseline=-1mm,scale=.9]
        \draw[to-] (-0.25,0) arc(180:-180:0.25);
        \node at (0,-.4) {\strandlabel{i}};
        \region{1.8,0}{\lambda};
        \node at (1,0) {$(\frac{u}{i}-1)$};
    \end{tikzpicture}
    &=
    \begin{tikzpicture}[KM,baseline=-1mm]
        \stretchedanticlockwisebubble{0,0}{u};
        \bubblelabel{0,0}{\frac{u}{i}-1};
        \node at (0,-.35) {\strandlabel{i}};
        \region{0.5,0}{\lambda};
    \end{tikzpicture}
    +
    \begin{tikzpicture}[KM,baseline=-1mm,scale=.9]
        \draw[to-] (-0.25,0) arc(180:-180:0.25);
        \node at (0,-.4) {\strandlabel{i}};
        \pin{.25,0}{1.7,0}{i(u-ix-i)^{-1}};
        \region{3.1,0}{\lambda};
    \end{tikzpicture},\\\label{sbgf-q2}
    \begin{tikzpicture}[KM,baseline=-1mm,scale=.9]
        \draw[to-] (-0.25,0) arc(-180:180:0.25);
        \node at (0,-.4) {\strandlabel{i}};
        \region{1.8,0}{\lambda};
        \node at (1,0) {$(\frac{u}{i}-1)$};
    \end{tikzpicture}
    &=
    \begin{tikzpicture}[KM,baseline=-1mm]
        \stretchedclockwisebubble{0,0}{u};
        \bubblelabel{0,0}{\frac{u}{i}-1};
        \node at (0,-.35) {\strandlabel{i}};
        \region{0.5,0}{\lambda};
    \end{tikzpicture}
    +
    \begin{tikzpicture}[KM,baseline=-1mm,scale=.9]
        \draw[to-] (0.25,0) arc(0:360:0.25);
        \node at (0,-.4) {\strandlabel{i}};
        \pin{-.25,0}{-1.7,0}{i(u-ix-i)^{-1}};
        \region{.55,0}{\lambda};
    \end{tikzpicture}.
\end{align}
\fi
For $i \neq j$ in $I$, we define the ``quantum'' counterclockwise and clockwise \emph{internal bubbles of color $j$}
\begin{align}\label{int1-q}
    \begin{tikzpicture}[KM,centerzero]
        \draw[-to] (0,-.5)\botlabel{i} to (0,.5);
        \anticlockwiseInternalbubbleR[j]{0,0};
        \region{0.4,0}{\lambda};
    \end{tikzpicture}
    &:= j^{1+h_j(\lambda)}\,
    \begin{tikzpicture}[KM,centerzero]
        \draw[-to] (0,-.5)\botlabel{i} to (0,.5);
        \draw[to-] (0.7,0) arc(360:0:0.2);
        \pinpin{0.3,0}{0,0}{-1.5,0}{(i-j+ix-jy)^{-1}};
        \strand{.5,-.38}{j};
        \region{0.95,0}{\lambda};
    \end{tikzpicture}
    + 
    \left[j^{h_j(\lambda)}
        \begin{tikzpicture}[KM,centerzero]
            \draw[-to] (0,-.5)\botlabel{i} to (0,.5);
            \pin{0,0}{-1.2,0}{(u-ix-i)^{-1}};
            \strand{.6,-.38}{j};
            \stretchedanticlockwisebubble{.6,0};
            \bubblelabel{.6,0}{\frac{u}{j}-1};
            \region{1.2,0}{\lambda};
        \end{tikzpicture} 
    \right]_{u:-1}\!\!\!,
    \\ \label{int2-q}
    \begin{tikzpicture}[KM,centerzero]
        \draw[-to] (0,-.5)\botlabel{i} to (0,.5);
        \clockwiseInternalbubbleL[j]{0,0};
        \region{-.4,0}{\lambda};
    \end{tikzpicture}
    &:= j^{1-h_j(\lambda)}
    \begin{tikzpicture}[KM,centerzero]
        \draw[-to] (0,-.5)\botlabel{i} to (0,.5);
        \draw[to-] (-0.7,0) arc(-180:180:0.2);
        \pinpin{-0.3,0}{0,0}{1.5,0}{(i-j+iy-jx)^{-1}};
        \strand{-.5,-.38}{j};
        \region{-.95,0}{\lambda};
    \end{tikzpicture}
    +
    \left[j^{-h_j(\lambda)}
        \begin{tikzpicture}[KM,centerzero]
            \draw[-to] (-1,-.5)\botlabel{i} to (-1,.5);
            \pin{-1,0}{.2,0}{(u-ix-i)^{-1}};
            \stretchedclockwisebubble{-1.6,0};
            \bubblelabel{-1.6,0}{\frac{u}{j}-1};
            \strand{-1.6,-.38}{j};
            \region{-2.15,0}{\lambda};
        \end{tikzpicture}\  
    \right]_{u:-1}\!\!\!.
\end{align}
The conventions here are similar to \cref{int1,int2}.
For similar reasons to \cref{daunting1,daunting2} in the previous section, the infinite vertical compositions in the next definitions
make sense:
\begin{align}\label{int3-q}
    \begin{tikzpicture}[KM,centerzero]
        \draw[-to] (0,-.5)\botlabel{i} to (0,.6);
        \anticlockwiseInternalbubbleR{0,0};
        \region{0.4,0.05}{\lambda};
    \end{tikzpicture}
    &:= tz^{-1} 
    \begin{tikzpicture}[KM,centerzero]
    \node at (-1,.2) {$\textstyle\prod_{i \neq j \in I}$};
        \draw[-to] (0,-.5)\botlabel{i} to (0,.6);
        \anticlockwiseInternalbubbleR[j]{0,.2};
        \pin{0,-.3}{.8,-.3}{(i+ix)^{-1}};
        \region{1,0.25}{\lambda};
    \end{tikzpicture}
    \ ,&
    \begin{tikzpicture}[KM,centerzero]
        \draw[-to] (0,-.5)\botlabel{i} to (0,.5);
        \clockwiseInternalbubbleL{0,0};
        \region{-.4,0}{\lambda};
    \end{tikzpicture}
    &:= -t^{-1}z^{-1}
    \begin{tikzpicture}[KM,centerzero]
     \node at (-1,.2) {$\textstyle\prod_{i \neq j \in I}$};
        \draw[-to] (0,-.5)\botlabel{i} to (0,.6);
        \clockwiseInternalbubbleL[j]{0,.2};
                \pin{0,-.3}{.8,-.3}{(i+ix)^{-1}};
        \region{-.8,-.3}{\lambda};
        \end{tikzpicture}\ .
\end{align}
Opening some parentheses and using \cref{trickconsequence} gives the following more explicit formulae:
\begin{align}\label{daunting1-q}
    \begin{tikzpicture}[KM,centerzero]
        \draw[-to] (0,-.5)\botlabel{i} to (0,.5);
        \anticlockwiseInternalbubbleR{0,0};
        \region{0.4,0}{\lambda};
    \end{tikzpicture}
    &=
    tz^{-1} \sum_{\substack{J, K\text{ with }|J|<\infty\\J \sqcup K=I-\{i\}}}
    \left[
        \begin{tikzpicture}[KM,centerzero]
            \draw[-to] (0,-.6)\botlabel{i} to (0,.5);
            \node at (-4,.2) {$\prod_{j \in J} j^{1+h_j(\lambda)}$};
            \draw[to-] (.7,.2) arc(360:0:0.2);
            \strand{.5,-.18}{j};
            \pinpin{0.3,.2}{0,.2}{-1.6,.2}{(i-j+ix-jy)^{-1}};
            \node at (1.95,.2) {$\prod_{k \in K} k^{h_k(\lambda)}$};
            \stretchedanticlockwisebubble{3.45,.2};
            \bubblelabel{3.45,.2}{\frac{u}{k}-1};
            \strand{3.45,-.15}{k};
            \pin{0,-.4}{-2,-.4}{(i+ix)^{-1} (u-ix-i)^{-1}};
            \region{1.2,-.35}{\lambda};
        \end{tikzpicture} 
    \right]_{u:-1}\!\!\!,
    \\ \label{daunting2-q}
    \begin{tikzpicture}[KM,centerzero]
        \draw[-to] (0,-.5)\botlabel{i} to (0,.5);
        \clockwiseInternalbubbleL{0,0};
        \region{-0.4,0}{\lambda};
    \end{tikzpicture}\ \ 
    &= -t^{-1}z^{-1}\!\!\!\!
    \sum_{\substack{J, K\text{ with }|J|<\infty\\J \sqcup K=I-\{i\}}}
    \left[
        \begin{tikzpicture}[KM,centerzero]
            \draw[-to] (0,-.6)\botlabel{i} to (0,.5);
            \node at (-1.95,.2) {$\prod_{j \in J}j^{1-h_j(\lambda)}$};
            \draw[to-] (-.7,.2) arc(-180:180:0.2);
            \strand{-.5,-.18}{j};
            \pinpin{-0.3,.2}{0,.2}{1.55,.2}{(i-j+iy-jx)^{-1}};
            \node at (-5.3,.2) {$\prod_{k \in K} k^{-h_k(\lambda)}$};
            \stretchedclockwisebubble{-3.7,.2};
            \bubblelabel{-3.7,.2}{\frac{u}{k}-1};
            \strand{-3.7,-.15}{k};
            \pin{0,-.4}{1.9,-.4}{(i+ix)^{-1}(u-ix-i)^{-1}};
            \region{-1.55,-.35}{\lambda};
        \end{tikzpicture} 
    \right]_{u:-1}\!\!\!.
\end{align}

\begin{theo}\label{average-q}
    In the quantum case, the leftward Heisenberg caps and cups are related to the leftward Kac-Moody caps and cups by 
    \begin{align}\label{leftwardbuggers-q}
        \begin{tikzpicture}[KM,anchorbase,scale=1.5]
            \draw (0.3,-0.4)\botlabel{j} to (0.3,-.2);
            \draw[-to] (-0.1,-.2) to (-0.1,-0.4)\botlabel{i};
            \draw[Q] (.3,-.2) to[out=90,in=90,looseness=1.5] (-.1,-.2);
            \notch{-.1,-.2};
            \notch{.3,-.2};
            \region{.5,-.3}{\lambda};
        \end{tikzpicture}
        &= \delta_{i,j}\ 
        \begin{tikzpicture}[KM,anchorbase,scale=1.5]
            \draw (0.3,-0.4) to (0.3,-.2);
            \draw[-to] (-0.1,-.2) to (-0.1,-0.4)\botlabel{i};
            \draw (.3,-.2) to[out=90,in=90,looseness=2] (-.1,-.2);
            \anticlockwiseInternalbubbleR{.27,-0.18};
            \region{.6,-.2}{\lambda};
        \end{tikzpicture}\ ,&
        \begin{tikzpicture}[KM,anchorbase,scale=1.5]
            \draw[-] (0.3,0.4)\toplabel{j} to (0.3,.2);
            \draw[-to] (-0.1,.2) to (-0.1,0.4)\toplabel{i};
            \draw[Q] (.3,.2) to[out=-90,in=-90,looseness=1.5] (-.1,.2);
            \notch{-.1,.2};
            \notch{.3,.2};
            \region{.5,.3}{\lambda};
        \end{tikzpicture}
        &= \delta_{i,j}\ 
        \begin{tikzpicture}[KM,anchorbase,scale=1.5]
            \draw[-] (0.3,0.4) to (0.3,.2);
            \draw[-to] (-0.1,.2) to (-0.1,0.4)\toplabel{i};
            \draw (.3,.2) to[out=-90,in=-90,looseness=2] (-.1,.2);
            \clockwiseInternalbubbleL{-.07,0.18};
            \region{.5,.2}{\lambda};
        \end{tikzpicture}
    \end{align}
    for $i,j \in I$ and $\lambda \in X$.  Hence, the Heisenberg bubble generating functions are related to the Kac-Moody bubble generating functions by
    \begin{align}\label{newbub-q}
        \begin{tikzpicture}[Q,baseline=-1mm]
            \draw[to-] (-0.25,0) arc(180:-180:0.25);
            \node at (.53,0) {\color{black}$(u)$};
            \region{.95,0}{\lambda};
        \end{tikzpicture}\!
        &= tz^{-1} {\textstyle\prod_{i \in I}} i^{h_i(\lambda)}\ 
        \begin{tikzpicture}[KM,baseline=-1mm]
            \draw[to-] (-0.25,0) arc(180:-180:0.25);
            \node at (.87,0) {\color{black}$(\frac{u}{i}-1)$};
            \region{1.55,0}{\lambda};
            \strand{0,-.4}{i};
        \end{tikzpicture},&
        \begin{tikzpicture}[Q,baseline=-1mm]
            \draw[-to] (-0.25,0) arc(180:-180:0.25);
            \node at (.53,0) {$\color{black}(u)$};
            \region{.95,0}{\lambda};
        \end{tikzpicture}\!
        &= - t^{-1}z^{-1} {\textstyle\prod_{i \in I}} i^{-h_i(\lambda)}\ 
        \begin{tikzpicture}[KM,baseline=-1mm]
            \draw[-to] (-0.25,0) arc(180:-180:0.25);
            \node at (.87,0) {$\color{black}(\frac{u}{i}-1)$};
            \strand{0,-.4}{i};
            \region{1.55,0}{\lambda};
        \end{tikzpicture}.
    \end{align}
    \end{theo}

\begin{proof}[Sketch proof]
The argument is similar to the proof of \cref{average}.
We define 
$\ \begin{tikzpicture}[Q,centerzero,scale=.8]
\draw[-,shadow] (-0.25,-0.2) to [out=90,in=90,looseness=3](0.25,-0.2);
\draw[-to] (-0.25,-0.22) to [out=-90,in=90,looseness=3](-0.25,-0.23);
\end{tikzpicture}\ $
and
$\ \begin{tikzpicture}[Q,centerzero,scale=.8]
\draw[-,shadow] (-0.25,0.2)  to [out=-90,in=-90,looseness=3](0.25,0.2);
\draw[to-] (-0.25,0.22) to [out=-90,in=90,looseness=3](-0.25,0.23);
\end{tikzpicture}\ $ so that
\begin{align*}
\begin{tikzpicture}[KM,anchorbase,scale=1.5]
\draw (0.3,-0.4)\botlabel{j} to (0.3,-.3);
\draw[-to] (-0.1,-.3) to (-0.1,-0.4)\botlabel{i};
\draw[Q,shadow] (.3,-.2) to[out=90,in=90,looseness=2] (-.1,-.2);
%\draw[Q,-to] (-0.1,-0.22) to [out=-90,in=90](-0.1,-0.23);
\draw[Q](-.1,-.2) to (-.1,-.3);\draw[Q](.3,-.2) to (.3,-.3);
\notch{-.1,-.3};
\notch{.3,-.3};
\region{.5,-.3}{\lambda};
\end{tikzpicture}
&=
\delta_{i,j}\ 
\begin{tikzpicture}[KM,anchorbase,scale=1.5]
\draw (0.3,-0.4) to (0.3,-.2);
\draw[-to] (-0.1,-.2) to (-0.1,-0.4)\botlabel{i};
\draw (.3,-.2) to[out=90,in=90,looseness=2] (-.1,-.2);
\anticlockwiseInternalbubbleR{.27,-0.18};
\region{.6,-.2}{\lambda};
\end{tikzpicture}\ ,&
\begin{tikzpicture}[KM,anchorbase,scale=1.5]
\draw (0.3,0.4)\toplabel{j} to (0.3,.3);
\draw[-to] (-0.1,.3) to (-0.1,0.4)\toplabel{i};
\draw[Q,shadow] (.3,.2) to[out=-90,in=-90,looseness=2] (-.1,.2);
%\draw[Q,-to] (-0.1,0.22) to [out=90,in=-90](-0.1,0.23);
\draw[Q](-.1,.2) to (-.1,.3);\draw[Q](.3,.2) to (.3,.3);
\notch{-.1,.3};
\notch{.3,.3};
\region{.5,.3}{\lambda};
\end{tikzpicture}
&= \delta_{i,j}\ 
\begin{tikzpicture}[KM,anchorbase,scale=1.5]
\draw[-] (0.3,0.4) to (0.3,.2);
\draw[-to] (-0.1,.2) to (-0.1,0.4)\toplabel{i};
\draw (.3,.2) to[out=-90,in=-90,looseness=2] (-.1,.2);
\clockwiseInternalbubbleL{-.07,0.18};
\region{.5,.2}{\lambda};
\end{tikzpicture}\ .
\end{align*}
Thus, to prove \cref{leftwardbuggers-q}, we must show that 
$\ \begin{tikzpicture}[Q,centerzero,scale=.8]
\draw[-,shadow] (-0.25,-0.2) to [out=90,in=90,looseness=3](0.25,-0.2);
\draw[-to] (-0.25,-0.22) to [out=-90,in=90,looseness=3](-0.25,-0.23);
\end{tikzpicture}\ 
=
\ \begin{tikzpicture}[Q,centerzero,scale=.8]
\draw[to-] (-0.25,-0.2) to [out=90,in=90,looseness=3](0.25,-0.2);
\end{tikzpicture}\ $
and
$\ \begin{tikzpicture}[Q,centerzero,scale=.8]
\draw[-,shadow] (-0.25,0.2)  to [out=-90,in=-90,looseness=3](0.25,0.2);
\draw[to-] (-0.25,0.22) to [out=-90,in=90,looseness=3](-0.25,0.23);
\end{tikzpicture}\ =
\ \begin{tikzpicture}[Q,centerzero,scale=.8]
\draw[to-] (-0.25,0.2)  to [out=-90,in=-90,looseness=3](0.25,0.2);
\end{tikzpicture}\ $.
To treat the case $\kappa \leq 0$,
we introduce the new counterclockwise 
bubble generating function
$\ \begin{tikzpicture}[Q,baseline=-1mm,scale=.9]
                \draw[-,shadow] (-0.25,0) arc(180:0:0.25);
                \draw[-to] (-.25,-.01) to (-.25,-.02);
                \draw[-] (0.25,0) arc(0:-180:0.25);
                \node at (.56,0) {\color{black}$(u)$};
            \region{1,0}{\lambda};
            \end{tikzpicture}$, which is the formal Laurent series in $Z(\catR)\llbracket u^{-1}\rrbracket$ with leading term $tz^{-1} u^\kappa$ such that
    \(
        \left[
            \begin{tikzpicture}[Q,baseline=-1mm,scale=.9]
                \draw[-,shadow] (-0.25,0) arc(180:0:0.25);
                                \draw[-to] (-.25,-.01) to (-.25,-.02);
                                \draw[-] (0.25,0) arc(0:-180:0.25);
                \node at (.55,0) {\color{black}$(u)$};
            \region{1,0}{\lambda};
            \end{tikzpicture}
        \right]_{u:<0}
        =
        \begin{tikzpicture}[Q,baseline=-1mm]
            \draw[-,shadow] (-0.25,0) arc(180:0:0.25);
                            \draw[-to] (-.25,-.01) to (-.25,-.02);
                        \draw[-] (0.25,0) arc(0:-180:0.25);
            \squared{.18,-.18}{u};
        \end{tikzpicture}
    \).
The shorthand $\beta_{i,r}$ now denotes the coefficients in the expansion
    \begin{equation}\label{how-q}
        i^{h_i(\lambda)}\,
        \begin{tikzpicture}[KM,baseline=-1mm]
            \stretchedanticlockwisebubble{0,0}{u};
            \bubblelabel{0,0}{\frac{u}{i}-1};
            \node at (0,-.35) {\strandlabel{i}};
            \region{0.5,0}{\lambda};
        \end{tikzpicture}
        =
        \sum_{r=0}^{h_i(\lambda)}
        \beta_{i,r} u^r,
    \end{equation}
so that $\beta_{i,h_i(\lambda)} = 1$ when $h_i(\lambda) \geq 0$. 
\iffalse
    Assuming that $h_i(\lambda) \geq 0$ so that it is defined, we have that $\beta_{i,h_i(\lambda)} = 1$.  We then have that
    \begin{equation}\label{abitmore-q}
        {\textstyle\prod_{k \in K}} k^{h_k(\lambda)}
        \begin{tikzpicture}[KM,baseline=-1mm]
            \stretchedanticlockwisebubble{0,0}{u};
            \bubblelabel{0,0}{\frac{u}{k}-1};
            \node at (0,-.35) {\strandlabel{k}};
            \region{0.5,0}{\lambda};
        \end{tikzpicture}
        =
        \sum_{\substack{(r_k)_{k \in K}\\0 \leq r_k \leq h_k(\lambda)}}
        u^{\sum_k r_k}
        {\textstyle\prod_{k \in K}} \beta_{k,r_k}.
    \end{equation}
    \fi
    Now the proof proceeds with analogues of the Claims 1--6 from the proof of \cref{average}.
    The homomomorphism
    $\widehat{\theta}_J:\widehat{A}_J \rightarrow Z(\catR_\lambda)$
    referred to as ``fancy notation'' in the proof of Claim 1 is now defined so that
%     $\widehat{A}_J$ is the completion of $A_J$ at the maximal ideal $\left(\frac{x_j}{j} - 1 \:\Big|\: j \in J\right)$, and 
$f(x_j \mid j \in J)$ goes to the natural transformation obtained by pinning $f(j x_j+j \mid j \in J)$ to \(
        \prod_{j \in J} 
        \begin{tikzpicture}[KM,baseline=-1mm,scale=.8]
            \draw[to-] (-0.25,0) arc(180:-180:0.25);
            \node at (0,-.4) {\strandlabel{j}};
            \opendot{.25,0};
            \region{.6,0}{\lambda};
        \end{tikzpicture}
    \)
    with the variable $x_j$ corresponding to the dot on the $j$th bubble.
\end{proof}

\begin{rem}
The formulae \cref{newbub-q} already appeared in \cite[(5.38)]{HKM}
(with a less direct proof).
\end{rem}